\documentclass[a4paper]{article}
 \usepackage{etoolbox}
 \newbool{PREPRINT}
 \booltrue{PREPRINT}

\ifbool{PREPRINT}{ 
\usepackage[colorlinks,bookmarksopen,bookmarksnumbered,citecolor=red,urlcolor=red,breaklinks=true]{hyperref}
\usepackage[affil-it]{authblk}
\usepackage[cm]{fullpage}
\usepackage{amsthm}

}{
\newcommand{\url}[1]{\texttt{#1}}
} 
\usepackage{amsmath, amssymb,mathtools}
\usepackage{times,microtype}
\usepackage{float}
\usepackage{graphicx,  amsthm, parskip, bbm}
\usepackage{enumerate,cleveref}
\usepackage{caption}
\usepackage{subcaption}
\renewcommand{\vec}[1]{\boldsymbol{#1}}
\usepackage{color}
\usepackage{setspace}
\usepackage[margin=3.0cm,includefoot]{geometry}
\usepackage{subcaption}

\newcommand{\figdir}{./}

\crefname{appendix}{}{}
\crefname{equation}{Eq}{Eqs}
\newtheorem*{theorem*}{Complexity Theorem}
\theoremstyle{plain}

\title{Multilevel Monte Carlo and Improved Timestepping Methods in Atmospheric Dispersion Modelling}
\ifbool{PREPRINT}{ 
\author[1]{Grigoris Katsiolides}
\author[*,1]{Eike H. M\"uller}
\author[1]{Robert Scheichl}
\author[1]{Tony Shardlow}
\author[2]{Michael B. Giles}
\author[3]{David J. Thomson}
\affil[1]{Department of Mathematical Sciences, University of Bath, Bath BA2 7AY, United Kingdom}
\affil[2]{Mathematical Institute, University of Oxford, Andrew Wiles Building, Radcliffe Observatory Quarter, Woodstock Road, Oxford OX2 6GG, United Kingdom}
\affil[3]{Met Office, Fitzroy Road, Exeter EX1 3PB, United Kingdom}
\affil[*]{Email: \texttt{e.mueller@bath.ac.uk}}

}{ 
\author[bath]{G. Katsiolides}
\author[bath]{E.H. M\"{u}ller\corref{cor1}\fnref{fn1}}
\ead{e.mueller@bath.ac.uk}
\author[bath]{R. Scheichl}
\author[bath]{T. Shardlow}
\author[oxford]{M.B. Giles}
\author[metoffice]{D.J. Thomson}
\cortext[cor1]{Corresponding author}
\fntext[fn1]{email: \texttt{e.mueller@bath.ac.uk}}
\address[bath]{Department of Mathematical Sciences, University of Bath, Claverton Down, Bath BA2 7AY, United Kingdom}
\address[oxford]{Mathematical Institute, University of Oxford, Andrew Wiles Building, Woodstock Road, Oxford OX2 6GG, United Kingdom}
\address[metoffice]{Met Office, FitzRoy Road, Exeter EX1 3PB, United Kingdom}
} 

\date{\today}
\begin{document}
\ifbool{PREPRINT}{ 
\maketitle
}{} 
\begin{abstract}
A common way to simulate the transport and spread of pollutants in the atmosphere is via stochastic Lagrangian dispersion models. Mathematically, these models describe turbulent transport processes with stochastic differential equations (SDEs). The computational bottleneck is the Monte Carlo algorithm, which simulates the motion of a large number of model particles in a turbulent velocity field; for each particle, a trajectory is calculated with a numerical timestepping method. Choosing an efficient numerical method is particularly important in operational emergency-response applications, such as tracking radioactive clouds from nuclear accidents or predicting the impact of volcanic ash clouds on international aviation, where accurate and timely predictions are essential. In this paper, we investigate the application of the Multilevel Monte Carlo (MLMC) method to simulate the propagation of particles in a representative one-dimensional dispersion scenario in the atmospheric boundary layer. MLMC can be shown to result in asymptotically superior computational complexity and reduced computational cost when compared to the Standard Monte Carlo (StMC) method, which is currently used in atmospheric dispersion modelling. To reduce the absolute cost of the method also in the non-asymptotic regime, it is equally important to choose the best possible numerical timestepping method on each level. To investigate this, we also compare the standard symplectic Euler method, which is used in many operational models, with two improved timestepping algorithms based on SDE splitting methods.
\end{abstract}
\ifbool{PREPRINT}{ 
\textbf{keywords}:
\newcommand{\sep}{, }
}{ 
\begin{keyword}
} 
Atmospheric dispersion modelling\sep Multilevel Monte Carlo\sep Stochastic differential equations\sep Numerical timestepping methods
\ifbool{PREPRINT}{ 
\\[1ex]
}{ 
\end{keyword}
} 
\maketitle

\section{Introduction}
Modelling the transport and dispersion of atmospheric pollutants is important in applied meteorological research, air-quality modelling for regulatory bodies and emergency-response applications. Numerical models have been used successfully to simulate the spread of radioactive clouds after the Fukushima nuclear accident in 2011 \cite{Draxler2015}, and for evaluating the impact of volcanic ash on aviation following the 2010 and 2011 eruptions of the Eyjafjallaj\"{o}kull and Gr\'{i}msv\"{o}tn volcanoes on Iceland \cite{Dacre2011,Webster2012}. The NAME atmospheric dispersion model \cite{Jones2004,Jones2007}, originally designed to predict the spread of radioactive plumes after the Chernobyl disaster \cite{Smith1989}, is a highly versatile numerical model developed and used by the UK Met Office in a wide range of applications, see e.g. \cite{Dacre2011,Webster2012,Redington2009}. Providing fast-yet-accurate and reliable predictions is crucial, in particular in emergency-response situations, where model outputs have to be passed on to decision makers in a timely fashion. At the heart of NAME, a Stochastic Differential Equation (SDE) describes the evolution of the distributions of position and velocity of the transported particles in the plume. This SDE is solved with a Monte Carlo method, which simulates the propagation of a large number of model particles in a turbulent velocity field. An efficient Monte Carlo algorithm and fast  timestepping methods for integrating the trajectories of the individual particles are crucial for fast and precise forecasts. To assess the performance of a model, it is desirable to minimise the total computational cost for a given level of accuracy (measured by the root-mean-square error) in the prediction of a particular quantity of interest, such as the mean particle position or the concentration field at a given time.

Recently, the Multilevel Monte Carlo (MLMC) method \cite{Heinrich2001,Giles2008} (see \cite{Giles2015} for a comprehensive review) has been used successfully in a wide range of applications from mathematical finance \cite{giles2012multilevel} to uncertainty quantification in subsurface flow \cite{graham2016mixed}. Let $\epsilon$ denote the upper bound (or tolerance) on the root-mean-square error of the numerical method. This error consists of a deterministic discretisation error (bias) and a statistical error due to the limited number of model particles. To reduce both error components, the Standard Monte Carlo (StMC) method usually requires an asymptotic cost $\propto N\cdot M$ since $N$ model particles are propagated over $M$ timesteps. For a first-order method, the number of timesteps (for a given simulation duration) is $M\propto \epsilon^{-1}$ and the central limit theorem implies $N\propto \epsilon^{-2}$, since the statistical error decreases $\propto 1/\sqrt{N}$. In total, this results in an asymptotic cost $\text{Cost}_{\text{StMC}}\le C_{\text{StMC}}\cdot\epsilon^{-3}$. In contrast, the MLMC algorithm can achieve significantly better scaling $\text{Cost}_{\text{MLMC}}\le C_{\text{MLMC}}\cdot\epsilon^{-2}$ by introducing a hierarchy of additional, coarser levels with $M/2,M/4,\dots$ timesteps. On each pair of subsequent levels, paths with two different timestep sizes but the same random noise are used to calculate the same approximation to a particular idealised trajectory. By summing up the contributions on all levels with a telescoping sum, the expected value on the original finest level with $M$ timesteps is recovered.
This reduces the overall runtime for two reasons. Firstly, the cost of one path calculation on a coarser level is proportional to the (smaller) number of timesteps on this level. In addition, only differences between subsequent levels are computed. Those differences have a much smaller variance and can be evaluated with a smaller number of samples. The overall reduction in computational cost is quantified in more detail in the complexity theorem from \cite{Giles2008}, which is reviewed later in this paper. The absolute runtime also depends on the constants $C_{\text{StMC}}$ and $C_{\text{MLMC}}$, which are very sensitive to the numerical timestepping method. This has recently been demonstrated for a set of simpler model problems in \cite{EikeRobTony}. The relative size of those constants is particularly important for relatively loose tolerances, which are typical in operational atmospheric dispersion applications.

Building on \cite{EikeRobTony}, the aim of this work is to investigate the performance enhancements that can be achieved with MLMC and improved timestepping algorithms in atmospheric modelling. This is a challenging problem for several reasons: Firstly, the SDE is nonlinear and the coefficients vary strongly or even diverge in parts of the domain. When simulating vertical dispersion in a neutrally stable atmospheric boundary layer, the ``profiles'' for the turbulent velocity fluctuations $\sigma_U$ and velocity decorrelation time $\tau$ have singularities near the ground and at the top of the domain and need to be regularised appropriately. In addition, suitable boundary conditions have to be applied to ensure that the model particles do not leave the atmospheric boundary layer. Finally, commonly used quantities of interest (in particular the particle concentration) are expressed as the expected value of a function which is not Lipschitz continuous and needs to be approximated correctly. All those issues are addressed in this paper.

While reflective boundary conditions for a one-dimensional diffusion process of the form $dX = a \,dt+\sigma\, dW$ have previously been treated by adding a non-decreasing ``reflection" function (local time of $X_t$ at ground level) to the particle position (see e.g.\ \cite{Skorokhod1961}), this method can not be directly applied to the system of two coupled SDEs for velocity and position that we study here. Instead, we use an approach similar to the \textit{impact law} described for piecewise-smooth dynamical systems in \cite{bernardo2008piecewise}. More specifically, to treat reflective boundary conditions for SDEs, we follow the work in \cite{Wilson1993}, who give an algorithm for elastic reflection at the top and bottom of the domain for a one-dimensional (vertical) dispersion model. The physical justification and implementation of various boundary conditions for atmospheric dispersion models is further discussed in \cite[Chap.~11]{Rodean} and the references given there. Care has to be taken when using reflective boundary conditions in a MLMC context since trajectories of particles on subsequent levels might diverge. To ensure the strong convergence of paths on subsequent levels, we consider an equivalent problem without boundary conditions, which can be obtained by periodically extending the domain. This new approach only requires the inclusion of factors of $\pm 1$ in the coupling between random numbers on the coarse and fine levels.

To predict particle concentrations, we replace the expected value of the indicator function by a suitable smoothing polynomial, matching the first $r$ moments of the indicator function. As discussed in \cite{GilesNagapetyanRitter}, any additional errors introduced by this process can be systematically controlled and we confirm numerically that they are of the same size as the bias and statistical error from the Monte Carlo method.

Several advanced timestepping methods for SDEs have been discussed in the literature (see e.g. \cite{Kloeden2011,Leimkuhler2015}). As a reference method, we use the simple symplectic Euler integrator, which is implemented in many operational atmospheric dispersion models. This integrator is very similar to the well-known Euler--Maruyama method \cite{Maruyama1955,Kloeden2011} and has the same restrictions on the maximal timestep size, which become more severe at the bottom of the boundary layer due to the small velocity decorrelation time. To avoid those stability issues, we also investigate the performance of two improved timestepping methods based on SDE splitting methods. The first one, called the geometric Langevin method \cite{bou2010long}, has previously been applied in the context of MLMC methods for simpler, linear model problems in \cite{EikeRobTony}. The second integrator is the Symmetric Langevin Velocity-Verlet method (BAOAB) \cite{Leimkuhler2015} which has been developed for molecular dynamics simulations. As we will argue below, the fact that it is almost second order leads to superior performance. Other related splitting methods for Langevin dynamics are discussed in \cite{Leimkuhler2015}.

Our numerical results confirm that, for tight tolerances, due to its asymptotic scaling with $\epsilon^{-2}$ the MLMC algorithm always outperforms a standard Monte Carlo method. This difference in performance is particularly pronounced if the particles spend less time close to the ground where the profiles $\tau$ and $\sigma_U$ have large variations. For looser tolerances, the standard Monte Carlo algorithm can be faster than MLMC. However, it is very sensitive to the timestepping method and using the geometric Langevin or BAOAB method can result in speedups of more than one order of magnitude, in particular if particles are released close to the ground.

An alternative way to avoid stability issues due to small velocity decorrelation times at the bottom of the boundary layer is adaptive timestepping.
Two adaptive MLMC algorithms are described and analysed in \cite{Hoel2011,Hoel2014}. The authors of \cite{Hoel2014} either use a sample averaged error to construct a mesh hierarchy, which is used for all realisations, or adapt the timestep size for each path independently, again via  an error estimator. In this work, we focus on stability issues and therefore employ a simple heuristic based on the local velocity decorrelation time. The implementation of the adaptive MLMC algorithm that we use is described in \cite{Giles2016}. We find that adaptive timestepping can further reduce the cost for the symplectic Euler integrator. However, for the geometric Langevin method, our numerical experiments seem to suggest that uniform timestepping methods lead to a shorter overall solution time.

To obtain realistic absolute performance numbers, all results in this paper were obtained with a freely available object-oriented C++ code \cite{MLMCLangevinCode2017}, which can be downloaded under the GPL license.

\paragraph{Structure} 
This paper is organised as follows: in Section \ref{sec:Models}, we introduce the stochastic model for simulating particle transport and dispersion in the stationary inhomogeneous turbulence in the boundary layer and describe the underlying SDE for particle transport and dispersion. Section \ref{sec:Algorithms} introduces the timestepping methods, in particular the geometric Langevin and BAOAB integrators, and contains a review of the MLMC algorithm with particular focus on its computational complexity. The proper treatment of reflective boundary conditions in the context of MLMC is discussed in Section \ref{sec:BoundaryConditions} and concentration estimates with suitable smoothed indicator functions are analysed in Section \ref{sec:SmoothingPolynomials}. In Section \ref{sec:Results}, we present numerical results for three quantities of interest: the mean particle position, the concentration at a particular location and the probability density function of particle position. This includes comparisons between the StMC and the MLMC methods so as to identify the most efficient combination. The impact of the regularisation height on the results and adaptive timestepping methods are discussed in Section \ref{sec:Adaptivity}. We finally conclude and outline ideas for further work in Section \ref{sec:Conclusion}.
Some more technical details are relegated to the appendices, where we analyse our methods based on the theory of modified equations \cite{shardlow2006modified,zygalakis2011existence}, present some mathematical derivations related to the approximation of indicator (box) functions with smoothing polynomials and provide further details on the number of samples per level for two representative MLMC runs.

\section{Lagrangian Models for Atmospheric Dispersion}
\label{sec:Models}
Dispersion models predict the spread and transport of atmospheric pollutants in a given velocity field, i.e.\  the dispersed material is treated as a passive non-interacting tracer. Usually, this velocity field is obtained from the output of a numerical weather prediction (NWP) model. Since the NWP model has a finite grid resolution, the total velocity field $u(x,t)$ can be split into two components. The evolution of the large-scale flow, which is resolved by the grid, is described by the Navier--Stokes equations and simulated as a deterministic velocity field $v(x,t)$ in the dynamical core of the model. In contrast, unresolved turbulence is described by an additional component $U(x,t)$, which 
is not known in detail and whose properties and effects can only be described statistically. In general, the velocity can be decomposed as
\begin{xalignat}{2}
  u(x,t) &= v(x,t) + U(x,t), &\text{where}\qquad \langle u(x,t) \rangle = v(x,t)
\end{xalignat}
and angular brackets $\langle\cdot\rangle$ denote the ensemble average. For an observable quantity $A(x,t)$, the ensemble average is defined by
\begin{equation}
\langle A(x,t) \rangle = \lim_{N \to \infty} \frac{1}{N} \sum_{n = 1}^{N} \alpha^{(n)}(x,t),
\end{equation}
where $\alpha^{(n)}(x,t)$ is the $n^{\text{th}}$ realisation of $A$, i.e.\ particle transport in one particular realisation of the turbulent velocity component $U(x,t)$.

Lagrangian dispersion models describe the motion of atmospheric tracer particles in turbulent atmospheric flow fields by a stochastic differential equation (SDE): each particle follows the mean flow $v(x,t)$ and turbulence is simulated by adding a random velocity component $U_t$. At equilibrium, the Lagrangian velocity $U_t$ of particles which pass through point $x$ at time $t$ agrees in distribution with the turbulent background velocity field $U(x,t)$. For simplicity, we assume that the background flow field vanishes and will set $v(x,t)=0$ in the following.

In \cite{Thomson1987}, several criteria for the selection of models for simulating particle trajectories in turbulent flows are discussed (see also \cite{Rodean} for an overview). For simplicity, we assume here that turbulence in each space direction can be treated independently and only consider the vertical component of the turbulent motion in this article; studying the more general three-dimensional case is the topic of further research. The following one-dimensional model for particle dispersion in a stationary inhomogeneous turbulent flow is derived in \cite{Thomson1987,Rodean}:
\begin{equation}
dU_t = -\frac{U_t}{\tau(X_t)}\,dt + \frac{1}{2}\left[ 1 + \left( \frac{U_t}{\sigma_U(X_t)} \right)^2 \right] \frac{\partial \sigma_U^2}{\partial X}(X_t)\,dt +  \sqrt{\frac{2\sigma_U^2(X_t)}{\tau(X_t)}}\,dW_t,\qquad
dX_t = U_t\,dt.
\label{generalised}
\end{equation}
Here $U_t$ is the Lagrangian vertical velocity fluctuation, $X_t$ is the particle's position and $W_t$ is a Brownian motion. To complete the description in Eq.~(\ref{generalised}), initial conditions for $U_t$ and $X_t$ at time $t=0$ need to be specified.
For stationary turbulence, $U(X,t)$ and $U(X,t_0)$ agree in distribution for any two times $t$ and $t_0$. The velocity variance is $\sigma_U^2(X)=\langle U^2(X,t_0)\rangle$ where $t_0$ is an arbitrary reference time. Using Kolmogorov similarity theory \cite{Kolmogorov1941}, 
the Lagrangian local velocity decorrelation timescale $\tau$ can be expressed in terms of the mean rate of dissipation of turbulent kinetic energy $\epsilon_{\text{turb}}$ as $\tau=2\sigma_U^2/(C_0\epsilon_{\text{turb}})$ \cite{Rodean}. Although for homogeneous turbulence, $\tau$ is identical to the Lagrangian velocity autocorrelation time of $U_t$, this is not true in the more general case of inhomogeneous turbulence considered here.

The term proportional to $\partial\sigma_U^2/\partial X(X_t)$ in Eq.~(\ref{generalised}) arises from the ``well-mixed condition'' \cite{Thomson1987}, which prevents the unphysical spatial agglomeration of particles in an initially homogeneous concentration field. 
The well-mixed condition requires the equilibrium $(X_t,U_t)$-distribution of particles in the model to be the same as the distribution of all air parcels. Mathematically this is equivalent to requiring that the equilibrium distribution (or invariant measure) $p(x,u)$ has the following two properties: the marginal distribution of the spatial component, i.e. $p(x)=\int_{-\infty}^{\infty} p(x,u)\;du$, is uniform in space and the conditional distribution $p(u|x=X)$ is the same as the distribution of $U(X,t)$. The latter is independent of $t$ because of our stationarity assumption and is here assumed to be Gaussian, $p(u|x=X)=\frac{1}{\sqrt{2\pi} \sigma_U(X)}\exp(-u^2/2\sigma^2_U(X))$. 

To simplify notation and to relate to molecular dynamics literature we introduce the following abbreviations
\begin{xalignat*}{3}
\lambda(X) &\equiv \frac{1}{\tau(X)},&
\sigma(X) &\equiv \sqrt{\frac{2\sigma_U^2(X)}{\tau(X)}}, &
V(X,U)  &\equiv - \frac{1}{2} \left[ \sigma_U^2(X) + U^2 \log \left( \sigma_U^2(X)/\sigma_U^2(X_{\text{ref}}) \right) \right]
\end{xalignat*}
(where $X_{\text{ref}}$ is a fixed reference height) so that Eq.~(\ref{generalised}) can be written more compactly as
\begin{equation}
dU_t= - \lambda(X_t) \,U_t\,dt - \frac{\partial V}{\partial X}(X_t,U_t) \, dt + \sigma(X_t) \,dW_t,\qquad dX_t = U_t\;dt.
\label{generalisedsim}
\end{equation}
Although only the derivative of $V(X,U)$ appears in the evolution equation (\ref{generalisedsim}), we make the dependence on the potential explicit since this allows to explore links to Langevin equations which can be written in a similar form (see also \cite{EikeRobTony}).
To complete this model, suitable vertical profiles for $\sigma_U(X)$ and $\tau(X)$ have to be chosen. For a neutrally stable surface layer and horizontally homogeneous turbulence, the authors of \cite{Wilson2009} derive the expression $\tau(X)\approx 0.5 X/\sigma_U(X)$ for the velocity decorrelation time in the case when $\sigma_U(X)$ is constant. If the turbulent eddy size is only constrained by the ground, then a typical eddy size is $\propto X$. Hence, this relationship expresses the fact that on average it takes a time $\propto X/\sigma_U(X)$ for a particle to decorrelate by going around one eddy. As we will see below, the fact that the decorrelation timescale $\tau(X)$ vanishes linearly for $X\rightarrow 0$ has important implications on the performance of the model. For the surface layer, the authors of \cite{Wilson2009} also suggest $\sigma_U=1.3 u_\ast$ where $u_\ast$ is the friction velocity. There is considerable uncertainty in the exact form of the profiles further away from the surface. Rather than exploring a wide range of different parametrisations, we concentrate on one representative case, which allows us to quantify the impact of variations in $\sigma_U(X)$ on the algorithms discussed here. In this paper, we use the form of $\sigma_U(X)$ suggested in \cite{Webster2003} for neutral and stable conditions. In summary, the following profiles are used to describe the propagation of particles in a boundary layer of depth $H$:
\begin{xalignat}{3}
\sigma_U(X) &= \kappa_{\sigma}u_\ast \left( 1 - \frac{X}{H} \right)^{\frac{3}{4}}, & \tau(X) &= \kappa_{\tau} \frac{X}{\sigma_U(X)}, & X \in [0,H],
\label{eqn:profiles}
\end{xalignat}
and the constants are set to $\kappa_{\sigma}=1.3$, $\kappa_{\tau}=0.5$, $u_\ast=0.2\operatorname{m/s}$.
\begin{figure}
  \begin{center}
    \includegraphics[scale = 0.35]{\figdir/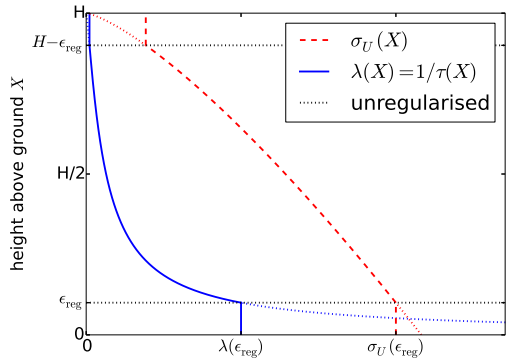}
    \caption{Vertical dependency of the profiles $\sigma_U(X)$ and $\lambda(X)$.}
    \label{fig:profiles}
  \end{center}
\end{figure}
Those profiles are shown in Fig.\ \ref{fig:profiles}; as discussed in Section \ref{sec:regularisation}, $\tau$ (or $\lambda$) needs to be regularised near the bottom boundary. Since $\sigma_U(X)$ and $\tau(X)$ are only defined within the boundary layer, suitable boundary conditions that restrict the particles to the interval $[0,H]$ need to be enforced.

Different boundary conditions for atmospheric dispersion models are discussed in \cite{Rodean}. Here we assume that when a particle hits the boundary, it is reflected elastically (see Fig.\ \ref{fig:reflection_elastic}). Mathematically this can be formulated as in \cite{bernardo2008piecewise}. If at some time $t^{(\text{refl})}$ a particle hits the lower ($X = 0$) or the upper ($X = H$) boundary, we use
\begin{equation}
\lim_{h \to 0, h > 0} U_{t^{(\text{refl})} + h} = - \lim_{h \to 0, h > 0} U_{t^{(\text{refl})} - h}.
\label{reflectLower}
\end{equation}
With initial conditions at $t=0$, it should be noted that, under the assumption that $\tau(X)$ is bounded from below by a non-zero positive constant and that $\sigma_U(X)$ and $\tau(X)$ are sufficiently smooth, it is possible to show existence and uniqueness of a solution for the inhomogeneous model for all $t>0$. The key idea is to find a suitable Lyapunov function and then follow the theory developed in \cite{khasminskii2011stochastic}; including the details of the proof is beyond the scope of this paper.

Given Eqs.\ (\ref{generalisedsim}) and an initial condition, it is now possible to calculate any functional of the path $U_t,X_t$. Any such functional $\phi(U_t,X_t)$ is referred to as a quantity of interest (QoI) in the following. Of particular interest in practical applications is the concentration field at a given time, which can be approximated by the probability of $X_t$ lying in one of the boxes of a predetermined grid. Mathematically for a box $B$, this probability can be expressed as the expected value of the indicator function $\mathbbm{1}_{B}(x)$ which is $1$ for $x\in B$ and $0$ otherwise. Indeed, such piecewise averaged concentration fields are one of the main outputs of widely used dispersion models such as NAME. 

In this paper, we concentrate on calculating the box-averaged concentration $\mathbb{E} \left[ \mathbbm{1}_{[a,b]} (X_T) \right]$ for some $T > 0$, i.e.\ the fraction of particles that land in a specified interval $[a,b] \subset [0,H]$. 
As we will show in Section \ref{sec:pdf}, it is possible to generalise this to calculate the piecewise constant concentration field at the final time by considering a vector of indicator functions. Since it is easier to predict, as a second quantity of interest, we will also consider the mean particle position $\mathbb{E}[X_T]$.
\subsection{Regularisation}
\label{sec:regularisation}
In Eq.~(\ref{generalised}), the profiles $\sigma_U(X)$ and $\tau(X)$ have singularities or vanish at the boundaries $X = 0$ and $X = H$. Consequently, $\partial{V}(X,U)/\partial X$ diverges for $X\rightarrow H$ and $\lambda(X)\rightarrow \infty$ for $X\rightarrow 0$. The latter restricts the timestep in numerical methods to an infinitely small value. It turns out that the issue is not merely the large value of $\lambda$, which could be treated by integrating a larger part of the equation exactly. The function $\lambda(X)$ also has an infinitely large gradient near the ground, which leads to more serious stability issues. In order to deal with these singularities, the functions $\tau(X)$ and $\sigma_U(X)$ have to be regularised. For this, we use the same approach as the Met Office NAME model by imposing a hard cutoff
\begin{xalignat}{2}
\tau(X) &= \begin{cases} 
\tau(\epsilon_{\text{reg}})&\text{if}\quad X < \epsilon_{\text{reg}},\\
\tau(H - \epsilon_{\text{reg}})&\text{if}\quad X > H-\epsilon_{\text{reg}},
\end{cases}\  &
\sigma_U(X) &= \begin{cases} 
\sigma_U(\epsilon_{\text{reg}})&\text{if}\quad X < \epsilon_{\text{reg}},\\
\sigma_U(H - \epsilon_{\text{reg}})&\text{if}\quad X > H-\epsilon_{\text{reg}},
\end{cases}\label{eqn:regularisation}
\end{xalignat}
where $\epsilon_{\text{reg}}\ll H$ is some small regularisation constant. This is shown for $\lambda(X)=1/\tau(X)$ in Fig.\ \ref{fig:profiles}. We will study the impact of $\epsilon_{\text{reg}}$ on the results in Section \ref{sec:sensitivity_regularisation}. Formally the sharp regularisation in Eq.~(\ref{eqn:regularisation}) might cause problems since the derivatives of $\sigma_U(X)$ and $\lambda(X)$ are not continuous at $X=\epsilon_{\text{reg}}$ and $X=H-\epsilon_{\text{reg}}$. As a consequence the drift function in Eq. (\ref{generalised}) has discontinuities at the same positions. This issue could be fixed by using a smooth regularisation, which would ensure that higher derivatives of $\sigma_U(X)$ and $\tau(X)$ vanish at $X=0$. Our numerical confirm that the sharp cutoff does not cause problems in practice.
\section{Algorithms}
\label{sec:Algorithms}
\subsection{Timestepping methods}
\label{sec:TimesteppingMethods}
Timestepping methods for SDEs and especially the Langevin equation are described in \cite{Kloeden2011,EikeRobTony,Leimkuhler2015}. To solve the SDE (\ref{generalisedsim}) numerically, we use three explicit timestepping methods that are of weak order one (i.e. the bias error for fixed travel time grows linearly with the timestep size), but we will argue below that one of the integrators (BAOAB) is nearly second order. The simulated time interval $[0,T]$ is split into $M$ intervals of size $h=T/M$ and, for any $n\in 0,\dots, M$, the quantity $(U_n, X_n)$ denotes the numerical approximation of $(U(t_n), X(t_n))$ at time $t_n = nh$.

All integrators considered in this paper are explicit. The first timestepping method is the Symplectic Euler method given by
\begin{equation}
\begin{aligned}
U_{n + 1} & = \left( 1 - \lambda(X_n) \,h \right) U_n - \frac{\partial V}{\partial X}(X_n,U_n) \,h + \sigma(X_n) \sqrt{h} \,\xi_n,
\\
X_{n + 1} & = X_n + U_{n + 1}\,h.
\label{eqn:MethodSymplecticEuler}
\end{aligned}
\end{equation}
For each timestep $n$, the random number $\xi_n \sim \text{Normal}(0,1)$, independently of all previous timesteps. Due to its simplicity it is very popular in operational models, such as the Met Office's NAME code, and it will serve as the benchmark in the following. It is referred to as ``symplectic'' because, when applied to a Hamiltonian system of ordinary differential equations, it inherits a symplectic structure from the Hamiltonian system.

The method in (\ref{eqn:MethodSymplecticEuler}) is only stable if $ \left| 1 - \lambda(X_n) h \right| < 1$, which leads to the constraint $h < 2/\lambda(X_n)$ on the timestep size. For large values of $\lambda(X)$, which occur close to the lower boundary $X\rightarrow 0$, the timestep size has to be reduced and this increases the computational cost of the method. This problem can be mitigated in adaptive timestepping algorithms, which will be discussed in Section \ref{sec:Adaptivity}.

As an alternative to Eq. (\ref{eqn:MethodSymplecticEuler}) we consider two splitting-based integrators. Those methods rely on the fact that the time evolution in Eq. (\ref{generalisedsim}) can be separated into three parts: a $B$-update, which refers to the velocity increment $dU=- \frac{\partial V}{\partial X}\;dt$; an $A$-update, which corresponds to the position update $dX=U\;dt$, and an Ornstein-Uhlenbeck process $dU=-\lambda\,dt + \sigma \,dW(t)$, which is referred to as the $O$-update. Based on this splitting, the $B$- and $A$- terms are approximated by explicit-Euler steps, and the $O$-update is realised as an exact OU process over a single step. The resulting geometric Langevin method described in \cite{bou2010long} and defined by
\begin{equation}
\begin{aligned}
U_{n + 1}^{\star} & =  e^{- \lambda(X_n) h }\,U_n + \sigma(X_n)
\sqrt{\frac{1 - e^{- 2 \lambda(X_n) h }}{2 \lambda(X_n)}}\,\xi_n,&&\text{($O$-step)}
\\
U_{n + 1} & = U_{n + 1}^{\star} - \frac{\partial V}{\partial X}(X_n,U_{n + 1}^{\star}) \,h,&&\text{($B$-step)}\\
X_{n + 1} & = X_n + U_{n + 1}\,h &&\text{($A$-step)}
\label{eqn:MethodSymplecticEulerOU}
\end{aligned}
\end{equation}
avoids the stability issue. In the notation of~\cite{Leimkuhler2015} this would be referred to as an OBA-integrator. As above, the random variables $\xi_n\sim \operatorname{Normal}(0,1)$ i.i.d. This method was introduced in the context of MLMC in \cite{EikeRobTony}, as a splitting method. In the case that $\sigma_U$ and $\lambda$ are independent of $X$ and $V$ is independent of $U$, it is the first-order splitting method found by taking the Symplectic Euler approximation to a separable Hamiltonian and taking the exact solution to the  Ornstein--Uhlenbeck process ($dU_t=-\lambda \,U_t \,dt + \sigma \,dW_t$). The formulation in Eq.~(\ref{eqn:MethodSymplecticEulerOU}) is generalised to allow application to Eq. (\ref{generalised}).

Following~\cite{Leimkuhler2015} we also consider the corresponding BAOAB method. Here the $B$- and $A$- updates are approximated by two explicit-Euler half-steps, and again the $O$-update is realised as an OU process over a single step. In full, the method is given by
\begin{equation}
\begin{aligned}
U_{n+\frac{1}{2}} & =  U_n - \frac{\partial V}{\partial X} (X_n,U_n) \frac{h}{2}, &&\text{($B$-half-step)}\\
X_{n+\frac{1}{2}} & = X_n + U_{n+\frac{1}{2}} \frac{h}{2},&&\text{($A$-half-step)}\\
U_{n+\frac{1}{2}}^{\star} & = e^{- \lambda (X_{n+\frac{1}{2}}) h} U_{n+\frac{1}{2}} + \sigma(X_{n+\frac{1}{2}}) \sqrt{\frac{1 - e^{- 2 \lambda(X_{n+\frac{1}{2}}) h }}{2 \lambda(X_{n+\frac{1}{2}})}} \xi_n,&&\text{($O$-step)}\\
X_{n+1} & = X_{n+\frac{1}{2}} + U_{n+\frac{1}{2}}^{\star} \frac{h}{2},&&\text{($A$-half-step)}\\
U_{n+1} & = U_{n+\frac{1}{2}}^{\star} - \frac{\partial V}{\partial X} (X_{n+1},U_{n+\frac{1}{2}}^{\star}) \frac{h}{2}&&\text{($B$-half-step)},
\label{eqn:MethodBAOAB}
\end{aligned}
\end{equation}
where $\xi_n\sim \operatorname{Normal}(0,1)$ i.i.d. as above. The method is still first order, but derivation of the second-order modified equation shows that only one new term appears and this term is bounded as $X\downarrow 0$; see \cref{sec:ModifiedEquations}. In other words, the BAOAB-integrator is unaffected by the singularity at $X=0$ and is nearly a second-order method. The alternative ABOBA, which has similar properties, is not discussed here. 

Using this sort of palindromic splitting with second-order integrators for $A$ and $B$,  a weak second-order  method can be derived. This would involve using implicit terms that require a nonlinear solve. We have chosen not to do this, as BAOAB proves itself to be very effective. The complexity of the MLMC method does not benefit from higher-order integrators and first-order accuracy is sufficient. For standard Monte Carlo, higher-order accuracy may be beneficial, though it can be hard to resolve the bias error (i.e., reduce the statistical error sufficiently to see the benefits). The use of second-order splitting methods is discussed further in \cite{EikeRobTony}.
 
\subsection{Standard Monte Carlo simulation}\label{sec:StMCcostanalysis}
In a standard Monte Carlo (StMC) simulation, $N$ independent trajectories $\{(U^{(i)}_n,X^{(i)}_n),i=1,\dots,N\}$ are computed, each associated with a particular i.i.d. realisation of the Gaussian increments $\xi^{(i)}_n$, $n= 0,\ldots,M$, in \eqref{eqn:MethodSymplecticEuler} or \eqref{eqn:MethodSymplecticEulerOU}. 
A quantity of interest (QoI) is evaluated for each trajectory. Possible QoIs are the final particle position or a variable which is one if the particle falls into a defined interval and zero otherwise. The expected value of the QoI can then be estimated by averaging the $N$ realisations of the QoI. In the first example the expectation value is the mean particle position; in the second case it is the fraction of particles which fall into the chosen interval.

We will only consider functionals of the solution $X_T$ at the final time $T$ as QoIs to analyse the cost of the StMC method. In particular, let $\mathcal{P} = \phi (X_T)$. 
Given integer constants $L$ and $M_0$, we approximate $\mathcal{P}$ by $\mathcal{P}_L=\phi(X_{M_L})$ where $X_{M_L}$ is the numerical approximation of the particle position $X_t$ at the final time $T = M_L h_L$, obtained with $M=M_L=M_0 2^L$ timesteps of size $h_L = \frac{T}{M_L} = \frac{T}{M_0}2^{-L}$. The StMC estimator with $N=N_L$ samples for $\mathbb{E}[\mathcal{P}]$ is given by
\begin{equation}
\hat{\mathcal{P}}^{(\text{StMC})}_L = \frac{1}{N_L} \sum_{i=1}^{N_L} \mathcal{P}_L^{(i)},\qquad
\text{where} \quad \mathcal{P}_L^{(i)} = \phi(X_{M_L}^{(i)}).
\label{eqn:StMCestimator}
\end{equation}
Due to the linearity of the expectation, $\mathbb{E}\left[\hat{\mathcal{P}}^{(\text{StMC})}_L\right]=\mathbb{E}\left[\mathcal{P}^{(i)}_L\right]=\mathbb{E}\left[\mathcal{P}_L\right]$. It is also easy to show that the total mean-squared error of the StMC estimator in Eq.~(\ref{eqn:StMCestimator}) is given by
\begin{equation}
\mathbb{E}\left[\left(\hat{\mathcal{P}}^{(\text{StMC})}_L-\mathbb{E}\left[\mathcal{P}\right]\right)^2\right] 
= 
\left(\mathbb{E}\left[\mathcal{P}_L\right]-\mathbb{E}\left[\mathcal{P}\right]\right)^2 + \frac{1}{N_L}\text{Var}\left[\mathcal{P}_L\right].\label{eqn:StMCerror}
\end{equation}
The first term on the right hand side of Eq.~(\ref{eqn:StMCerror}) is the square of the numerical bias and, for the first-order methods described in Eqs.\ (\ref{eqn:MethodSymplecticEuler}) and (\ref{eqn:MethodSymplecticEulerOU}) above, this bias grows linearly with the timestep size: $\mathbb{E}\left[\mathcal{P}_L-\mathcal{P}\right]\propto h_L\propto M_L^{-1}$. The second term is the (mean square) statistical error, which is proportional to $N_L^{-1}$. To limit the total mean-squared error to $\epsilon^2$ (where $\epsilon$ is a fixed tolerance), it is sufficient to choose $M_L$ and $N_L$ such that each of the two terms on the right hand side of Eq.~(\ref{eqn:StMCerror}) is smaller than $\epsilon^2/2$. As already discussed in the introduction, this choice leads to $N_L\propto \epsilon^{-2}$ and $M_L\propto \epsilon^{-1}$, resulting in a total cost of the StMC method 
\begin{equation}
\text{Cost}^{(\text{StMC})} \le C_{\text{StMC}}\cdot\epsilon^{-3}\propto N_L\cdot M_L.
\end{equation}
The size of the constant $C_{\text{StMC}}$ depends on the timestepping method.
As we will argue in the following, the MLMC approach can lead to improved asymptotic scaling with asymptotic cost proportional to $\epsilon^{-2}$. As can be proven with the help of the Central Limit Theorem (see e.g. \cite{Lord2014} for a definition), this is the optimal computational complexity for a Monte Carlo method. 
\subsection{Multilevel Monte Carlo simulation}\label{sec:MLMCcostanalysis}
We generalise the notation from the previous section by considering an arbitrary level $\ell\le L$ with $M_\ell=M_02^\ell$ timesteps of size $h_\ell=\frac{T}{M_\ell}=\frac{T}{M_0}2^{-\ell}$ and $N_\ell$ Monte Carlo samples. On each level, we define the random variable $\mathcal{P}_\ell=\phi(X_{M_\ell})$ where $X_{M_\ell}$ is the approximation of the particle position $X_T$ at time $T=M_\ell h_\ell$, now calculated with a timestep size $h_\ell$. $\mathcal{P}_\ell^{(i)}=\phi(X_{M_\ell}^{(i)})$ is the value of this random variable for the $i$-the sample.

The key idea for the MLMC approximation \cite{Heinrich2001,Giles2008,Giles2015} is that $\mathbb{E} \left[ \mathcal{P}_L \right]$ can be written as the following telescoping sum:
\begin{equation}
\mathbb{E}\left[ \mathcal{P}_L \right] = \sum\limits_{\ell=1}^L \mathbb{E}\left[ \mathcal{P}_\ell - \mathcal{P}_{\ell-1} \right] + \mathbb{E}\left[ \mathcal{P}_0 \right].
\label{eqn:telescoping_sum}
\end{equation}
The MLMC estimator is then defined by
\begin{equation}
\hat{\mathcal{P}}^{(\text{MLMC})}_L = \sum_{\ell=0}^L \hat{Y}_{\ell,N_\ell},
\label{eqn:mlmcestimator}
\end{equation}
where
\begin{equation}
\hat{Y}_{0,N_0} = \hat{\mathcal{P}}^{(\text{StMC})}_0 \quad \text{and} \quad \hat{Y}_{\ell,N_\ell} = \frac{1}{N_\ell} \sum_{i=1}^{N_\ell} Y_\ell^{(i)}\quad\text{with}\quad Y_\ell^{(i)} = \mathcal{P}_\ell^{(i)} - \mathcal{P}_{\ell-1}^{(i)}\, , \quad \text{for} \quad \ell \geq 1.\label{eqn:Yhatell}
\end{equation}
Since the expected value is linear, the MLMC estimator does not introduce any additional bias and
\begin{equation*}
\mathbb{E}\left[ \hat{\mathcal{P}}^{(\text{MLMC})}_L \right] = \sum_{\ell=1}^L \mathbb{E}\left[ \mathcal{P}_\ell - \mathcal{P}_{\ell-1} \right] + \mathbb{E}\left[ \mathcal{P}_0 \right] = \mathbb{E}\left[ \mathcal{P}_L \right] = \mathbb{E}\left[ \hat{\mathcal{P}}^{(\text{StMC})}_L \right].
\end{equation*}
The estimators $\hat{Y}_{\ell,N_{\ell}}$ on two subsequent levels $\ell$ and $\ell+1$ are independent: even though they are evaluated with the same timestep size, the samples $\mathcal{P}_\ell^{(i)}$ required to calculate $\hat{Y}_{\ell,N_\ell}$ and $\hat{Y}_{\ell+1,N_{\ell+1}}$ in Eq. (\ref{eqn:Yhatell}) are generated with different random numbers. Due to this independence of the estimators $\hat{Y}_{\ell,N_{\ell}}$ the mean-square error of $\hat{\mathcal{P}}^{(\text{MLMC})}_L $ is given by
\begin{equation}
\mathbb{E}\left[\left(\hat{\mathcal{P}}^{(\text{MLMC})}_L-\mathbb{E}\left[\mathcal{P}\right]\right)^2\right] 
= 
\left(\mathbb{E}\left[\mathcal{P}_L\right]-\mathbb{E}\left[\mathcal{P}\right]\right)^2 + \sum_{\ell=0}^{L}\frac{1}{N_\ell}\text{Var}\left[Y_\ell\right],\label{eqn:MLMCerror}
\end{equation}
which should be compared to the corresponding expression for StMC in Eq.~(\ref{eqn:StMCerror}). To calculate $\hat{Y}_{\ell,N_\ell}$ on each level $\ell>0$ and to evaluate the MLMC estimator in Eq.~(\ref{eqn:mlmcestimator}), $N_{\ell}$ samples of random numbers $\xi^{(\ell,i)}_n,\;i=1,\dots,N_\ell$, for each timestep $n=1,\dots,M_\ell$ are generated. For the $i$th sample, both a path with $M_{\ell}$ timesteps and a path with $M_{\ell}/2$ timesteps are calculated for the \textit{same} driving Brownian motion $W_t$. The difference $Y_\ell^{(i)}$ between the QoI on the two paths is averaged according to Eq.~(\ref{eqn:Yhatell}). On the coarsest level $\ell=0$, the standard MC estimator with $N_0$ samples and $M_0$ timesteps is used.

The two main advantages of the MLMC algorithm are:
\begin{itemize}
\item Since $\ell<L$ and hence $M_\ell < M_L$, paths on the coarse levels require fewer timesteps and are cheaper to evaluate. In the case we consider here, it can be shown that the computational cost is concentrated on the coarsest level with $M_0=\mathcal{O}(1)$ timesteps, independent of $L$ and of the fine level timestep size.
\item Instead of calculating the quantity of interest itself, only differences $Y_\ell = \mathcal{P}_\ell-\mathcal{P}_{\ell-1}$ of the QoI between subsequent levels are evaluated. Those differences have a smaller variance than $\mathcal{P}_\ell$. As can be seen from Eq.~(\ref{eqn:MLMCerror}), this has a direct impact on the total error.
\end{itemize}
As discussed in detail in \cite{Giles2008}, for a given root-mean-square error $\epsilon$, it is now possible to choose the $N_\ell$ such that the total cost $\propto \sum_{\ell=0}^L N_\ell M_\ell$ is minimised. This is quantified in the following simplified version of the complexity theorem proved in \cite{Giles2008}.\pagebreak
\begin{theorem*}
With the MLMC setting as above and if there exist estimators $\hat{Y}_{\ell,N_\ell}$ and positive constants $c_1, c_2, c_3$ such that
\begin{enumerate}[(i)]
\item $\left|\mathbb{E}[\hat{\mathcal{P}_\ell}-\mathcal{P}]\right| \leq c_1 h_\ell,$
\item $\mathbb{E}[\hat{Y}_{\ell,N_\ell}] = \left \{
\begin{array}{lcc}
\mathbb{E}[\mathcal{P}_0]& &\text{ if } \quad  \ell=0,\\
\mathbb{E}[\mathcal{P}_\ell-\mathcal{P}_{\ell-1}]& &\text{ if } \quad \ell>0,
\end{array}
\right .$
\item $\operatorname{Var}[\hat{Y}_{\ell,N_\ell}] \leq c_2 N_\ell^{-1}h_\ell^2,$
\item $C_\ell$, the computational complexity of $\hat{Y}_{\ell,N_\ell}$, is bounded by $c_3N_\ell h_\ell^{-1}$,
\end{enumerate}
then there exists a positive constant $C_{\text{MLMC}}$ such that for all $\epsilon < e^{-1}$ there are values $L$ and $N_l$ for which $\hat{\mathcal{P}}^{(\operatorname{MLMC})}_L$ satisfies
\begin{center}
$\mathbb{E}\left[ \left( \hat{\mathcal{P}}^{(\operatorname{MLMC})}_L - E\left[ \mathcal{P} \right] \right)^2 \right] < \epsilon^2,$
\end{center}
with a computational complexity $\operatorname{Cost}^{(\operatorname{MLMC})}$ with bound
\begin{center}
$\operatorname{Cost}^{(\operatorname{MLMC})} \leq C_{\text{MLMC}}\cdot \epsilon^{-2}.$
\end{center}
\end{theorem*}
Condition (ii) is obviously true with our definition of $\hat{Y}_{\ell,N_\ell}$. Conditions (i) and (iv) are easy to satisfy and they hold for the two timestepping methods considered in the article. However, the quadratic variance decay in condition (iii) requires the correct coupling of paths between subsequent levels. The random variables $\xi_n$ in Eqs.\ (\ref{eqn:MethodSymplecticEuler}) and (\ref{eqn:MethodSymplecticEulerOU}) that are used for calculating the paths on two coupled levels, have to be chosen correctly. For the symplectic Euler method in Eq.~(\ref{eqn:MethodSymplecticEuler}), this can be achieved by using a simple Brownian bridge, i.e.\ combining the random variables $\xi_n$ and $\xi_{n+1}$ on one level to obtain a random variable $\xi_n^{(\text{coarse})}$ on the next coarser level as
\begin{equation}
\xi_n^{(\text{coarse})} = \frac{1}{\sqrt{2}}(\xi_n + \xi_{n + 1}).
\label{EulerCoupled}
\end{equation}
As shown in \cite{EikeRobTony}, for the geometric Langevin method in Eq.~(\ref{eqn:MethodSymplecticEulerOU}), this is achieved by using 
\begin{equation}
\xi_n^{(\text{coarse})} = \frac{e^{- \lambda(X_n) h }\xi_n + \xi_{n + 1}}{\sqrt{e^{- 2\lambda(X_n) h} + 1}}.
\label{SymEulerCoupled}
\end{equation}
Since the Ornstein-Uhlenbeck step for the BAOAB method is the same as the one for the geometric Langevin method, we can also use equation Eq.~(\ref{SymEulerCoupled}) to couple the random variables of the MLMC coarse step for the BAOAB integrator.

As is common in the MLMC literature, in the numerical section below, we carefully check condition (iii) by plotting the variance $\text{Var}[Y_\ell] = N_\ell \cdot \text{Var}[\hat{Y}_{\ell,N_\ell}]$ as a function of the step size $h_\ell$.
\section{MLMC with Reflective Boundary Conditions}
\label{sec:BoundaryConditions}
\begin{figure}
\centering
\begin{subfigure}{0.49\linewidth}
\centering
\includegraphics[width=0.8\linewidth]{\figdir/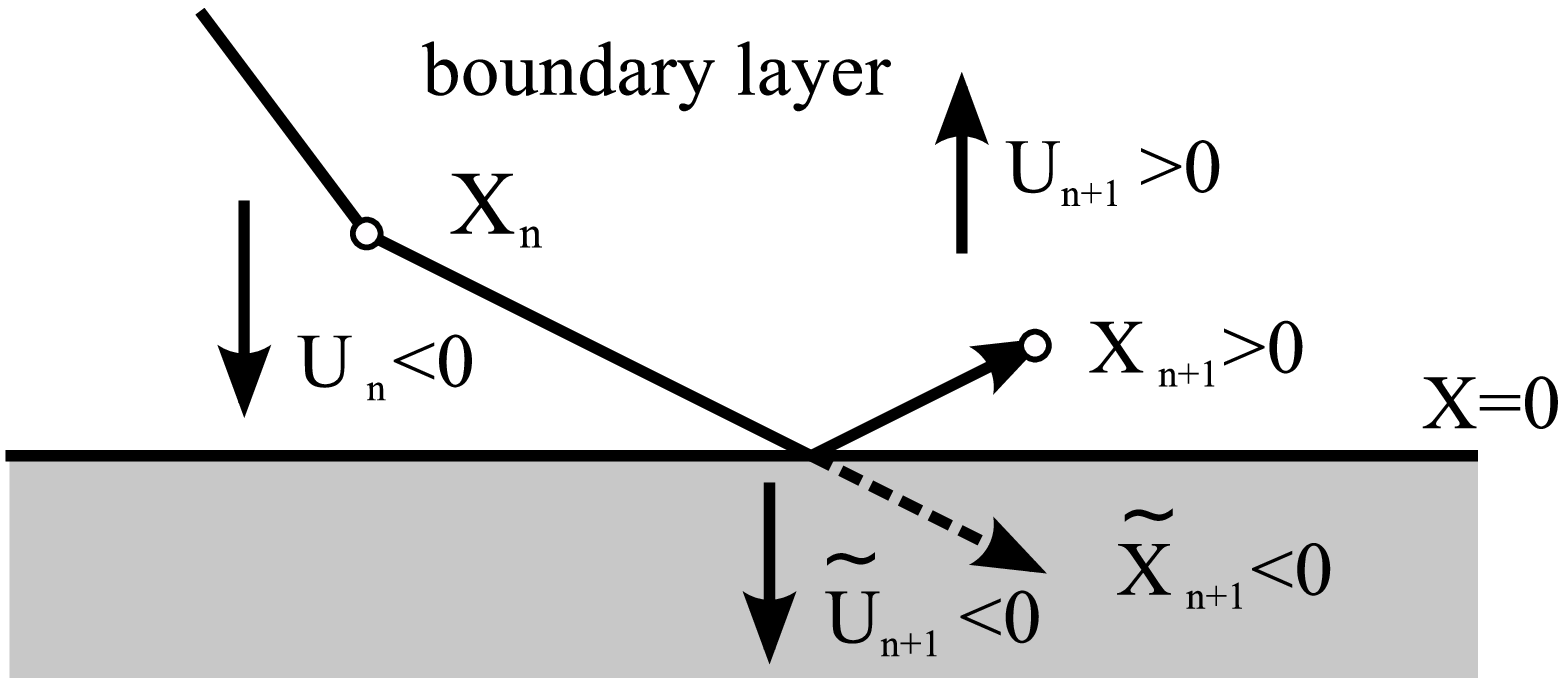}
\caption{Elastic reflection at the lower boundary}
\label{fig:reflection_elastic}
\end{subfigure}
\hfill
\begin{subfigure}{0.49\linewidth}
\centering
\includegraphics[scale = 0.4]{\figdir/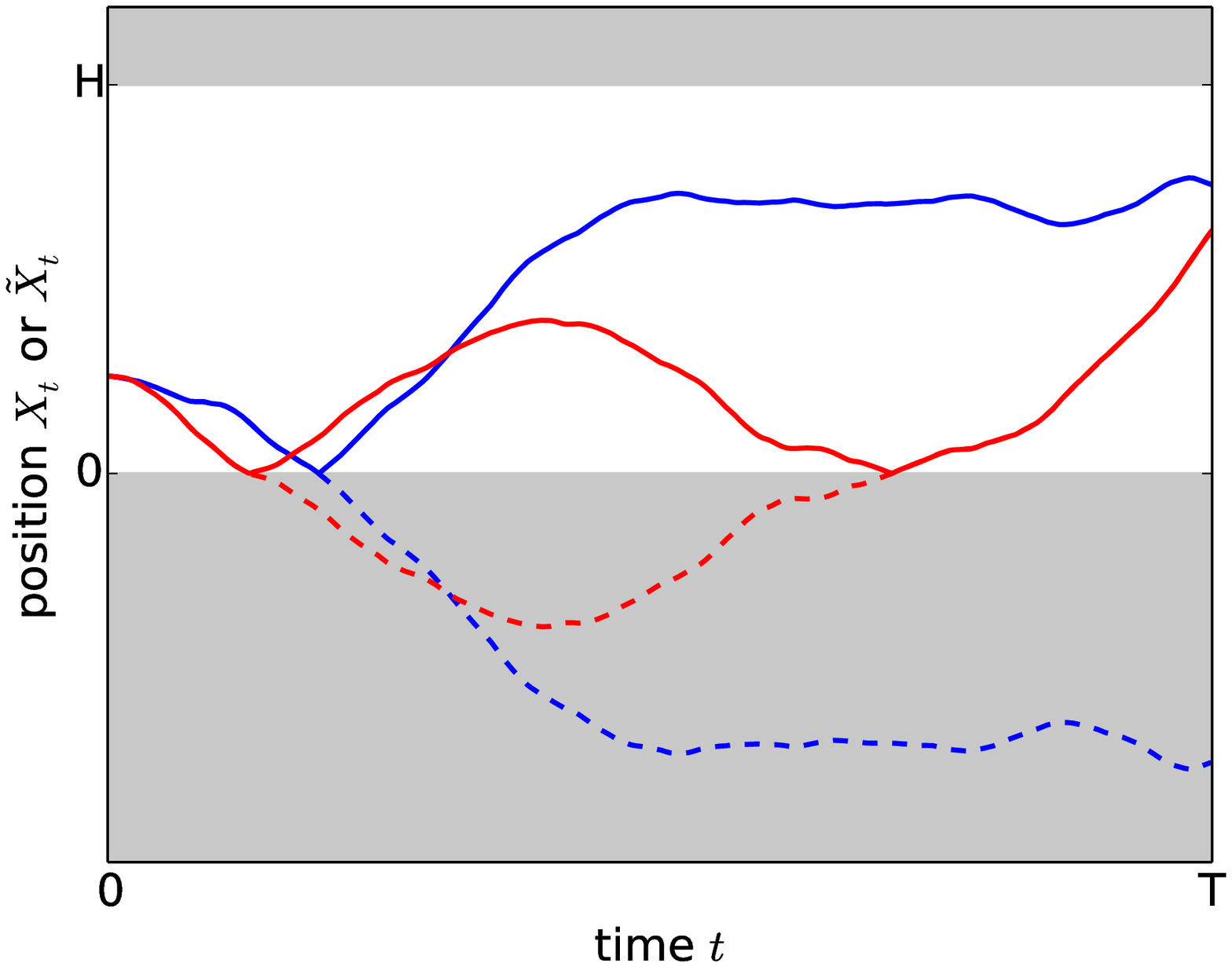}
\caption{Two paths in the original random variable $X_t$ (solid lines) and the extended random variable $\tilde{X}_t$ (dashed lines)}
\label{fig:reflection_schematic}
\end{subfigure}
\caption{Treatment of elastic boundary conditions}
\end{figure}

As described in Section \ref{sec:regularisation}, we assume that particles are reflected elastically at the top and bottom of the domain. As shown in \cite{Wilson1993}, this can be accounted for in the timestepping method by simply checking if a particle has left the domain after each timestep, and, if necessary, replacing $X_{n+1}\mapsto -X_{n+1}$ (or $X_{n+1}\mapsto 2H-X_{n+1}$ at the top boundary) and $U_{n+1}\mapsto-U_{n+1}$ in accordance with Eq.~(\ref{reflectLower}). This is shown schematically in Fig.\ \ref{fig:reflection_elastic}.

However, since the paths on subsequent levels are coupled in MLMC, care has to be taken to prevent the fine and coarse paths from diverging, which would result in violation of condition (iii) in the complexity theorem. This can be achieved by modifying Eqs.\ (\ref{EulerCoupled}) and (\ref{SymEulerCoupled}) appropriately. The key idea is illustrated by the following thought experiment, which we first explain for one reflective boundary (in Section \ref{sec:top_reflection} we show that the same treatment works for both boundaries by considering periodic extensions): imagine extending the domain of $X$ to the whole real line and making it equivalent to the problem without boundary conditions. For this, an extended SDE in the variables $\tilde{X}$ and $\tilde{U}$ is formulated, which is related to the initial one by
\begin{equation}
\begin{pmatrix}
U \\
X
\end{pmatrix}
 = S(\tilde{X})
\begin{pmatrix}
\tilde{U} \\
\tilde{X}
\end{pmatrix},
\label{eqn:reflection}
\end{equation}
where $S(\tilde{X}_t) = (-1)^{n_{\text{refl}}(t)} = \operatorname{sign}(\tilde{X}_t)$ and $n_{\text{refl}}(t)$ is the number of reflections up to time $t$. Let $a(X,U)$ and $b(X)$ be some functions which are defined for positions $X\ge 0$ and any velocity $U$. For the  SDE 
\begin{xalignat}{2}
dU_t &= a(X_t,U_t)\,dt + b(X_t) \,dW_t, &
dX_t &= U_t \,dt,
\label{eqn:problem_nonextended}
\end{xalignat}
with reflection at $X = 0$, the equivalent extended SDE is
\begin{xalignat}{2}
d\tilde{U}_t &= A(\tilde{X}_t,\tilde{U}_t)\,dt + B(\tilde{X}_t) \,dW_t,&
d\tilde{X}_t &= \tilde{U}_t \,dt.
\label{eqn:problem_extended}
\end{xalignat}
The functions $A$ and $B$ in Eq.~(\ref{eqn:problem_extended}) are related to $a$ and $b$ in Eq.~(\ref{eqn:problem_nonextended}) by
\begin{xalignat}{2}
A(\tilde{X}_t,\tilde{U}_t) &=
\begin{cases}
a(\tilde{X}_t,\tilde{U}_t) & \text{if}\quad \tilde{X}_t \geq 0,\\
- a(- \tilde{X}_t,- \tilde{U}_t) & \text{if}\quad \tilde{X}_t < 0,
\end{cases}\  &
B(\tilde{X}_t) &=
\begin{cases}
b(\tilde{X}_t) & \text{if}\quad \tilde{X}_t \geq 0,\\
b(- \tilde{X}_t) & \text{if }\quad \tilde{X}_t < 0.
\end{cases}
\label{eqn:drift_diffusion}
\end{xalignat}
In Fig.\ \ref{fig:reflection_schematic}, both the solution $\tilde{X}_t$ of the extended SDE (\ref{eqn:problem_extended}) and the solution $X_t$ of the original SDE (\ref{eqn:problem_nonextended}) are plotted. The solutions are related by Eq.~(\ref{eqn:reflection}).

One possible approach would now be to use the extended SDE  (\ref{eqn:problem_extended}), which does not have any boundary conditions, and then reconstruct the variables $U_t$, $X_t$ by Eq.~(\ref{eqn:reflection}). However, it turns out that this is not necessary, and we can instead use the original set of equations if we carefully adapt the way the random variables are coupled in the MLMC algorithm.

\subsection{Coupling between paths on different levels}
\label{sec:reflection_analysis}
To achieve this, we derive some relations between the extended and the original variables and the SDEs in Eqs.\ (\ref{eqn:problem_nonextended}) and (\ref{eqn:problem_extended}). Discretising Eq.~(\ref{eqn:reflection}) gives
\begin{equation}
\begin{pmatrix}
U_n \\
X_n
\end{pmatrix}
 = S_n
\begin{pmatrix}
\tilde{U}_n \\
\tilde{X}_n
\end{pmatrix},\label{eqn:XnUn}
\end{equation}
where $S_n = (-1)^{n_{\text{refl}}(t_n)}$. By using this relation in Eq.~(\ref{eqn:drift_diffusion}), we obtain
\begin{xalignat}{2}
A(\tilde{X}_n,\tilde{U}_n) &= S_n a(S_n \tilde{X}_n, S_n \tilde{U}_n) = S_n a(X_n,U_n),&
B(\tilde{X}_n) &= b (S_n \tilde{X}_n) = b (X_n).
\end{xalignat}
Now consider the fine path on a MLMC level with timestep size $h$. Discretising the extended model in Eq.~(\ref{eqn:problem_extended}) with the symplectic Euler method gives
\begin{xalignat}{2}
\tilde{U}_{n + 1} &= \tilde{U}_n + A(\tilde{X}_n,\tilde{U}_n) \,h + B(\tilde{X}_n) \sqrt{h}\,\tilde{\xi}_n,&
\tilde{X}_{n + 1} &= \tilde{X}_n + \tilde{U}_{n + 1} \,h,
\end{xalignat}
where $\tilde{\xi}_n\sim \text{Normal}(0,1)$ i.i.d. Replacing $(\tilde{U}_n,\tilde{X}_n)$ by $(U_n,X_n)$ according to Eq.~(\ref{eqn:XnUn}) results in an equation for the original variables $U_n$ and $X_n$
\begin{equation}
\begin{aligned}
S_{n + 1}^{(f)} U_{n + 1} &= S_n^{(f)} U_n + S_n^{(f)} 
a(X_n,U_n) \,h + b(X_n) \,\sqrt{h}\,\tilde{\xi}_n,\\
S_{n + 1}^{(f)} X_{n + 1} &= S_n^{(f)} X_n + S_{n + 1}^{(f)} \,U_{n + 1} \,h.\label{speedExtended}
\end{aligned}
\end{equation}
To distinguish fine and coarse levels, we use the superscripts ``$(f)$'' and ``$(c)$'' on the reflection factors. Multiplying both sides by $\left(S_{n+1}^{(f)}\right)^{-1}=S_{n+1}^{(f)}$ and defining $s_{n + 1}^{(f)} = S_{n + 1}^{(f)} S_n^{(f)}$, the timestepping scheme can be expressed as
\begin{equation}
\begin{aligned}
U_{n + 1} &= s_{n + 1}^{(f)} \left( U_n + a(X_n,U_n) \,h + b(X_n) \, \sqrt{h} \,S_n^{(f)} \,\tilde{\xi}_n\right),\\
X_{n + 1} &= s_{n + 1}^{(f)} \left( X_n + s_{n + 1}^{(f)} \,U_{n + 1}\, h \right). \label{eqn:reflected_equation}
\end{aligned}
\end{equation}
The same equation can be written down on the next coarser level, where we replace $h\mapsto 2h$, $S_n\mapsto S_n^{(c)}$ etc. 
Note that $s_{n+1}^{(f)}$ equals $-1$ if there is a reflection at time $t_{n + 1}$ or $1$ otherwise. Hence the factors $s_{n+1}^{(f)}$ in Eq.~(\ref{eqn:reflected_equation}) take care of the elastic reflection described above: if $X_{n+1}<0$, velocity and position $(U_{n+1},X_{n+1})$ get replaced by $(-U_{n+1},-X_{n+1})$. In addition, the random forcing term proportional to $b(X_n)$ contains an additional factor $S_n^{(f)}$. This could be taken into account by sampling $\tilde{\xi}_n\sim\text{Normal}(0,1)$ and multiplying each sample by the correct factor $S_n^{(f)}$ as in Eq.~(\ref{eqn:reflected_equation}). However, since the normal distribution is symmetric and $S_n^{(f)}=\pm 1$, on each individual level, the same result is obtained by directly using a random variable $\xi_n\sim \text{Normal}(0,1)$ and replacing $S_n^{(f)}\tilde{\xi}_n\mapsto \xi_n$ to obtain
\begin{equation}
U_{n + 1} = s_{n + 1}^{(f)} \left( U_n + a(X_n,U_n) \,h + b(X_n)  \,\sqrt{h}\,\xi_n\right).
\end{equation}
Formally the random variables $\xi_n$ and $\tilde{\xi}_n$ are related by
\begin{equation}
  \xi_n=S_n^{(f)}\tilde{\xi}_n
\label{eqn:coupled_samples}
\end{equation}
and agree in distribution since $S_n^{(f)}=\pm 1$. To couple the fine and coarse path correctly, it is necessary to require that the random variable $\tilde{\xi}_n^{(\text{coarse})}$ on the \textit{coarse} paths fulfils the conditions in Eqs.\ (\ref{EulerCoupled}) and (\ref{SymEulerCoupled}). For the Euler scheme, this means that
\begin{equation*}
  \tilde{\xi}_n^{(\text{coarse})} = \frac{1}{\sqrt{2}}\left(\tilde{\xi}_n+\tilde{\xi}_{n+1}\right).
\end{equation*}
Using Eq.~(\ref{eqn:coupled_samples}) and the corresponding relation $\xi_n^{(\text{coarse})}=S_n^{(c)}\tilde{\xi}_n^{(\text{coarse})}$, this implies that reflection can be taken into account by replacing Eq.~(\ref{EulerCoupled}) by
\begin{equation}
\xi_n^{(\text{coarse})} = \frac{S_n^{(c)} }{\sqrt{2}} \left(S_n^{(f)} \xi_n + S_{n + 1}^{(f)}  \xi_{n + 1} \right).
\label{eqn:euler_coupled_treatment}
\end{equation}
A similar argument shows that for the geometric Langevin and BAOAB methods Eq.~(\ref{SymEulerCoupled}) needs to be replaced by
\begin{equation}
\xi_n^{(\text{coarse})} = S_n^{(c)} \frac{e^{- \lambda(X_n) h } S_n^{(f)} \xi_n + S_{n + 1}^{(f)} \xi_{n + 1}}{\sqrt{e^{- 2\lambda(X_n) h} + 1}}.
\label{eqn:splitting_coupled_treatment}
\end{equation}
Since this treatment of the reflective boundary conditions does not change the distribution of the normal random variables used in the coarse step, no extra bias is added and the telescoping sum (\ref{eqn:telescoping_sum}) is preserved.

\subsubsection{Reflection at the top boundary}\label{sec:top_reflection}
To treat elastic reflection at both the top and bottom boundary, periodic extensions of the domain can be used. For $X \in \mathbb{R}$, let $n(X) = 2n$ if $X \in \left[(2n - 1)H,(2n + 1)H \right)$, $n \in \mathbb{Z}$ and $\eta(X) = X - n(X)H$ be a function that maps $X$ to $[- H, H)$. Similarly to Eq.~(\ref{eqn:drift_diffusion}), the extended coefficients can be defined as
\begin{xalignat}{2}
A(\tilde{X}_t,\tilde{U}_t) &=
\begin{cases}
a(\eta(\tilde{X}_t),\tilde{U}_t) & \text{if}\quad \eta(\tilde{X}_t) \geq 0,\\
- a(- \eta(\tilde{X}_t),- \tilde{U}_t) & \text{if}\quad \eta(\tilde{X}_t) < 0,
\end{cases}\  &
B(\tilde{X}_t) &=
\begin{cases}
b(\eta(\tilde{X}_t)) & \text{if}\quad \eta(\tilde{X}_t) \geq 0,\\
b(- \eta(\tilde{X}_t)) & \text{if}\quad \eta(\tilde{X}_t) < 0.
\end{cases}
\end{xalignat}
It is easy to verify that the function $\eta(X)$ reduces the problem to the case of reflecting only at $X = 0$ and the same analysis as in Section \ref{sec:reflection_analysis} applies. Therefore, no further changes are required in Eqs.\ (\ref{eqn:euler_coupled_treatment}) and (\ref{eqn:splitting_coupled_treatment}) and $S_n^{(\cdot)}$ counts reflections at both boundaries.

\begin{figure}
\centering
\captionsetup{justification=centering}
\includegraphics[scale = 0.4]{\figdir/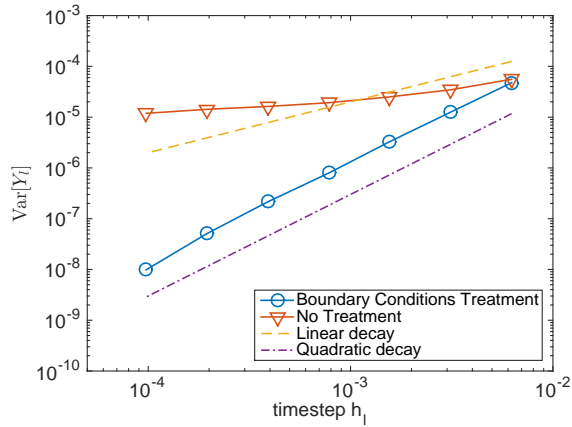}
\caption{Geometric Langevin method - correct boundary condition treatment for
mean particle position}
\label{SymplEulerOUBoundSpeTreat}
\end{figure}

Fig.\ \ref{SymplEulerOUBoundSpeTreat} shows the variance $\operatorname{Var}[\hat{Y}_{\ell,N_\ell}]$ of the difference process as a function of the timestep size. As can be seen from this plot, the quadratic decay required for condition (iii) in the complexity theorem is recovered. In the same figure, we also show the variance decay with naive (and incorrect) coupling of the paths at the reflective boundaries, i.e.\ when setting $S_n^{(\cdot)}=1$. Due to the symmetry properties of the Gaussian increments, this still leads to the same unbiased estimators on each level, satisfying Assumptions (i), (ii) and (iv), but the variance of the difference decays with a rate that is not even linear and violates condition (iii).

Although the function $A(X,U)$ is discontinuous at the regularisation height due to the hard cutoff introduced in Section~\ref{sec:regularisation}, Fig.~\ref{SymplEulerOUBoundSpeTreat} confirms numerically that the method works for the model considered here even with the sharp regularisation in Eq.~(\ref{eqn:regularisation}). While this is not the case for our model, in general the treatment of boundary conditions described here could introduce discontinuities in the functions $A(X,U)$ and $B(X)$ at $X=0$ and $X=H$ which would require careful consideration.
\section{Smoothed indicator functions for concentration estimates}
\label{sec:SmoothingPolynomials}
\label{conc_sub_sec}
The concentration of particles which land in an interval $[a,b]\subset [0,H]$ can be expressed as the expected value of the indicator function defined in Section~\ref{sec:Models}, 
i.e.\ $\mathbb{E} \left[ \mathbbm{1}_{[a,b]} (X_T) \right]$. However, since the indicator function is not continuous, this leads to problems in the MLMC method: it manifests itself in the fact that the variance on subsequent levels only decays linearly with the step size and condition (iii) in the complexity theorem is violated (this is consistent with the fact that the modified-equation analysis in \cite{EikeRobTony}, which is used to prove the complexity theorem there, also assumes that the quantity of interest is Lipschitz continuous). The numerical results below (Figs.\ \ref{varianceratessymplecticeuler01} and \ref{varianceratessymplecticeulerOU01}) confirm that this is indeed the case. 

As we will demonstrate now, this issue can be fixed by approximating the indicator function by a suitable smooth function; any errors introduced by this approximation can be quantified and minimised systematically.
\subsection{Smoothing polynomials}
In order to recover the quadratic variance decay rates and satisfy condition (iii), we construct a smooth function that approximates the indicator function as described in \cite{GilesNagapetyanRitter}. By linearity, it is possible to write
\begin{equation}
\mathbb{E} \left[ \mathbbm{1}_{[a,b]} (X_T) \right] = \mathbb{E} \left[ \theta(X_T-b) \right] - \mathbb{E} \left[ \theta(X_T-a) \right] \qquad \text{with}\quad \theta(X) = \begin{cases}1 & \text{for}\quad X\le0,\\0&\text{otherwise}.\end{cases}
\label{conc}
\end{equation} 
In \cite{GilesNagapetyanRitter}, to approximate the step function $\theta(X)$, the authors define 
\begin{equation}
g_r(X) =
\begin{cases}
1 & \text{if}\quad X < - 1,\\
p_r(X) & \text{if}\quad - 1 \leq X \leq 1, \\
0 & \text{if}\quad X > 1.
\end{cases}
\label{eqn:gr}
\end{equation}
Here, $p_r(X)$ is a polynomial of degree at most $r + 1$ such that $p_r(- 1) = 1$, $p_r(1) = 0$ and the first $r$ moments of $p_r(X)$ over the interval $[-1,+1]$ agree with the moments of the step function $\theta(X)$ on the same interval. This is guaranteed if 
$$
\int_{- 1}^{1} s^j p_r(s)\, ds = \frac{(- 1)^j}{j + 1}\, , \qquad \text{for all} \ \ j = 0, \ldots, r - 1.
$$ 
For an interval of interest $[S_0,S_1]$,  it can be shown \cite{GilesNagapetyanRitter} that, for some $c>0$ and all $\delta$ sufficiently small, 
\begin{equation}
\sup_{s \in [S_0,S_1]} \left| \mathbb{E} \left[ \theta(X_T-s) \right] - \mathbb{E} \left[ g_r \left( \frac{X_T - s}{\delta} \right) \right] \right| \leq c\cdot \delta^{r + 1}, 
\end{equation}
assuming the density of $X_T$ is $r$-times continuously differentiable on an open set containing $[S_0,S_1]$. By choosing $\delta$ and $r$ appropriately, it is possible to approximate any QoI based on the indicator function to arbitrary accuracy.

\begin{figure}
\begin{subfigure}{.5\textwidth}
\centering
\includegraphics[scale = 0.4]{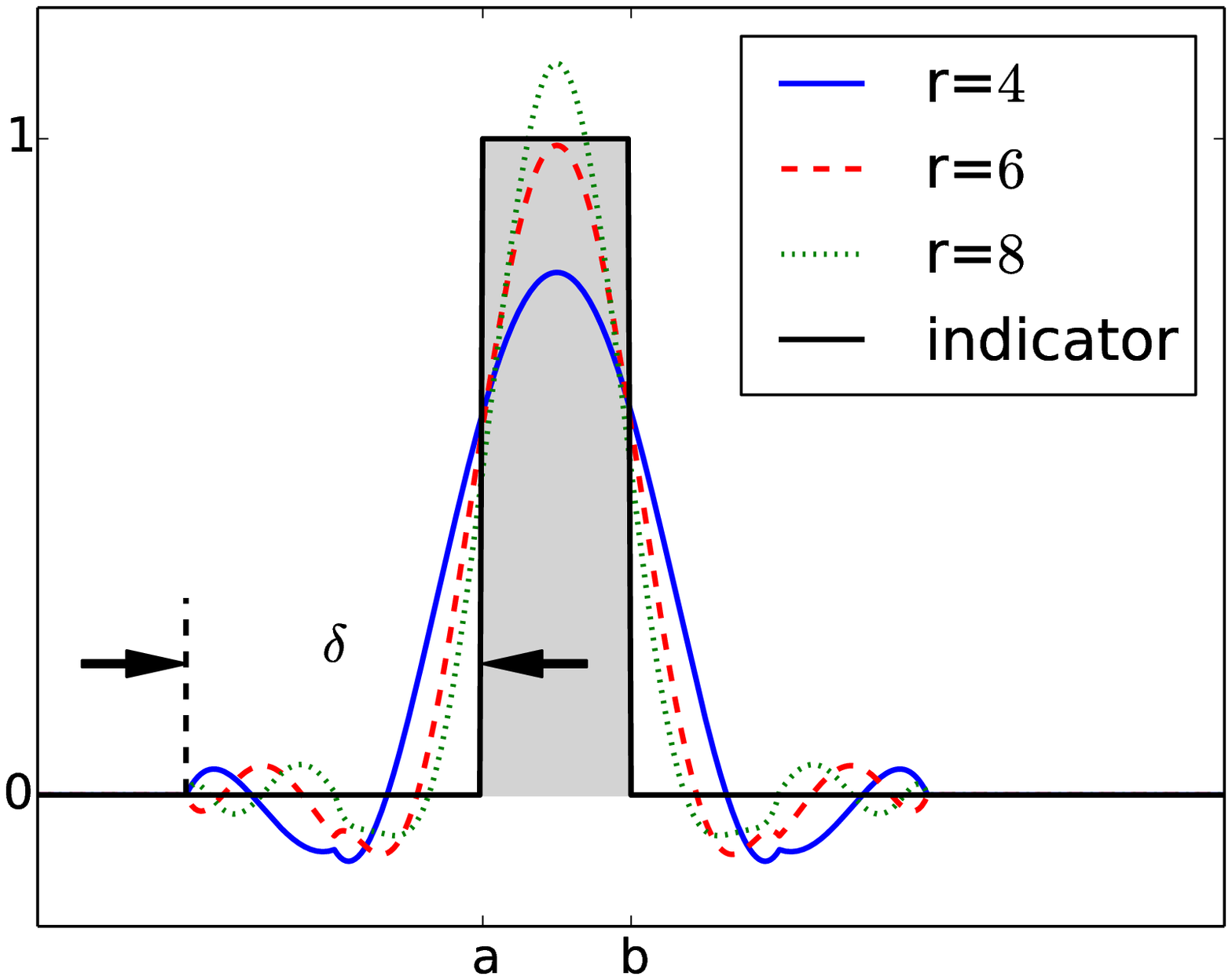}
\caption{$\delta = 0.1$}
\label{Fig:0.1}
\end{subfigure}%
\begin{subfigure}{.5\textwidth}
\centering
\includegraphics[scale = 0.4]{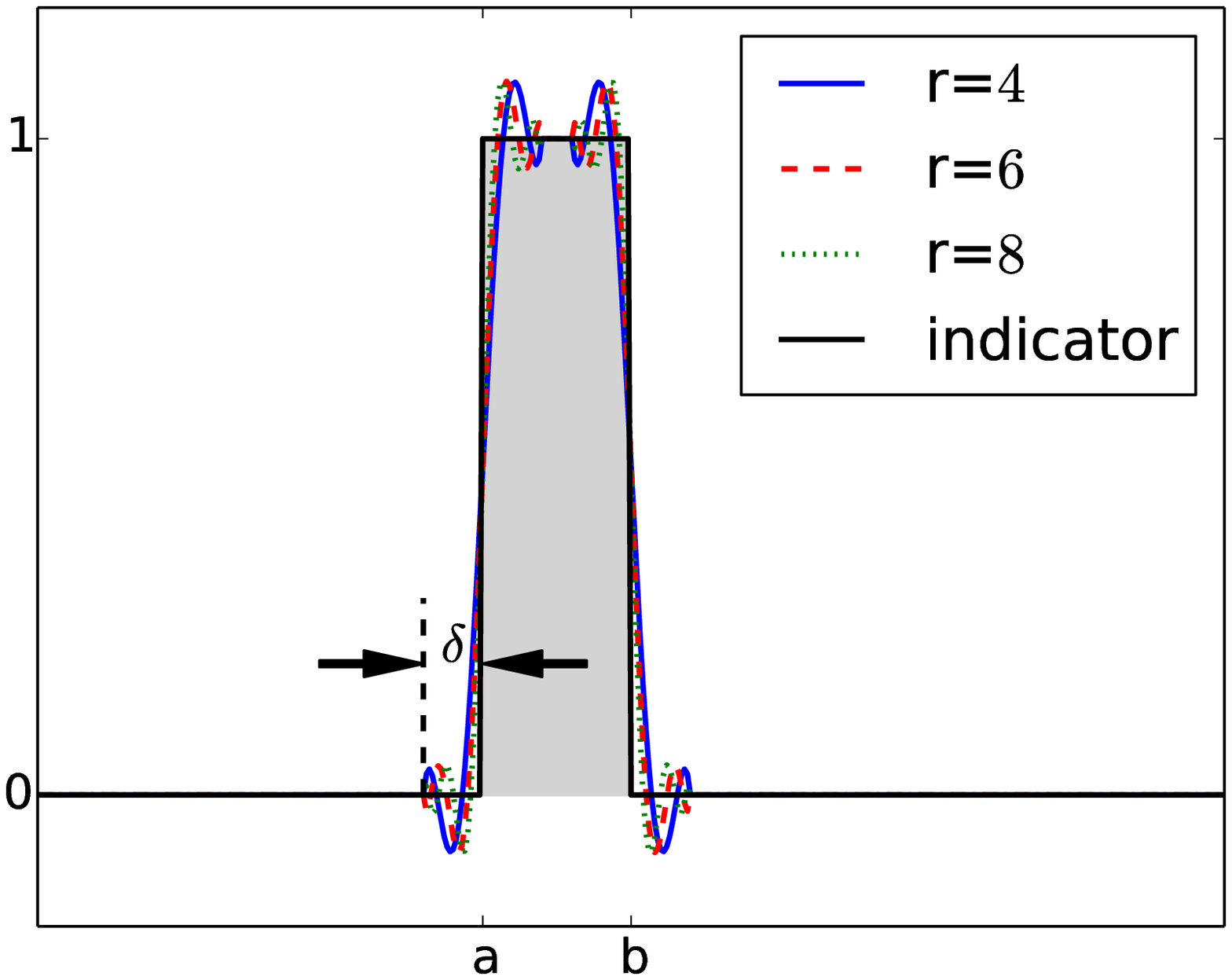}
\caption{$\delta = 0.02$}
\label{0.02}
\end{subfigure}%
\caption{Smoothing polynomials $P^{r,\delta}_{[a,b]}(X)$ for different values of $r$, $\delta$ and for $a=0.1$, $b=0.15$. The indicator function is also shown.}
\end{figure}

More specifically, $\mathbbm{1}_{[a,b]} (X_T)=\theta(X_T-b)-\theta(X_T-a)$ is approximated by
\begin{equation}
P^{r,\delta}_{[a,b]}(X_T) \equiv
g_r \left( \frac{X_T - b}{\delta} \right)  - g_r \left( \frac{X_T - a}{\delta} \right).
\label{eqn:SmoothedIndicator}
\end{equation}
This function is plotted for $\delta = 0.1$ and $\delta = 0.02$ and different polynomial degrees $r$ in Figs.\ \ref{Fig:0.1} and \ref{0.02}.

Since (in contrast to the indicator function) $P^{r,\delta}_{[a,b]}$ is Lipschitz continuous, we expect to recover quadratic variance decay rates in MLMC. This is confirmed numerically in Section \ref{sec:smoothing_sensitivity}.
\subsection{Error analysis}\label{eqn:IndicatorErrorAnalysis}
To quantify the error introduced by the smoothing process described in the previous section, we analyse the weak and strong errors of the problem with the smoothed indicator function. For ease of notation, we write $P(X)\equiv \mathbbm{1}_{[a,b]}(X)$ and $P^{r,\delta}(X)\equiv P^{r,\delta}_{[a,b]}(X)$ and denote by $\hat{P}_L$ and $\hat{P}_L^{r,\delta}$ the Monte Carlo estimates of $\mathbb{E}[P]$ and $\mathbb{E}[P^{r,\delta}]$, which can be obtained either using StMC or MLMC estimators with stepsize $h_L = \frac{T}{M_L} = \frac{T}{M_0 2^L}$. After some straightforward algebra (see \ref{weakApp}) we obtain a bound on the bias (or weak) error
\begin{equation}
\left| \mathbb{E}[P] - \mathbb{E}[\hat{P}^{r,\delta}_L] \right| \leq 2 c\cdot \delta^{r + 1} + \left| \mathbb{E}[P^{r,\delta}] - \mathbb{E}[P^{r,\delta}_L] \right|.
\label{weak}
\end{equation}
In this equation, the second term describes the bias error introduced by the timestepping method. It is $\propto h_L$ for the first-order methods used in this work. The total mean-square error is
\begin{gather}
\begin{split}
\mathbb{E} \left[ \left( \hat{P}^{r,\delta}_L - \mathbb{E}[P] \right)^2 \right] &\leq \text{Var}( \hat{P}^{r,\delta}_L) + \left( \mathbb{E}[P^{r,\delta}_L] -\mathbb{E}[P^{r,\delta}] \right)^2 + 4 c^2 \cdot\delta^{2(r + 1)}\\&\qquad+ 4 c\cdot \delta^{r + 1} \left| \mathbb{E}[P^{r,\delta}_L] -\mathbb{E}[P^{r,\delta}] \right|, \end{split}\label{strong}
\end{gather}
where $c, \delta$ and $r$ are as described before and the quantity $\big|\mathbb{E}[P^{r,\delta}_L] -\mathbb{E}[P^{r,\delta}]\big|\propto h_L$ is the timestepping error.

From the strong-error bound \eqref{strong}, it can be observed that it suffices to let $\delta \to 0$ or $r \to \infty$ to ensure that the last two terms are dominated by the first two terms which represent the sampling error and the discretisation error as in the case of the indicator function without smoothing. Therefore, it is possible to systematically eliminate the smoothing error for practical calculations by choosing suitable $r$ and $\delta$. However, care has to be taken since the variance of $Y_\ell$ increases as the approximation to the indicator function improves. This can be observed numerically, for $r \to \infty$, in Fig.~\ref{varianceratesconc01}. The behaviour of $\text{Var}[Y_\ell]$ for $\delta \to 0$ is similar.
\section{Numerical Results}\label{sec:Results}
We now describe several numerical experiments to systematically quantify errors and to compare the efficiencies of the single and multilevel estimators with all three timestepping methods. Results for three quantities of interest are presented: the particle position $X_T$ at the final time $T$, the concentration in an interval $[a,b]$ and a piecewise constant approximation of the concentration field over the entire boundary layer. The latter is of particular interest in practical applications. In all cases, we first check the quadratic variance decay required in condition (iii) of the complexity theorem and then calculate the total computational cost as a function of the tolerance $\epsilon$ for the root-mean-square error. For the choice of the number of time steps $M_0$ on the coarsest level, we ensure in all cases a priori that $\text{Var}[Y_1]$ is less than $\frac12 \text{Var}[Y_0]$ (if this was not the case, it would be more efficient to reduce the number of levels). For Symplectic Euler, extra care has to be taken to ensure that $h_0$ satisfies the stability constraint discussed in Section \ref{sec:TimesteppingMethods}. The number of levels~$L$ is also chosen in advance by first approximating the bias constants $\alpha$ and $c_1$ in the bias error expression $\mathbb{E}[\hat{\mathcal{P}_\ell}-\mathcal{P}] \approx c_1 h_\ell^{\alpha}$ and then making sure that this bias error is of size $\epsilon/\sqrt 2$ (see \cite{EikeRobTony} for details). The time required for calculating $\alpha$, $c_1$ and choosing the number of levels is not included in the computational costs reported below. See \cite{Giles2008,Collier2015} for cheap ways to estimate these constants adaptively on-the-fly. To limit the total mean square error to $\epsilon^2$, the second term in (\ref{eqn:MLMCerror}) is also reduced below $\epsilon^2/2$ by adaptively choosing the number of samples $N_\ell$ using on-the-fly estimators for $\text{Var}[Y_\ell]$, see Eqn. (12) in \cite{Giles2008}.

All results were generated with a freely available object-oriented C++ code developed by the authors that can be downloaded under the GPL license from {\small \url{https://bitbucket.org/em459/mlmclangevin}}.
The numerical results reported here were obtained with the version which has been archived as \cite{MLMCLangevinCode2017}. The code was compiled with version~4.8.4 of the GNU C++ compiler and run sequentially on an Intel i5-4460 CPU with a clock speed of 3.20GHz. All costs reported below refer to measured CPU times.

\begin{figure}
  \centering
  \includegraphics[width=0.75\linewidth]{\figdir/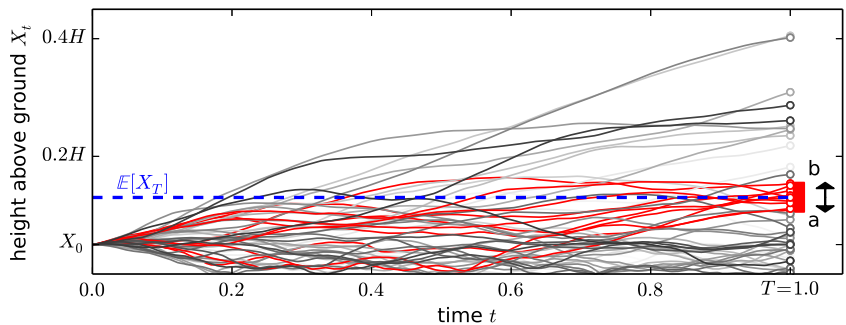}
  \caption{Sample trajectories from a typical model run. The expected value $\mathbb{E}[X_T]$ and the interval $[a,b]$ used for the concentration calculation are shown on the right. Trajectories which end in $[a,b]$ are highlighted in red.}
  \label{fig:trajectories}
\end{figure}

In the numerical experiments, we use $\kappa_{\sigma} = 1.3, \kappa_{\tau} = 0.5, u_\ast = 0.2 \operatorname{m/s}$ and a boundary-layer depth of $H = 1 \operatorname{km}$ in Eq.~(\ref{eqn:profiles}). All times, positions and velocities are expressed in units of the reference distance $X_{\text{ref}}=H=10^3\operatorname{m}$, reference velocity $U_{\text{ref}}=1\operatorname{m/s}$ and reference time $t_{\text{ref}}=X_{\text{ref}}/U_{\text{ref}}=10^3 s\approx 17$~minutes. Particles are released at $X_0 = 0.05$ (corresponding to a height of $50\operatorname{m}$ above ground) with a velocity of $U_0 = 0.1$ (i.e., an upward release speed of $0.1\operatorname{m/s}$), but we also study the impact of the release height on performance in Section \ref{sec:pdf}. Some sample trajectories for this set-up are shown in Fig.\ \ref{fig:trajectories}.

Unless explicitly stated otherwise, the regularisation height is $\epsilon_{\text{reg}}=0.01$. In this case, $\tau(X)$ does not exceed a value of $20s$, which is consistent with the choice in the Met Office NAME model. The impact of $\epsilon_{\text{reg}}$ on the results is quantified in Section 
\ref{sec:sensitivity_regularisation}.
\subsection{Particle Position}\label{sec:ParticlePosition}
We consider first the expected value of the particle position at the final time, $\mathbb{E}\left[X_T\right]$. For all experiments, we use $M_0 = 40$ timesteps on the coarsest level and a final time $T = 1$ (corresponding to $\sim 17$ minutes).

\begin{figure}
\begin{subfigure}{.5\textwidth}
\centering
\captionsetup{justification=centering}
\includegraphics[scale = 0.4]{\figdir/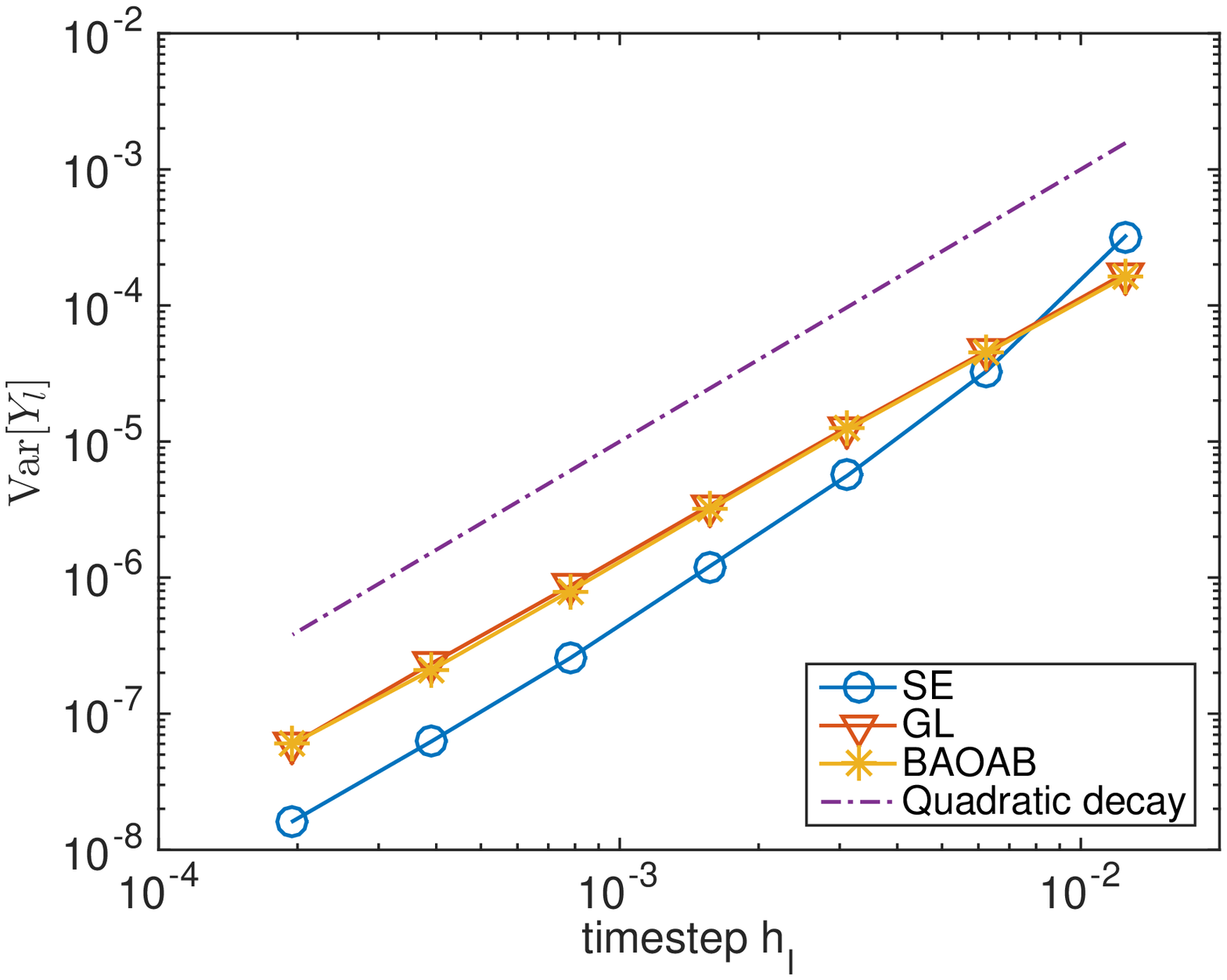}
\caption{Variance for different timestep sizes}
\label{VarianceRatesPosition}
\end{subfigure}%
\begin{subfigure}{.5\textwidth}
\centering
\captionsetup{justification=centering}
\includegraphics[scale = 0.4]{\figdir/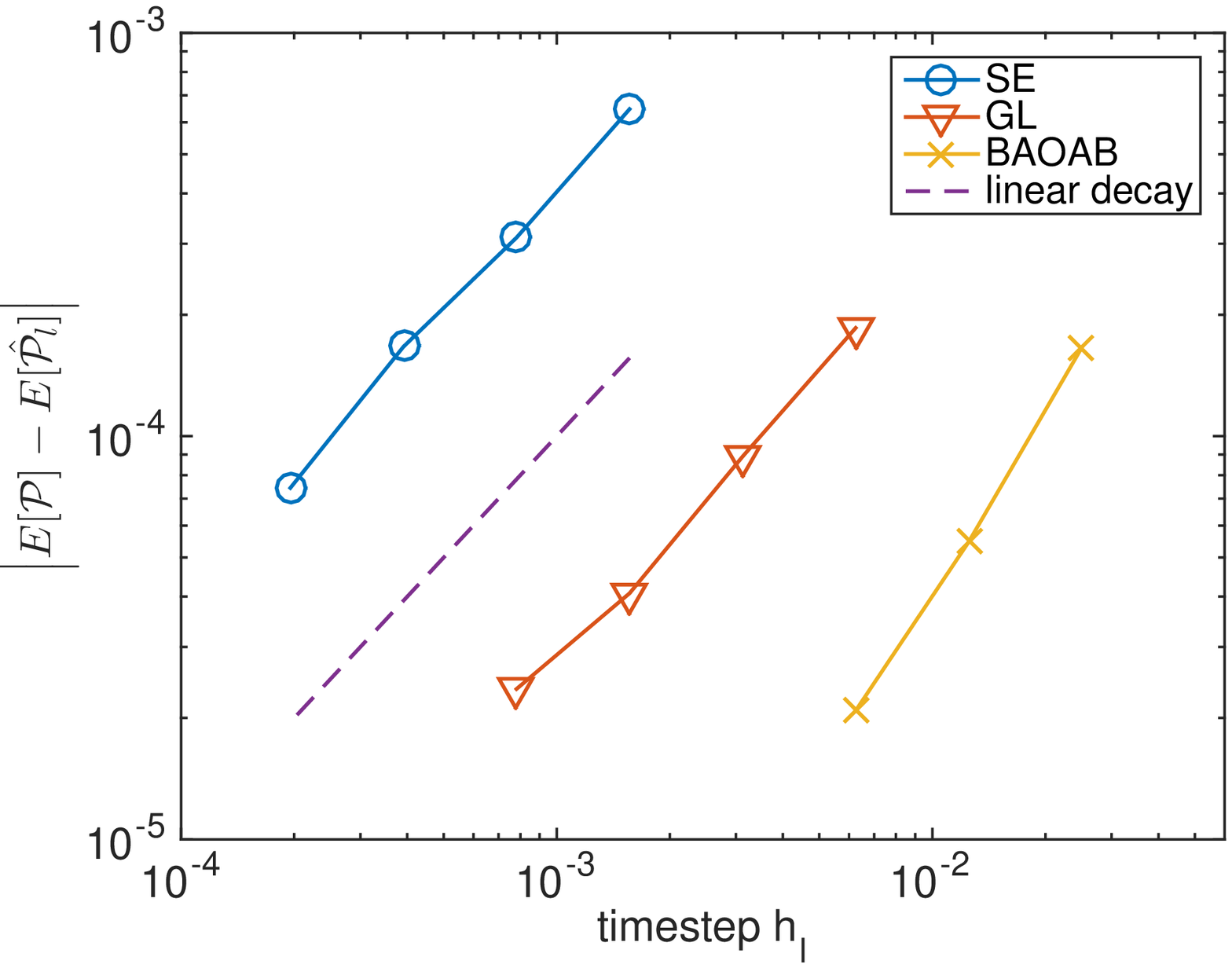}
\caption{Bias for different timestep sizes}
\label{BiasRatesPosition}
\end{subfigure}%
\caption{Variance and total computational cost - Mean particle position}
\end{figure}

The validity of condition (iii) in the complexity theorem is confirmed numerically in Fig.\ \ref{VarianceRatesPosition} where we plot the variance of $Y_{\ell}=X_{T,h_\ell}-X_{T,h_{\ell-1}}$ against $h_\ell$. For all numerical methods, the rate is quadratic, $\text{Var}[Y_\ell]\propto h_\ell^2$. We observe that the geometric Langevin (GL) and BAOAB methods have essentially an identical variance in absolute terms. The variance of the Symplectic Euler (SE) method is slightly smaller, except for the largest timestep sizes. As we can see  in Fig.\ \ref{BiasRatesPosition}, condition (i) is also satisfied for all methods. We note, however, that the splitting methods lead to a dramatic reduction in bias error compared to SE, for any fixed value of $h_\ell$. The bias error for the GL method is about $13\times$ smaller than that of SE. The BAOAB integrator reduces this by an additional factor of $4$ leading to a total reduction of $52\times$ relative to SE.

\begin{figure}
\begin{subfigure}{.5\textwidth}
\centering
\captionsetup{justification=centering}
\includegraphics[scale = 0.4]{\figdir/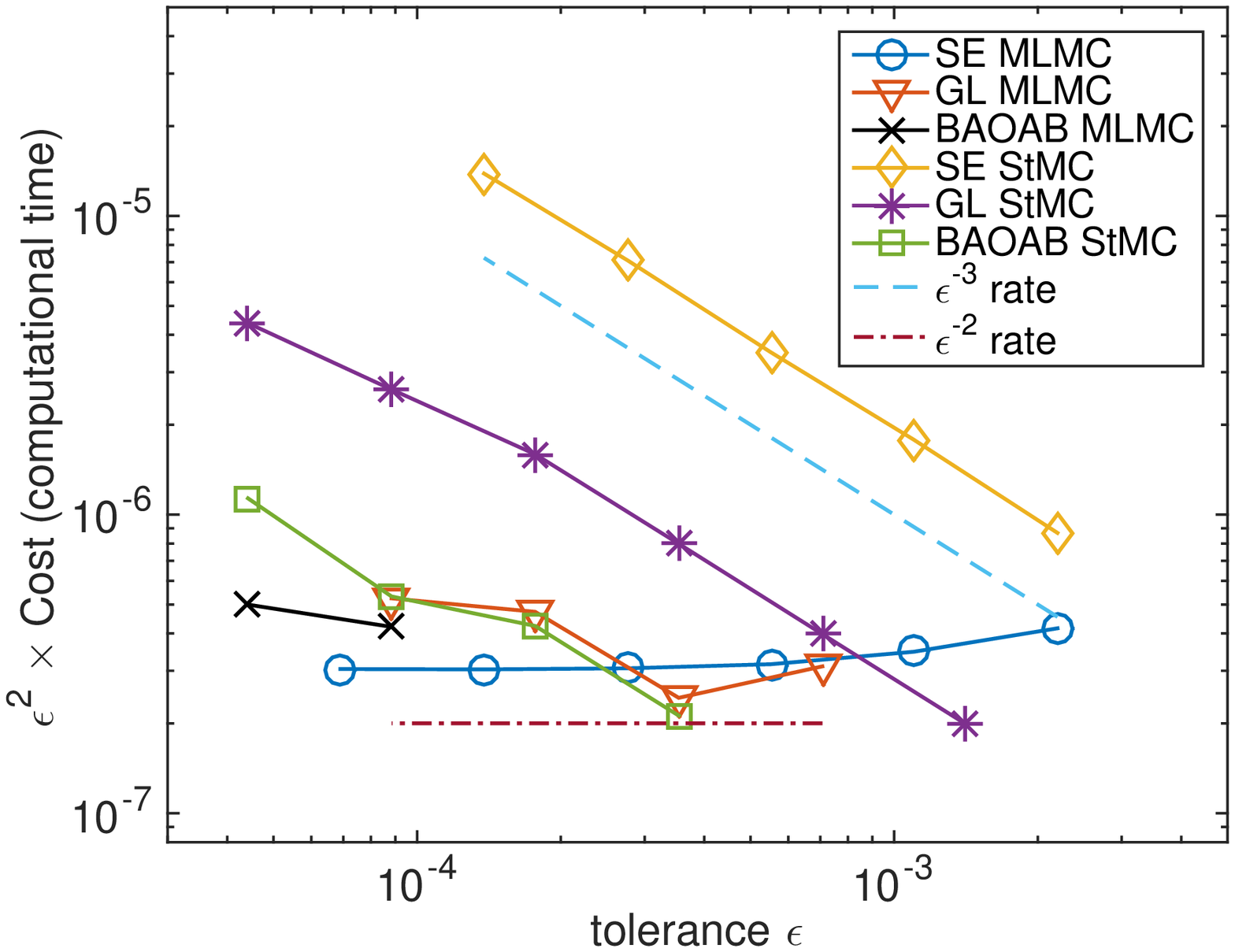}
\caption{Mean particle position}
\label{CostRatesPosition}
\end{subfigure}%
\begin{subfigure}{.5\textwidth}
\centering
\captionsetup{justification=centering}
\includegraphics[scale = 0.4]{\figdir/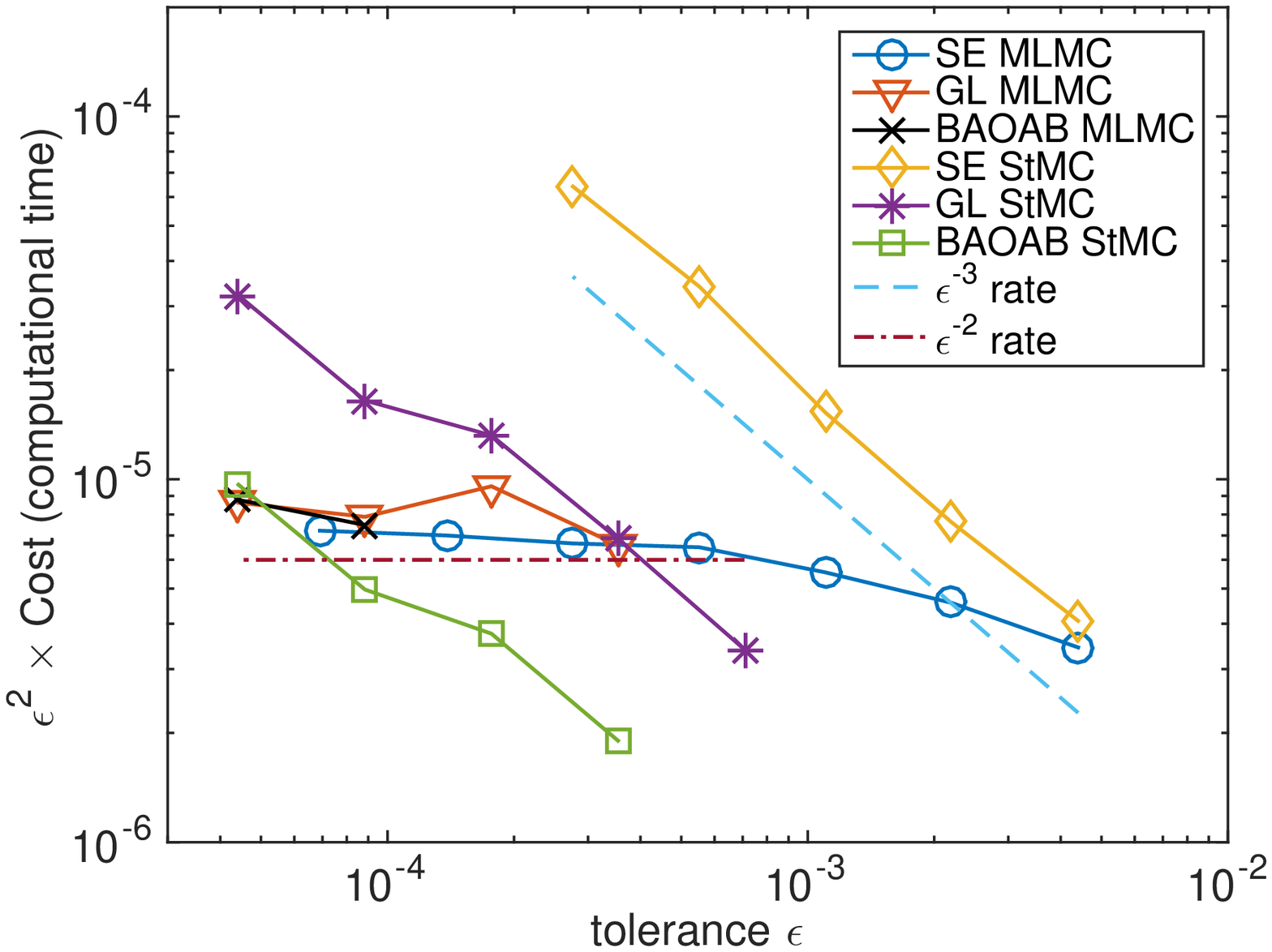}
\caption{Concentration using $P^{r=4, \delta=0.1}_{[a,b]}(X_T)$}
\label{CostRatesConcentration}
\end{subfigure}%
\caption{Computational cost for varying tolerance}
\end{figure}

In Fig.\ \ref{CostRatesPosition}, the total computational costs of the StMC and MLMC estimators for all discretisation methods are plotted as functions of the root-mean-square error tolerance $\epsilon$. The StMC method shows an asymptotic scaling proportional to $\epsilon^{-3}$ for all discretisation schemes, as discussed in Section \ref{sec:StMCcostanalysis}. The cost of the MLMC method is proportional to $\epsilon^{-2}$ in all cases, as predicted by the complexity theorem in Section \ref{sec:MLMCcostanalysis}. However, for a fixed $\epsilon$ the balance between bias and statistical error is different for the three discretisation methods: for fixed $h_\ell$, the GL and BAOAB methods have a smaller bias but larger variance error, so that, for a given $\epsilon$, the SE method will require a smaller timestep. For the SE method, the maximal allowed timestep size is restricted by stability constraints, whereas for the GL and BAOAB algorithms there are no such restrictions.

To give the reader an idea of the associated timestep sizes used in Fig.\ \ref{CostRatesPosition}, note that for the SE method the timestep size varies from $h_L = 3.1\cdot10^{-3} (\epsilon=2.2\cdot10^{-3})$ to $h_L = 9.8\cdot10^{-5} (\epsilon=6.9\cdot10^{-5})$. For the GL method the corresponding range in timestep sizes is $h_L = 2.5\cdot10^{-2} (\epsilon=1.4\cdot10^{-3})$ to $h_L = 7.8\cdot10^{-4} (\epsilon=4.4\cdot10^{-5})$. For BAOAB the largest timestep is $h_L = 2.5\cdot10^{-2} (\epsilon=3.5\cdot10^{-4})$ and the smallest $h_L = 3.1\cdot10^{-3} (\epsilon=4.4\cdot10^{-5})$. For each method the StMC stepsize and the MLMC stepsize on the finest level agree for a given root mean square error $\epsilon$. In \ref{DataOnLevels} we provide further details on the number of samples $N_\ell$ on each level $\ell$ for two representative MLMC runs.

Comparing the MLMC results among themselves, all discretisation schemes perform similarly with only small differences and SE being slightly cheaper. For the StMC algorithm, on the other hand, the splitting methods are clearly superior to the reference SE implementation. The GL integrator is approximately $5\times$ faster across the range of tolerances $\epsilon$ considered, while the BAOAB method is more than $20\times$ faster. We conclude that for small tolerances MLMC (with any discretisation scheme) leads to significant speed-ups, while for larger tolerances ($> 5 \cdot 10^{-4}$) the StMC method with BAOAB timestepping is the most efficient.

\subsection{Concentration}\label{sec:Concentration}
To estimate the concentration in a given interval $[a,b]\subset [0,1]$, we calculate the expected value of the indicator function at the final time, i.e. $\mathbb{E}\left[\mathbbm{1}_{[a,b]}(X_T)\right]$. The following results were obtained with an interval of width $0.05$ centred at $0.1305$ (i.e.\ $a=0.1055$ and $b=0.1555$; see Fig.\ \ref{fig:trajectories}). The centre of the interval roughly matches the expected value of $X_T$ calculated in the previous section, $\mathbb{E}[X_T]=0.1301\pm 4\cdot10^{-4}$.
\subsubsection{Sensitivity to smoothing parameters}\label{sec:smoothing_sensitivity}

\begin{figure}
\begin{subfigure}{.5\textwidth}
\centering
\includegraphics[scale = 0.4]{\figdir/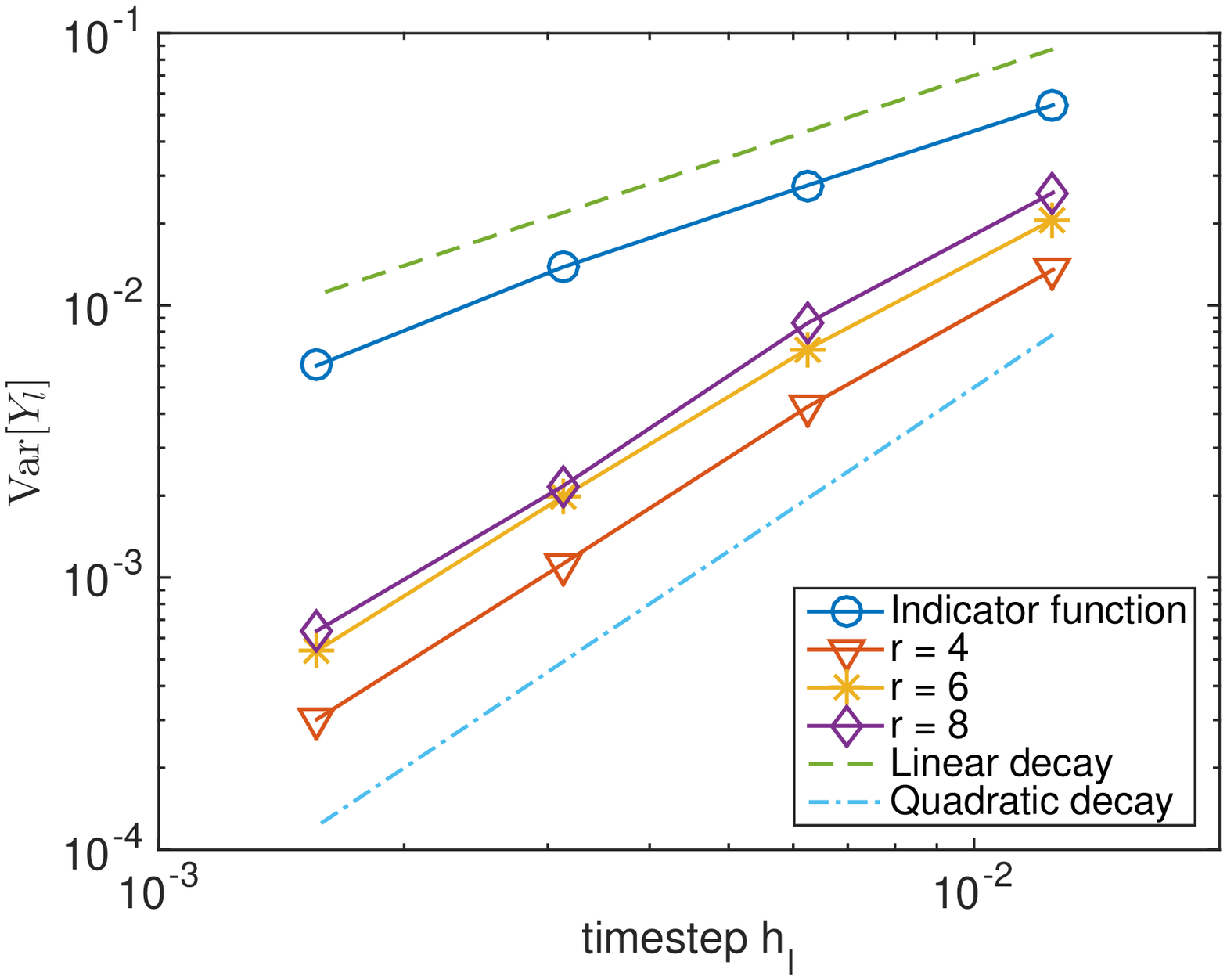}
\caption{Geometric Langevin}
\label{varianceratessymplecticeulerOU01}
\end{subfigure}%
\begin{subfigure}{.5\textwidth}
\centering
\includegraphics[scale = 0.4]{\figdir/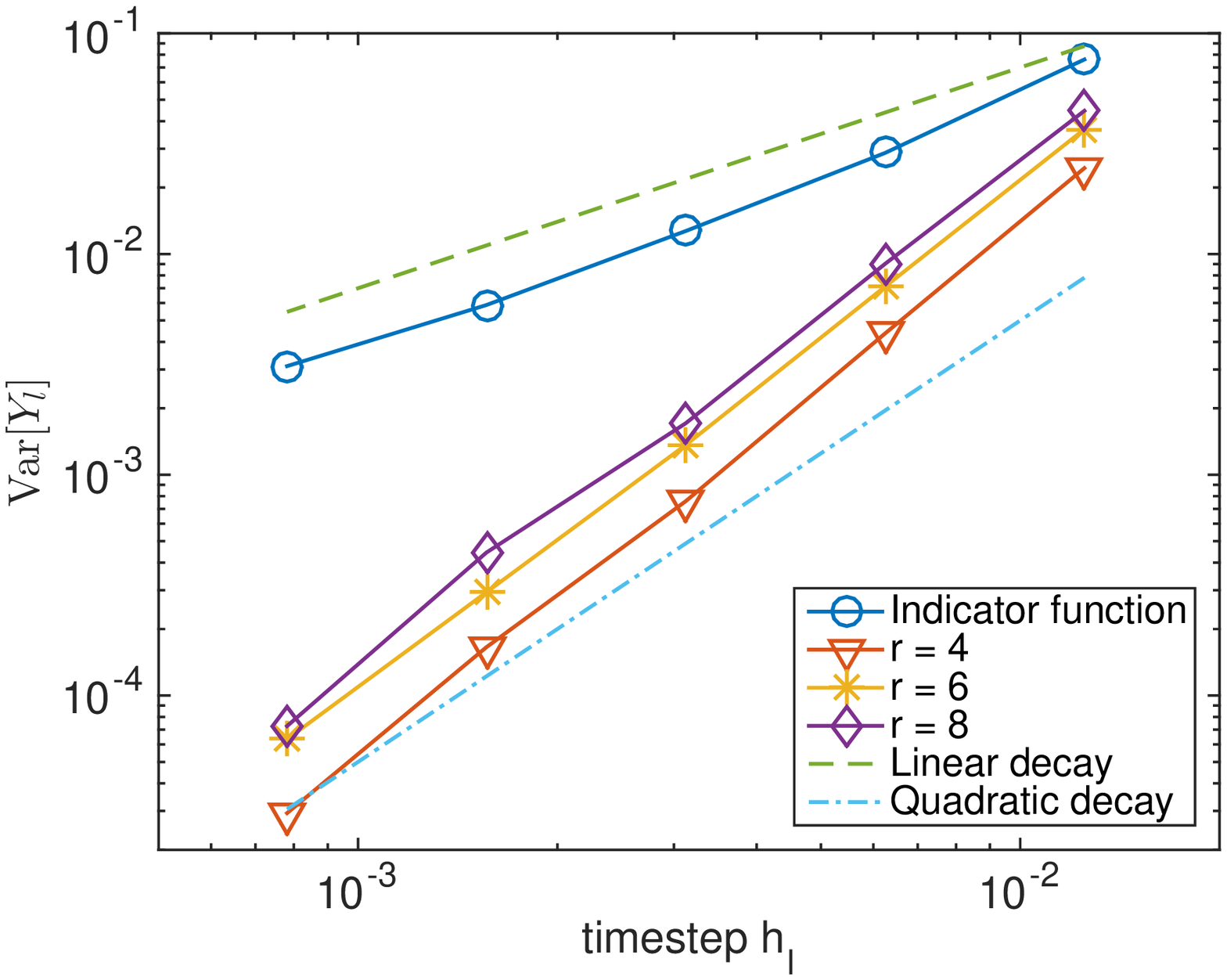}
\caption{Symplectic Euler}
\label{varianceratessymplecticeuler01}
\end{subfigure}%
\caption{Variance - Concentration $P_{[a,b]}^{r,\delta}$ with $\delta = 0.1$}
\label{varianceratesconc01}
\end{figure}
As discussed in Section \ref{sec:SmoothingPolynomials}, to estimate concentrations with the MLMC algorithm it is useful to replace the discontinuous indicator function $\mathbbm{1}_{[a,b]}$ by the smooth function $P^{r,\delta}_{[a,b]}$ defined in Eq.\  (\ref{eqn:SmoothedIndicator}). In the following, we numerically quantify the errors introduced by this approximation and confirm the results from the theoretical analysis in Section \ref{eqn:IndicatorErrorAnalysis}. To verify condition (iii) in the complexity theorem and show that smoothing indeed leads to the correct variance decay, we plot  in Figs.~\ref{varianceratessymplecticeulerOU01} and \ref{varianceratessymplecticeuler01} the variances of $Y_{\ell}$, for the GL and the SE methods, as functions of the timestep size $h_\ell$. As discussed in Section \ref{sec:ParticlePosition}, the SE method has a larger bias error compared to GL and therefore we use a smaller finest timestep size. The variance rates for BAOAB are not included since they are very similar to those of GL. The quantity of interest is $P^{r,\delta}_{[a,b]}(X_T)$ with $\delta=0.1$, for three different values of $r$, as well as the indicator function $\mathbbm{1}_{[a,b]}$. These are the functions plotted in Fig.~\ref{Fig:0.1}. For both plots the tolerance on the root-mean-square is fixed to $\epsilon = 10^{-3}$. We observe quadratic dependence $\text{Var}[Y_\ell]\propto h_\ell^2$ for the smoothed QoIs, while for the non-smooth indicator function the variance depends linearly on $h_\ell$. The absolute value of the variance grows as $r$ increases and this is intuitively expected since the function $P_{[a,b]}^{r,\delta}$ converges to the non-smooth indicator function as $r \to \infty$.

To show that the error introduced by replacing the indicator function by a smooth approximation is indeed small, we tabulate estimates of the expected values of the indicator function and of $P^{r,\delta}_{[a,b]}$, for different values of $r$, in Tab.\ \ref{tab:comparison_to_indicatorfunction}. We choose the time step size $h$ such that the discretisation error for each method is less than $\epsilon/\sqrt{2}$ with $\epsilon = 10^{-3}$, and choose a large number of samples, such that the standard deviation of each estimator is less than $3.5\cdot 10^{-5}$. As can be seen from the last three columns of Tab.\ \ref{tab:comparison_to_indicatorfunction}, the difference between the indicator function and the smoothed approximations $P^{r,\delta}_{[a,b]}$ is of the order $5\cdot 10^{-5}$ in all cases. We conclude that (compared to the bias) the smoothing error is clearly negligible for the choices of $\delta$ and $r$ in Tab.\ \ref{tab:comparison_to_indicatorfunction}. 

\begin{table}
  \centering
  \begin{tabular}{lccccccc}
  Quantity of interest & $r$ & \multicolumn{3}{c}{value} & \multicolumn{3}{c}{difference}\\
   & & GL & SE & BAOAB & GL & SE & BAOAB\\
  \hline\hline
  & $4$ & 0.16678 & 0.16705 & 0.16687 & $1.2\cdot10^{-5}$ & $7.4\cdot 10^{-5}$ & $3.1\cdot10^{-5}$ \\
  $\mathbb{E}\left[P^{r,\delta}_{[a,b]}(X_T)\right]$ & $6$ & 0.16669 & 0.16709 & 0.16677 & $7.2\cdot10^{-5}$ & $4.1\cdot10^{-5}$ & $7.5\cdot10^{-5}$ \\
  & $8$ & 0.16682 & 0.16714 & 0.16683 & $5.1\cdot 10^{-5}$ & $9.9\cdot 10^{-6}$ & $1.6\cdot 10^{-5}$ \\
  \hline
  $\mathbb{E}\left[\mathbbm{1}_{[a,b]}(X_T)\right]$ & & 0.16677 & 0.16713 & 0.16684 \\
  \end{tabular}
  \caption{Estimates of the expected values for the indicator function and the smoothed function $P^{r,\delta}$ for different polynomial degrees $r$; in all cases $\delta=0.1$ and the tolerances on the bias and sampling errors were $\epsilon/\sqrt{2}=7.1\cdot 10^{-4}$ and $3.5\cdot 10^{-5}$, respectively. The last three columns show the difference between the smoothed functional and the indicator function, $\big|\mathbb{E}\left[P^{r,\delta}_{[a,b]}(X_T)\right]-\mathbb{E}\left[\mathbbm{1}_{[a,b]}(X_T)\right]\big|$; this quantity has a sampling error which is bounded by $5.0\cdot 10^{-5}$.}
  \label{tab:comparison_to_indicatorfunction}
\end{table}

The cost for evaluating $P^{r,\delta}_{[a,b]}(X_T)$ is also negligible. It is less than the cost of one time step and only depends very mildly on $r$. In addition, we have also confirmed that the behaviour observed in Tab.\ \ref{tab:comparison_to_indicatorfunction} is independent of the choice of the interval $[a,b]$ and thus, we will choose $\delta=0.1$ and $r=4$ in all subsequent experiments.

\subsubsection{Cost comparison between Standard and Multilevel Monte Carlo}
Neglecting the smoothing error, we now compare the cost of the StMC and MLMC methods as a function of the tolerance on the root-mean-square error. Based on the results in the previous section, we fix $\delta = 0.1$ and $r = 4$. For the SE and the GL methods, we choose $M_0 = 80$; for BAOAB, we choose $M_0 = 40$.

In Fig.\ \ref{CostRatesConcentration}, the total computational cost for the StMC and MLMC method is plotted for different tolerances on the root-mean-square error $\epsilon$. We see a similar picture as in Section \ref{sec:ParticlePosition}. While for larger tolerances, the StMC with GL or BAOAB are cheaper than the MLMC algorithms, for tighter tolerances the MLMC algorithms become superior with an asymptotic growth of the cost at the expected rate of $\epsilon^{-2}$. We observe the expected growth $\propto \epsilon^{-3}$ for the three StMC methods. For the MLMC methods, there is again little dependence on the discretisation method. On the other hand, comparing the StMC algorithms, for a given $\epsilon$, there is again a large difference in the efficiency between the three discretisation methods, with the BAOAB method being the fastest.

The reason is again, as in the previous section, that the bias error is significantly smaller for the splitting methods. As for the particle position, the bias error for the GL method is around $13\times$ smaller than the bias for the SE algorithm, for a given timestep size $h_\ell$, while the bias for BAOAB is even $52\times$ smaller. Details for the number of samples on each MLMC level can be found in \ref{DataOnLevels}.

\subsection{Probability Density Function}\label{sec:pdf}
Since the main output of operational dispersion models like NAME is a box-averaged concentration field, we finally calculate an approximation to the probability density function (p.d.f.) at the final timestep. For this a set of equidistant points $0=a_0<a_1<a_2<\dots<a_{k-1}<a_k=H$ is chosen and the particle concentration is estimated for each interval (or bin) $[a_i,a_{i+1}]$. In the Monte Carlo algorithm, this is achieved by considering the following vector-valued quantity of interest:
\begin{equation}
  \vec{\mathcal{P}} = (\mathcal{P}_0,\mathcal{P}_1,\dots,\mathcal{P}_k)=\left(P^{r,\delta}_{[a_0,a_1]}(X_T),P^{r,\delta}_{[a_1,a_2]}(X_T),\dots,P^{r,\delta}_{[a_{k-1},a_k]}(X_T)\right),
\label{vector_QoI}
\end{equation}
which is related to the p.d.f. or particle concentration $\rho(X)$ by
\begin{equation*}
\int_a^b \rho(X)\;dX = \mathbb{E}[\mathbbm{1}_{a,b}] \approx\mathbb{E}[P_{[a,b]}^{r,\delta}]\qquad \text{for all intervals}\;[a,b]\subset [0,H].
\end{equation*}
As above, the indicator functions are replaced by smoothed polynomials $P^{r,\delta}_{[a,b]}$ with $\delta = 0.1$ and $r = 4$ to guarantee a quadratic variance decay in the MLMC algorithm. The QoI is now vector-valued, but a scalar quantity is required to estimate the optimal number of samples $N_\ell$ on each level of the MLMC algorithm (see \cite{Giles2008,Giles2015} for details); we choose the maximal absolute variance over all intervals for this scalar. More specifically, let $\mathcal{P}_{\ell,i}=P_{[a_i,a_{i+1}]}^{r,\delta}(X_{T,\ell})$ and $Y_{\ell,i} = \mathcal{P}_{\ell,i}-\mathcal{P}_{\ell-1,i}$, for $\ell=0,\ldots,L$ and $i=0,\ldots,k$. Then, we use $ \max_{i=0,k}|\text{Var}[Y_{\ell,i}]|$ in place of $\text{Var}[Y_\ell]$ to estimate $N_\ell$ on level $\ell$ .

\begin{figure}
\centering
\includegraphics[scale = 0.25]{\figdir/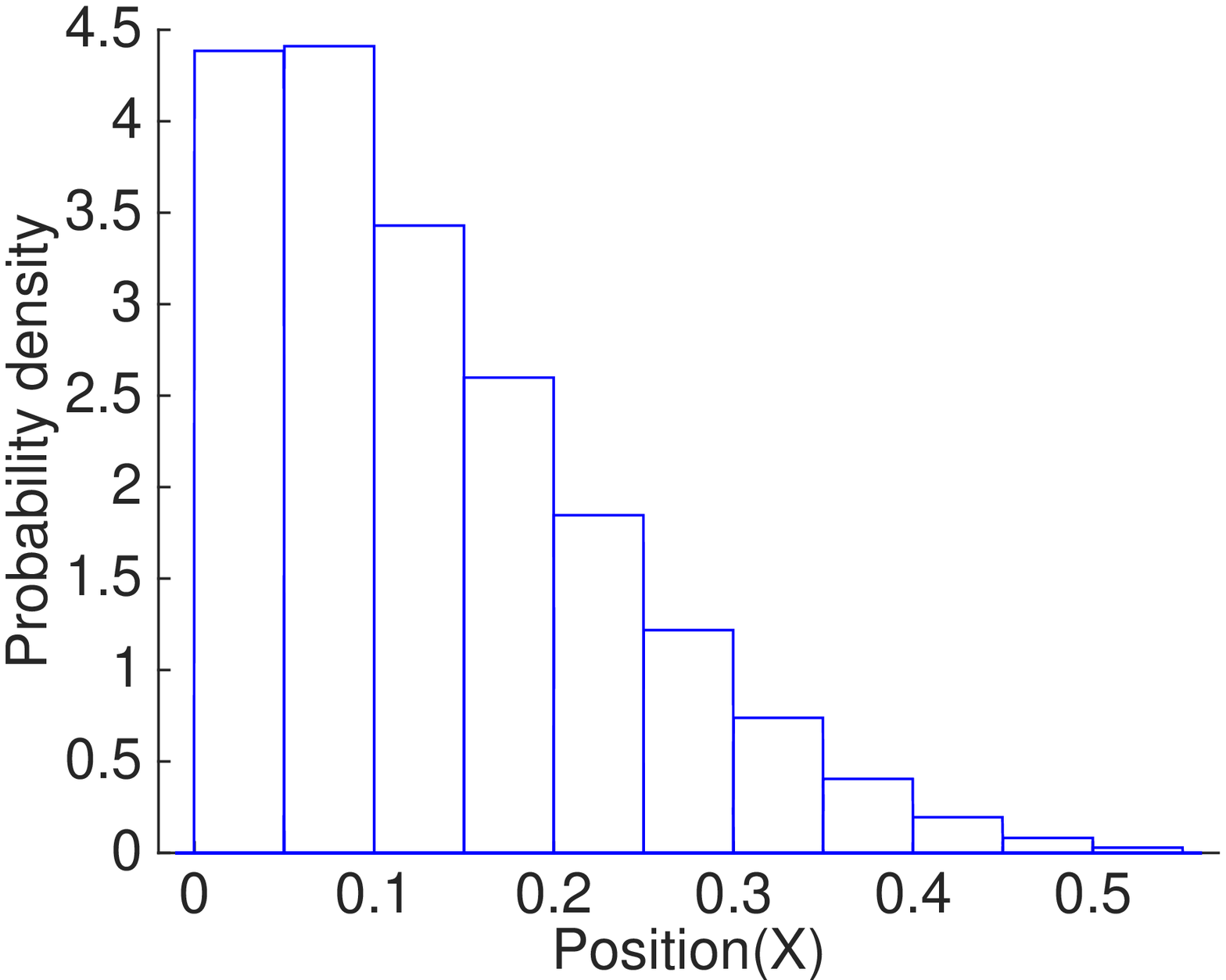}
\includegraphics[scale = 0.25]{\figdir/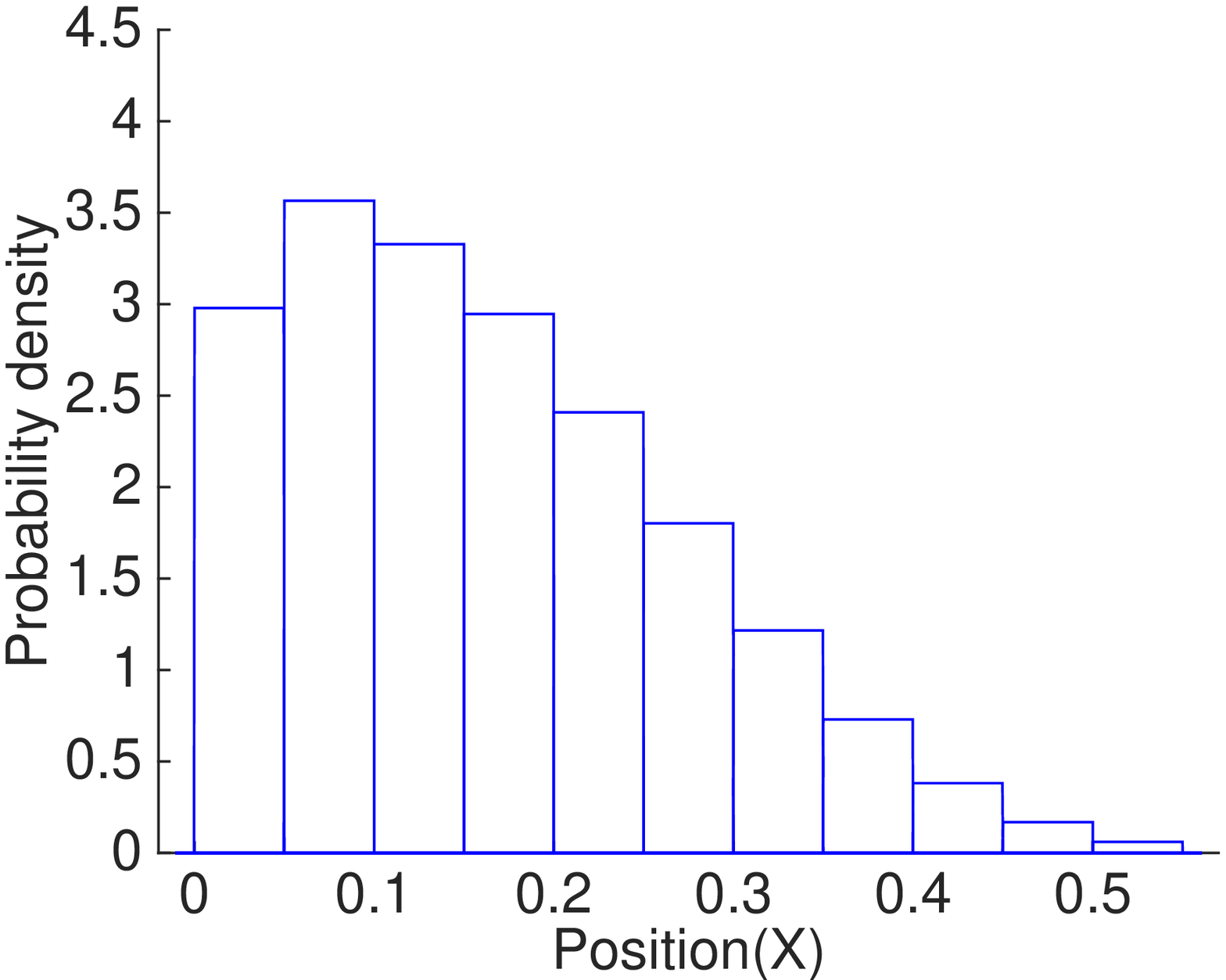}
\includegraphics[scale = 0.25]{\figdir/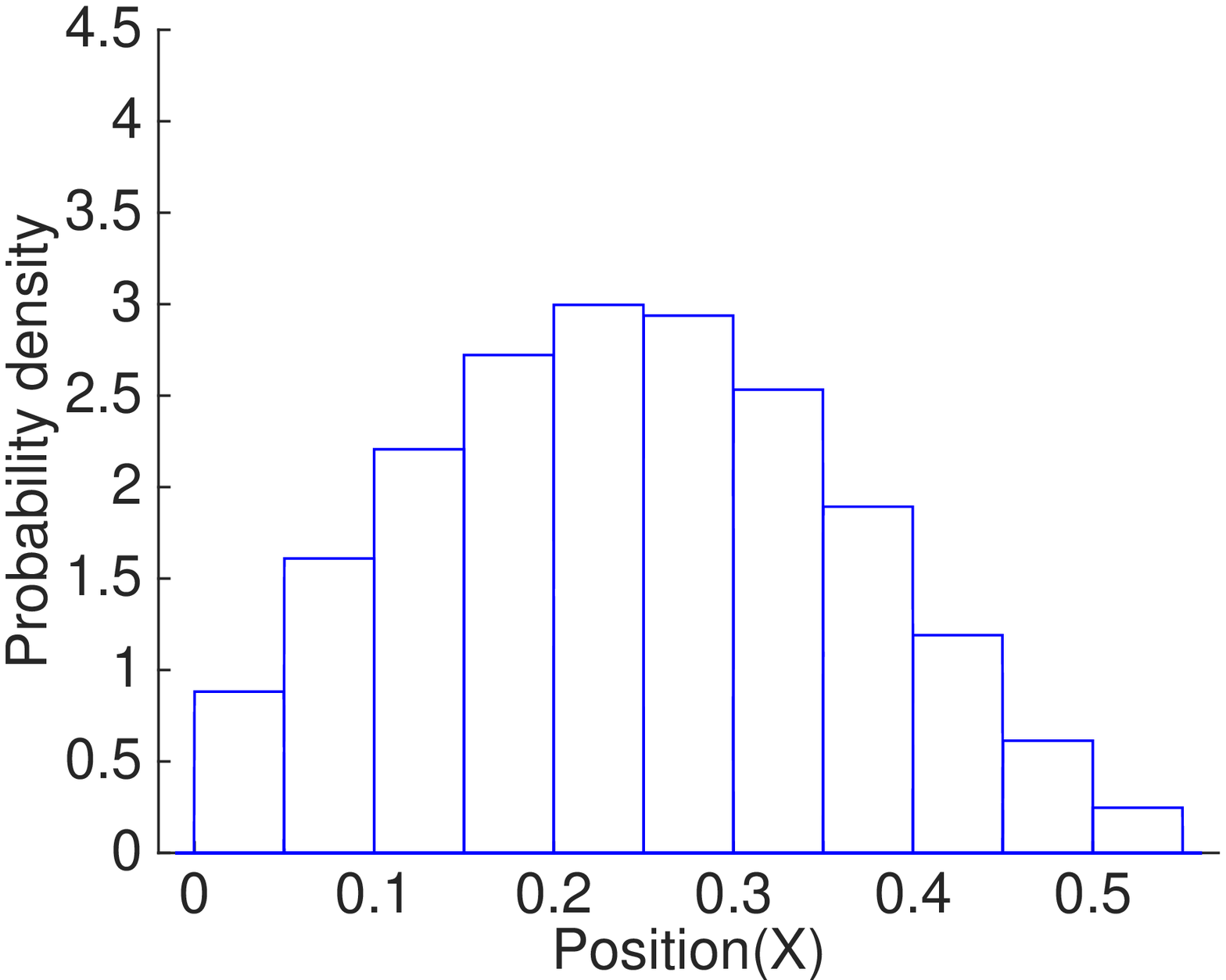}
\caption{Binned concentration field for release height $0.05$ (left), $0.10$ (center) and $0.20$ (right). The results were obtained with the MLMC method and symplectic Euler timestepping.}
\label{fig:concentration}
\end{figure}
Fig.\ \ref{fig:concentration} shows the binned probability density function of $X_T$ obtained with the SE  method for a box size of 0.05. The leftmost plot is obtained with $X_0=0.05$, corresponding to a release height of $50\operatorname{m}$ in physical units; this is the same initial condition which was used to obtain all previous results in Sections \ref{sec:ParticlePosition} and \ref{sec:Concentration}. However, since we found that the computational cost is very sensitive to the release height, we also report results for $X_0=0.1$ and $X_0=0.2$, corresponding to $100\operatorname{m}$ and $200\operatorname{m}$ above ground, respectively.

We require that the root-mean-square error for the particle concentration $\rho(X)$ does not exceed  a tolerance of $\epsilon = 4 \cdot 10^{- 3}$. As before the number of timesteps on the coarsest level is $M_0=80$ for SE and GL and $M_0=40$ for BAOAB.
We observe that the bias error of SE decreases with increasing release height, since the particles spend less time near the ground where $\tau(X)$ varies strongly. Hence as $X_0$ increases, the size of the timestep $h_L$ that is required to reduce the bias error below the tolerance $\epsilon/\sqrt{2}$ on the finest level also increases. For GL and BAOAB the bias error stays approximately the same, and so the release height has almost no influence on the cost for those two methods. The computational costs of all the different estimators to reduce the root-mean-square error below a tolerance of $\epsilon = 4 \cdot 10^{- 3}$ are shown in Tab.\ \ref{tab:ReleaseHeightResults}.
\begin{table}
  \centering
  \begin{tabular}{l|rrr|rr}
    release height  & \multicolumn{3}{|c|}{StMC} & \multicolumn{2}{|c}{MLMC} \\
    & SE & GL & BAOAB & SE & GL \\
    \hline\hline
    $0.05$ ($50\operatorname{m}$) & 2561.41 & 321.77 & 91.21 & 218.52 & 265.84 \\
    $0.10$ ($100\operatorname{m}$) & 867.32 & 277.72 & 79.25 & 120.17 & 125.92 \\
    $0.20$ ($200\operatorname{m}$) & 316.33 & 280.04 & 78.15 & 64.78 & 103.69 \\ 
  \end{tabular}
  \caption{Computational cost (measured runtimes in seconds) for different release heights. Results are shown for all timestepping methods and both Standard- and (for SE and GL) Multilevel- Monte Carlo.}
  \label{tab:ReleaseHeightResults}
\end{table}

Looking at the StMC algorithms, we note again that the performance is improved dramatically by using GL or BAOAB methods instead of SE. This improvement is particularly dramatic for small release heights where the SE method is about $28\times$ slower than BAOAB method and about $8\times$ slower than GL method, which agrees with the speed-ups observed for other quantities of interest in the previous sections. The reason for this is the reduction in bias error due to better treatment of the large and strongly varying profile $\lambda(X)$ by the GL or BAOAB methods which enables us to use larger timesteps. For higher release points the particles spend less time near the ground where $\lambda(X)$ varies strongly so that the effect on the SE method is less pronounced. As stated above, the GL and BAOAB methods are fairly insensitive to this change.

For the MLMC algorithms, the total runtime is also reduced as the release height increases. We did not include any results for BAOAB, since the number of levels is too small to achieve any speed-up for the chosen tolerance. This would be different for smaller tolerances. For the SE and the GL methods, the MLMC methods lead to significant gains over their single level counterparts in all cases. Not surprisingly, the biggest speed-up is achieved in the case of the SE method for low release heights, where MLMC is about $12\times$ faster than StMC. Overall, the BAOAB method with StMC estimator performs best for this tolerance, especially for low release heights $X_0$, but for the higher release height of $X_0=0.2$, it is in fact the SE method with MLMC estimator that is the fastest.

In summary, compared to the reference StMC-SE algorithm used in many Lagrangian atmospheric dispersion models, both the new GL and BAOAB timestepping methods and the MLMC estimator reduce the runtime significantly when computing concentration fields.

\section{Adaptive timestepping}\label{sec:Adaptivity}
Since the velocity decorrelation time $\tau(X)$ varies strongly over the boundary layer and limits the timestep size near the ground, we also investigate potential improvements from adaptive timestepping methods. By using small timesteps only in those parts of the domain where $\tau(X)$ is large, those methods can potentially reduce the overall runtime. For StMC adaptive timestepping is easy to implement and is currently used in the NAME dispersion model. Adaptive timestepping methods can also be used in the MLMC algorithm. We briefly review the technique proposed in \cite{Giles2016} and adapt it for the timestepping methods used here.
\subsection{Sensitivity to regularisation height}\label{sec:sensitivity_regularisation}
\begin{figure}
\begin{subfigure}{.5\textwidth}
\centering
\includegraphics[scale = 0.4]{\figdir/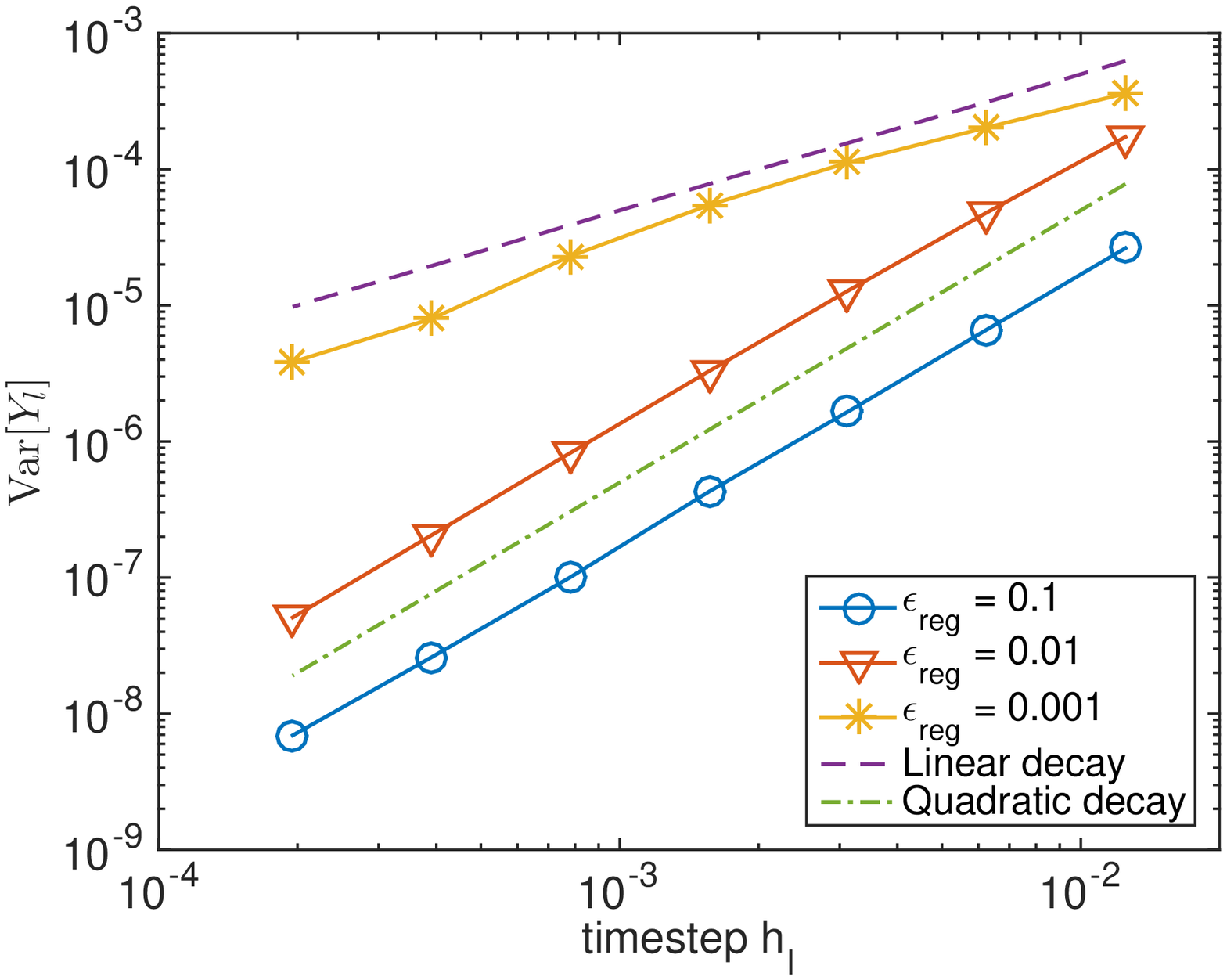}
\caption{Geometric Langevin}
\label{SymEulerOUReg}
\end{subfigure}%
\begin{subfigure}{.5\textwidth}
\centering
\includegraphics[scale = 0.4]{\figdir/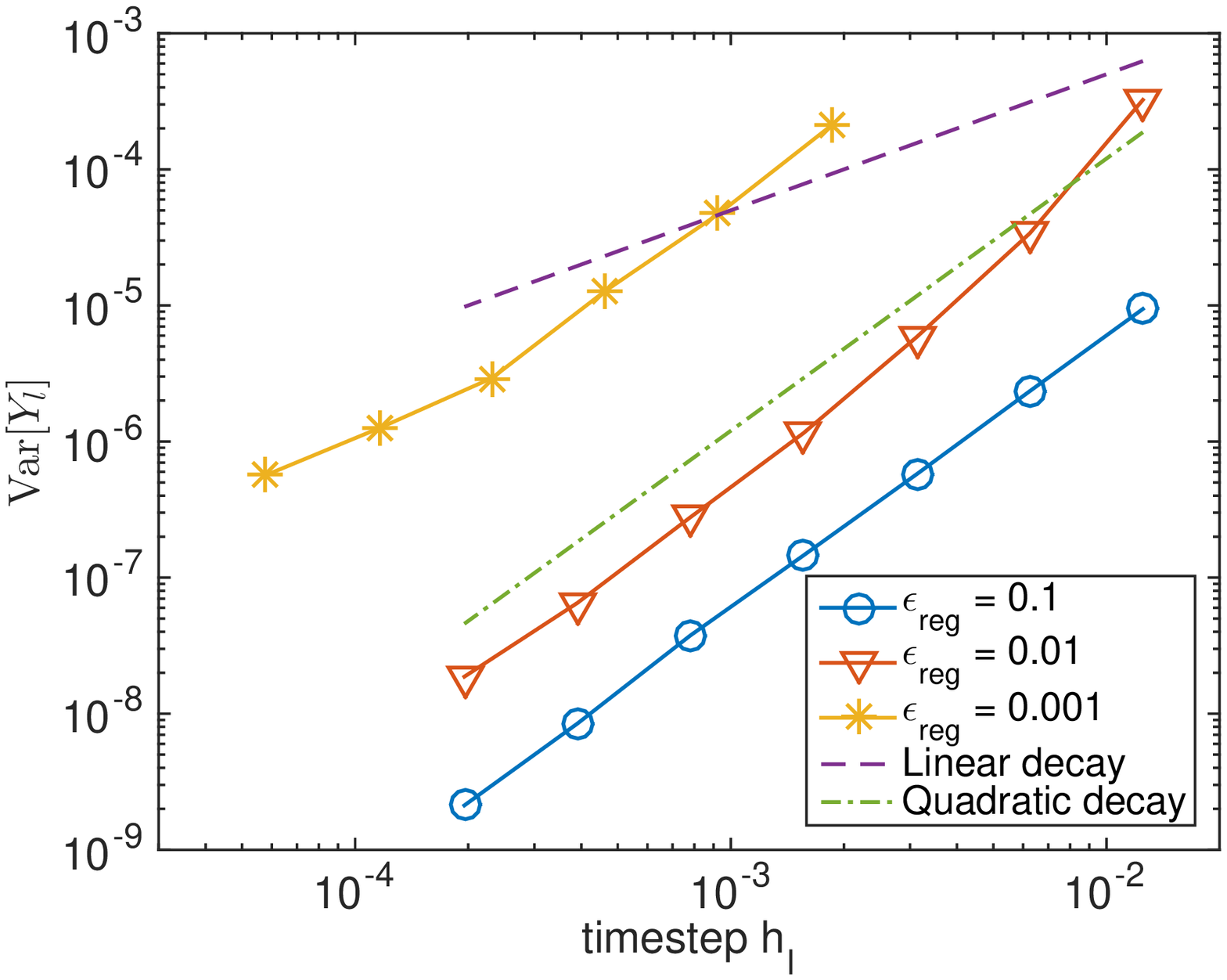}
\caption{Symplectic Euler}
\label{SymEulerReg}
\end{subfigure}%
\caption{Variance rates - Mean particle position}
\label{e_reg_sensitivity}
\end{figure}

As the largest possible value of $\tau(X)$ depends on the regularisation height $\epsilon_{\text{reg}}$ introduced in Section \ref{sec:regularisation}, we first investigate the impact of this parameter on our results. As before, the first thing to check is the quadratic variance decay $\text{Var}[Y_\ell]\propto h_\ell^2$ required for condition (iii) of the complexity theorem.
In Figs.\ \ref{SymEulerOUReg} and \ref{SymEulerReg}, the variance is shown as a function of timestep size $h_\ell$ for SE and GL timestepping methods as $\epsilon_{\text{reg}}$ is varied. In both cases, the quantity of interest is the final particle position $X_T$ and the same set-up as in Section \ref{sec:ParticlePosition} is used. For GL, the number of timesteps on the coarsest level is $M_0=40$ and up to $L=7$ levels are used. 
For SE, the timestep sizes are limited by stability constraints, i.e.\ $h < C/ \max_{X\in[0,H]} \left\{ \lambda(X) \right\} = C/ \lambda (\epsilon_{\text{reg}})$ for some constant $C$. This restriction becomes more severe as the size of the regularisation parameter decreases. For $\epsilon_{\text{reg}} = 0.1$ and $\epsilon_{\text{reg}} = 0.01$, we use the same set-up as for the GL methods, but for $\epsilon_{\text{reg}} = 0.001$ a larger number of $M_0 = 270$ on the coarsest level and a reduced number of $L = 6$ is used to ensure that the method is stable on all levels.

From those plots, it can be observed that the variance decay is indeed quadratic for $\epsilon_{\text{reg}}\ge0.01$ for the values of $h_\ell$ considered, but $\text{Var}[Y_\ell]\propto h_\ell$ for the smallest $\epsilon_{\text{reg}}$. This deterioration in the variance decay rate can be explained using the theory of modified equations \cite{shardlow2006modified,zygalakis2011existence,EikeRobTony} in \ref{sec:ModifiedEquations}: since $\lambda \rightarrow \infty$ for $X\rightarrow 0$, the variance decay rate is quadratic only if $h_\ell\cdot\epsilon_{\text{reg}}^{-2}$ is not too large. This is because for larger $h_\ell$ such that $h_\ell\cdot\epsilon_{\text{reg}}^{-2} \gtrsim 1$, heuristically we can only hope for linear variance decay with $\text{Var}[Y_\ell] \propto h_\ell$ since the size of the term $h_\ell v_1^{(U)}$ in Eq. (\ref{eqn:ModifiedEquations}) is $\mathcal{O}(1)$ and not $\mathcal{O}(h_\ell)$ and therefore the modified equation expansion breaks down. However, asymptotically, as $h_\ell$ goes to $0$ (for fixed $\epsilon_{\text{reg}})$, we recover the expected quadratic variance decay rate $\text{Var}[Y_\ell] \propto h_\ell^2$ as soon as $h_\ell\cdot\epsilon_{\text{reg}}^{-2}\ll 1$. We stress, however, that for the regularisation height used in all experiments above, which roughly agrees with what is used in the Met Office NAME model, the variance indeed behaves like $\text{Var}[Y_\ell]\propto h_\ell^2$ for practically relevant values of $h_\ell$ and hence condition (iii) in the complexity theorem is satisfied.

For both discretisation methods the main reason for the increase in the variance is the large value of $\lambda$ relative to the timestep size when the particle is close to the lower boundary. By decreasing the timestep size, the variance decreases but at the same time the total cost increases. We therefore now look at how this effect can be mitigated with adaptive timestepping.
\subsection{Adaptive timestepping}
Various methods for adaptive timestepping in the MLMC method have been proposed in the literature \cite{Giles2016,Hoel2012,hoel2014implementation,Fang2016}.
Instead of using a fixed timestep $h$ for all timesteps, in an adaptive algorithm the timestep size is adjusted at every time $t_n$. Here we use
\begin{equation}
  h^{(\text{adaptive})}_n \coloneqq \min\left\{h,\frac{\lambda(X_{\text{adapt}})}{\lambda(X_n)}h\right\},
\end{equation}
where $X_n$ is the current particle position and $X_{\text{adapt}}$ is a reference height. Since $\lambda(X)$ decreases with height, this choice guarantees that the product of the adaptive timestep size and $\lambda(X_n)$ never exceeds $h\,\lambda(X_{\text{adapt}})$, which can then be chosen such that the SE method is stable. Since the telescoping sum in Eq.~(\ref{eqn:telescoping_sum}) is preserved, this choice of adaptive timestep is independent of the MLMC level and therefore does not introduce any additional bias \cite{Giles2016}. Since $\lambda(X)$ is a decreasing function of $X$, it is only necessary to adapt the timestep size for particles below the reference height $X_{\text{adapt}}$. However, since $\lambda\rightarrow\infty$ close to the ground, it is still necessary to regularise $\lambda(X)$ as described in Section \ref{sec:regularisation}. Without regularisation, the timestep size can become arbitrarily small, and it can take an infinite amount of time to calculate even a single trajectory.

\begin{figure}
 \begin{center}
   \includegraphics[width=0.4\linewidth]{\figdir/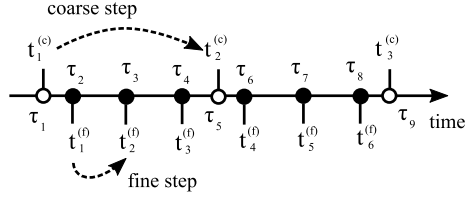}
   \caption{Non-nested timesteps in the adaptive MLMC algorithm.}
   \label{fig:AdaptiveTimesteps}
 \end{center}
\end{figure}

When using adaptivity within MLMC, the coarse and fine timesteps are not necessarily nested. However, it is easy to adapt the MLMC algorithm to account for this \cite{Giles2016}. For this, the time interval $[0,T]$ is divided into a number of intervals $[\tau_j,\tau_{j+1}]$ such that $\tau_0=0$, $\tau_N=T$ and each of the points $\tau_j$ is either a fine or a coarse time point (see Fig.\ \ref{fig:AdaptiveTimesteps}, which has been adapted from Fig.\ 1 in \cite{Giles2016}).
More explicitly, define
\begin{equation}
  \tau_j = \begin{cases}
    t_i^{(c)} & \text{for coarse time points,}\\
    t_{i'}^{(f)} & \text{for fine time points.}
  \end{cases}
\end{equation}
A normal random variable $\xi_j\sim \operatorname{Normal}(0,1)$ is associated with each interval $[\tau_j,\tau_{j+1}]$. Then, for each coarse or fine time interval $[t_i^{(x)},t_{i+1}^{(x)}]$ (where $x=c,f$), those random variables can be combined to construct a random increment $\Delta W(t^{(x)}_i,t^{(x)}_{i+1}) \sim \text{Normal}(0,t^{(x)}_{i+1}-t^{(x)}_{i})$ with the correct distribution. For SE with reflection, the expression for the coarse path is
\begin{equation}
 \Delta W(t^{(c)}_i,t^{(c)}_{i+1}) = S_i^{(c)}
\sum_{j=j_-}^{j_+-1} S_j^{(f)}
\xi_j \sqrt{\tau_{j+1}-\tau_{j}},\quad\text{where} \ \ \tau_{j_-}=t_i^{(c)} \ \text{and} \  \tau_{j_+}=t_{i+1}^{(c)}\,.\label{eqn:AdaptivityIncrementEuler}
\end{equation}
Here, the reflection factor is $S_j^{(\cdot)} = (-1)^{n_{\text{refl}}(\tau_j)}$, where, as before, $n_{\text{refl}}(\tau_j)$ counts the number of reflections up to time $\tau_j$. The corresponding expression $\Delta W(t^{(f)}_i,t^{(f)}_{i+1})$ for the fine path is the same as in (\ref{eqn:AdaptivityIncrementEuler}), with the only difference that the reflection factors $S_i^{(c)}$ and $S_j^{(f)}$ do not appear in this case. Since $\Delta W(t_i^{(x)},t_{i+1}^{(x)})$ is a sum of Gaussian random variables, it is easy to see that this variable has a variance of $t_{i+1}^{(x)}-t_i^{(x)}$ and can be used instead of $\sqrt{h}\,\xi_n$ for the velocity increment in (\ref{eqn:MethodSymplecticEuler}). The crucial observation is that if a fine interval overlaps with a coarse interval, they share some of the random increments and this ensures the correct coupling between fine and coarse paths. The expression in (\ref{eqn:AdaptivityIncrementEuler}) reduces to the standard expression for the case that the coarse points coincide with fine points.
For GL, the expression in (\ref{eqn:AdaptivityIncrementEuler}) has to be modified to
\begin{equation}
\begin{aligned}
 \Delta W(t^{(c)}_i,t^{(c)}_{i+1}) &= S_i^{(c)}\sum_{j=j_-}^{j_+-1}\Biggl\{
S_j^{(f)} \xi_j 
\sqrt{\frac{1-\exp[-2\lambda(X(\tau_j))\cdot (\tau_{j+1}-\tau_j)]}{2\lambda(X(\tau_j))}}\\
&\qquad\qquad \cdot \exp\Biggl[
-\sum_{k=j+1}^{j_+-1}\lambda(X(\tau_k))\cdot(\tau_{k+1}-\tau_k)
\Biggr] \Bigg\}.
\end{aligned}
\label{eqn:AdaptivityIncrementSplitting}
\end{equation}
This increment replaces $\sqrt{\frac{1 - \exp[- 2 \lambda(X_n) h ]}{2 \lambda(X_n)}}\xi_n$ in (\ref{eqn:MethodSymplecticEulerOU}). Again the corresponding expression for the fine path is obtained from Eq.\ (\ref{eqn:AdaptivityIncrementSplitting}) by removing the reflection factors $S_i^{(c)}$ and $S_j^{(f)}$. It can be verified that this is equivalent to (\ref{SymEulerCoupled}) in the case that the coarse timesteps coincide with fine timesteps.
\subsection{Numerical results}
\begin{figure}
\begin{subfigure}{.5\textwidth}
\centering
\captionsetup{justification=centering}
\includegraphics[scale = 0.4]{\figdir/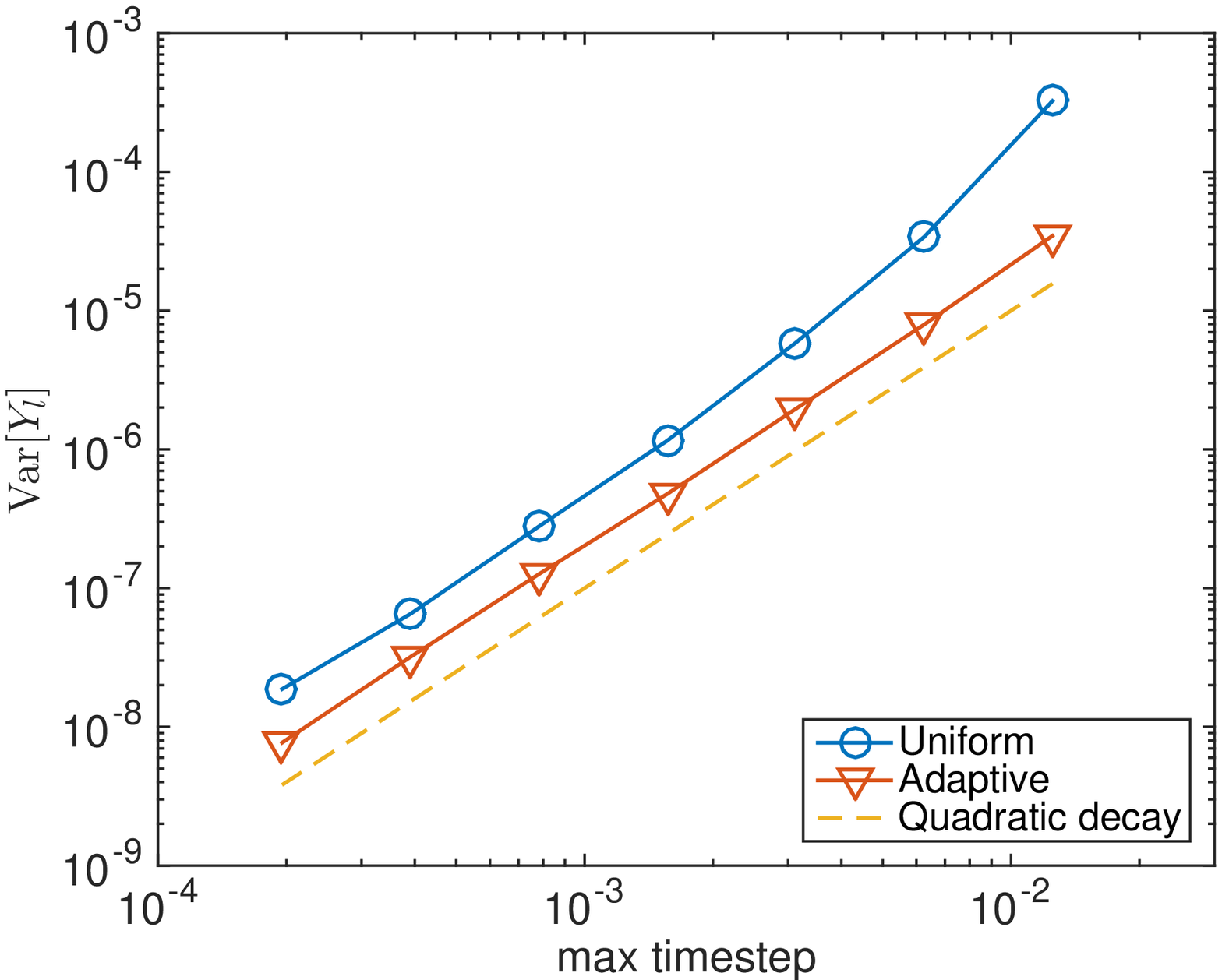}
\caption{Symplectic Euler}
\label{fig:variance_adaptivitySE}
\end{subfigure}%
\begin{subfigure}{.5\textwidth}
\centering
\captionsetup{justification=centering}
\includegraphics[scale = 0.4]{\figdir/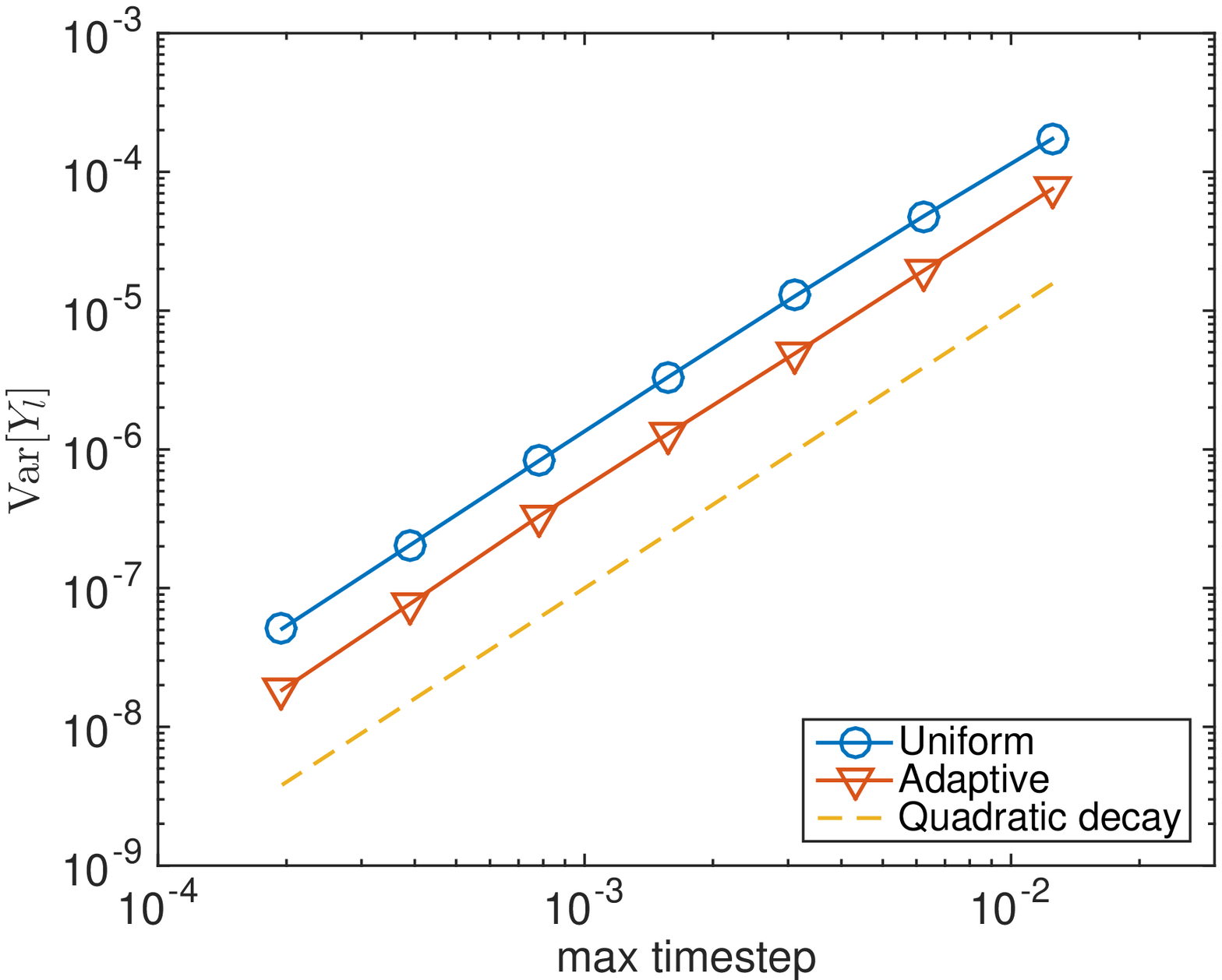}
\caption{Geometric Langevin}
\label{fig:variance_adaptivityGL}
\end{subfigure}%
\caption{Variance $\text{Var}[Y_\ell]$ as function of the timestep size}
\end{figure}

We now present some  numerical results to quantify potential improvements from adaptive timestepping for two integrators (SE and GL). The adaptive MLMC algorithm described in the previous section is used to calculate the expected value of the particle position $X_T$ at the final time with the same set-up as in Section \ref{sec:ParticlePosition} and $X_{\text{adapt}} = 0.05$. As before, we first confirm that condition (iii) in the complexity theorem is satisfied. The variance $\text{Var}[Y_\ell]$ as a function of the timestep size $h_\ell$ is shown for both SE and GL integrators in Figs.~\ref{fig:variance_adaptivitySE} and \ref{fig:variance_adaptivityGL}. As expected, both for the uniform and adaptive timestepping methods, the variance decays quadratically with the (maximal) timestep size, $\text{Var}[Y_\ell]\propto h_\ell^2$. In both cases, adaptive timestepping reduces the variance for a given $h_\ell$.

\begin{figure}
\begin{subfigure}{.5\textwidth}
\centering
\captionsetup{justification=centering}
\includegraphics[scale = 0.4]{\figdir/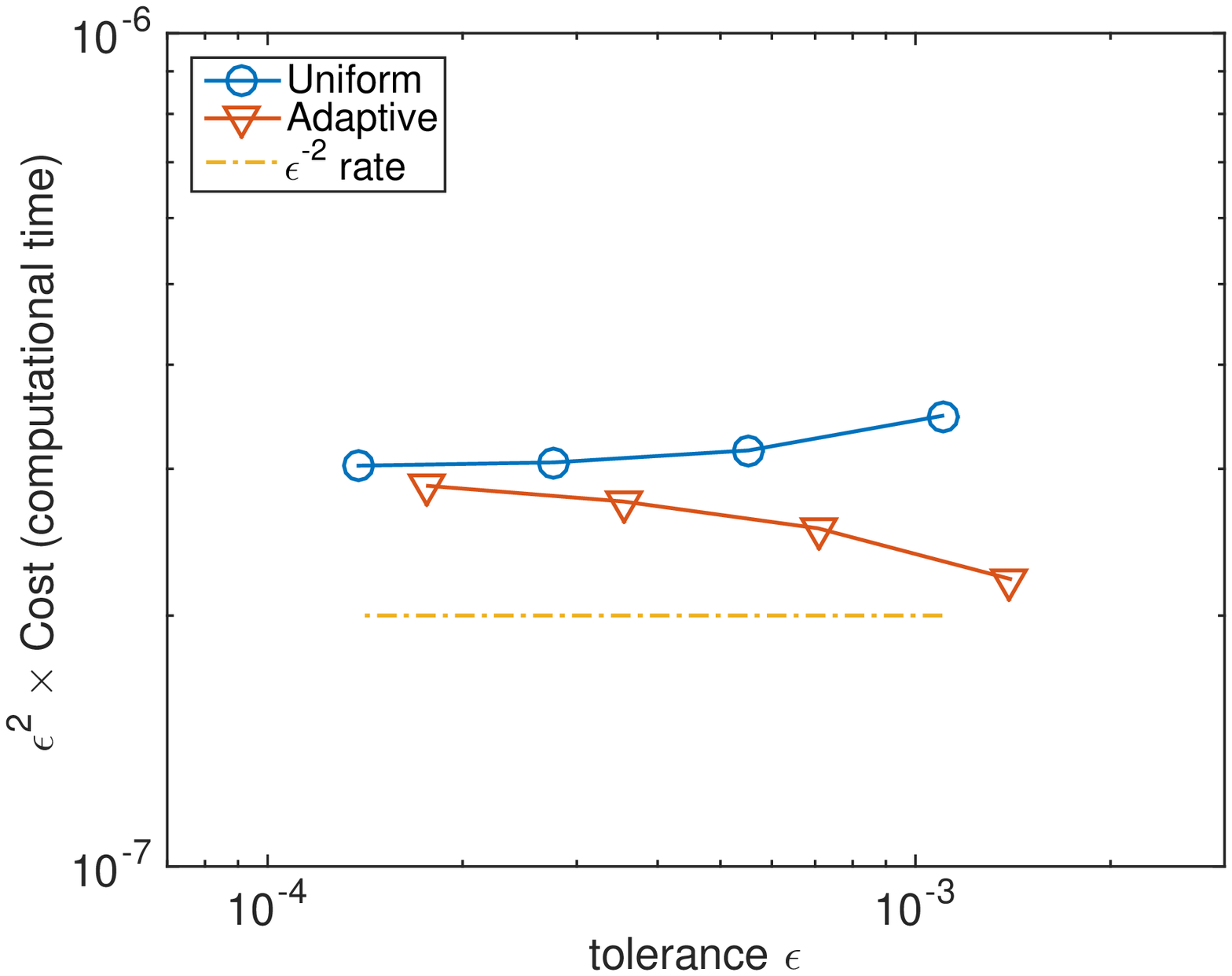}
\caption{Symplectic Euler}
\label{fig:cost_adaptivitySE}
\end{subfigure}%
\begin{subfigure}{.5\textwidth}
\centering
\captionsetup{justification=centering}
\includegraphics[scale = 0.4]{\figdir/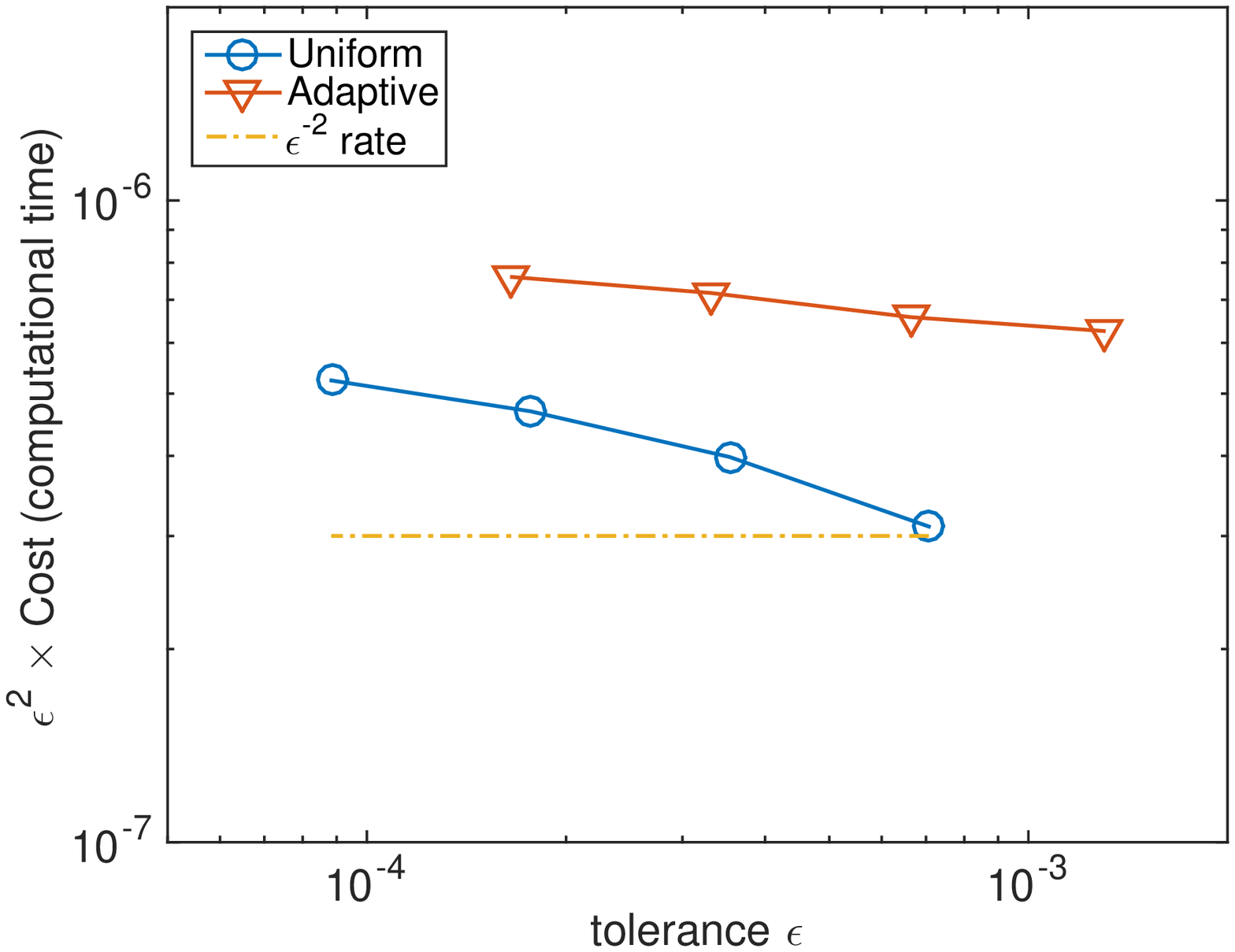}
\caption{Geometric Langevin}
\label{fig:cost_adaptivityGL}
\end{subfigure}%
\caption{Computational cost of adaptive- and non-adaptive timestepping algorithms for different tolerances $\epsilon$ on the root-mean-square error.}
\end{figure}

The total computational cost for both methods is compared in Figs.~\ref{fig:cost_adaptivitySE} and \ref{fig:cost_adaptivityGL}. While the results for SE are as expected and adaptive timestepping reduces the overall cost, the opposite can be observed for GL; here uniform timestepping reduces the runtime (albeit not by much as $\epsilon\rightarrow 0$). Looking at the bias error for a given (maximal) timestep size, we find that for the SE method this bias is about $6\times$ smaller with adaptive timestepping. This and the fact that the variance is reduced (see Fig.\ \ref{fig:variance_adaptivitySE}) explains the results in Fig.\ \ref{fig:cost_adaptivitySE}. For the GL method, on the other hand, adaptive timestepping \textit{increases} the bias by about $3\times$. This is a counter-intuitive result, since all individual timesteps in the adaptive algorithm are smaller than for the corresponding uniform implementation. For a (deterministic) ODE, it would then be easy to prove that adaptivity reduces the discretisation error. For the SDE considered here, however, the theoretical analysis is significantly harder and beyond the scope of this work. To exclude any bugs in the code, we have carried out several tests. For example, we have confirmed that for the corresponding ODE equations ($\sigma_U=0$) adaptivity reduces the bias error also for the GL method. The phenomenon appears to be linked to the divergence of $\lambda(X)$ close to the ground and might require careful tuning of both $X_{\text{adapt}}$ and $\epsilon_{\text{reg}}$ to see any benefits 
from adaptive timestepping.

\section{Conclusions and future work}\label{sec:Conclusion}
In this paper, we presented two improvements to Monte Carlo algorithms for the solution of SDEs encountered in atmospheric dispersion modelling: the geometric Langevin (GL) and BAOAB methods use a splitting approach to avoid numerical instabilities due to very small velocity decorrelation times close to the ground. Multilevel Monte Carlo reduces the algorithmic complexity as a function of the tolerance $\epsilon$ on the root-mean-square error from $\epsilon^{-3}$ to $\epsilon^{-2}$. Both methods were applied to the simulation of the vertical spread and diffusion of particles in a turbulent boundary layer. The intricacies of the application considered here requires careful adaptation of the algorithms to account for suitable boundary conditions at the top and bottom of the atmosphere and to treat the divergence of the inverse velocity decorrelation time near the ground. When predicting particle concentrations, a smoothed indicator function has to be used and we confirmed that the additional errors introduced by the smoothing are small.

The proposed splitting methods GL and BAOAB significantly reduce the computational cost as compared to Symplectic Euler, with gains of up to a factor 13 and 52, respectively, for all the quantities of interest considered in our numerical experiments. The main reason for this is a significant reduction in the bias error since the GL and BAOAB algorithms are able to deal better with large values of the inverse velocity decorrelation time $\lambda(X)$. As a modified equation analysis shows, the BAOAB integrator is almost second order and - in contrast to both SE and GL - the remaining $\mathcal{O}(h)$ term does not diverge near the lower boundary. For small tolerances $\epsilon$, the MLMC algorithm leads to further speed-ups, which are particularly pronounced if the particles are released at a larger height. When calculating piecewise constant approximations to the particle concentration field, the MLMC algorithm is nearly three times faster than the standard Monte Carlo algorithm even when used in combination with the improved GL and BAOAB timestepping methods.

In a preliminary study we also explored the use of adaptive timestepping and find that this can lead to a further, but less dramatic, reduction in runtime for the symplectic Euler method. The unexpected results observed for the GL algorithm require further investigation and careful tuning of parameters.

Of central importance for the application of the MLMC algorithm is the correct dependence of the multilevel variance on the timestep size. Numerically we observed the correct behaviour, but a theoretical proof is more challenging due to the non-smooth regularisation of the profile functions. Although this problem can be addressed by using an equivalent smooth regularisation instead, a proof relies on further technical details such as boundedness of all moments of $U_t$, which would then allow the application of the modified equation analysis in \cite{EikeRobTony}. The treatment of reflective boundary conditions introduces discontinuities, which do not have any adverse effects on our numerical results but require a careful treatment in a rigorous mathematical analysis. Since the focus of this work is the introduction of new numerical methods to atmospheric dispersion modelling, this would be beyond the scope of this article.

In addition to investigating those mathematical details more rigorously, there are various ways of extending our work. One obvious avenue to explore is the use of quasi Monte Carlo methods (QMCs), which potentially allows a reduction of the computational complexity beyond the standard Monte Carlo limit of $\epsilon^{-2}$, set by the central limit theorem (see e.g. \cite{Kuo2015}). By using a randomly shifted point lattice, any additional bias from the QMC sampling can be removed.

So far we have treated turbulence in all space directions independently and focussed on motion in the vertical direction where the atmospheric conditions vary strongly. Although similar approximations are made in operational models, ultimately we would like to apply our approach to a full three-dimensional scenario, for example to study the global spread of ash clouds, which would require the use of a more complicated, 3d turbulence model described in \cite{Thomson1987}. One particular question to be addressed is as to whether the separation of particle pairs by the resolved flow field will have any impact on the efficiency of MLMC.

 NAME models a plethora of additional physical effects such as gravitational settling, chemical interactions, radioactive decay chains and gamma radiation from non-local radioactivity, which are computationally expensive. Incorporating more of the physics into our methodology is the subject of future work. For example, radioactive particles require functionals that depend on the full path to establish the total exposure to radiation, rather than the final-time functionals in this paper. In another direction,  this paper describes the multilevel approach  in terms of a hierarchy of levels of increasing timestep size and the first coarsening step might consist of using the same timestep size but not including all physical processes. Conceptually this is somewhat similar to the difference between $p$- and $h$- refinement in finite element methods. The multilevel approach is very flexible and can also be adapted to computing the probability of rare events by use of the Girsanov transform and importance sampling \cite{Kebaier2015-ie}. This could for example be used to calculate the particle concentration in boxes at the edge of the plume.
It is our hope that eventually operational models such as NAME will make use of the Multilevel Monte Carlo method and the improved GL- and BAOAB- integrators described in this work.
\section*{Acknowledgements}
The PhD of Grigoris Katsiolides is supported by an EPSRC CASE studentship, which is partially funded by the Met Office.

\appendix
\section{Modified equations and regularisation}\label{sec:ModifiedEquations}

Our theoretical tools for proving the complexity theorem rely on the theory of modified equations \cite{shardlow2006modified,zygalakis2011existence} and their application to the MLMC method as described in \cite{EikeRobTony}; this approach is an alternative to the usual strong convergence proofs. We discuss this analysis here since it can be used to justify the use of a regularisation cut-off at the upper and lower boundary as outlined in Section \ref{sec:BoundaryConditions}.

The key idea is as follows: under suitable conditions, the numerical solution $X_n$ of the initial SDE (\ref{generalisedsim}) is the weak second-order approximation of the solution of a set of modified equations -- recall that $X_n$ converges only with order one to the original equations. If those modified equations exist, it can be shown that condition (iii) in the complexity theorem holds \cite{EikeRobTony}.

For our model, the modified equations have the form
\begin{equation}
\begin{aligned}
dU &= \left( - \lambda(X) U - \frac{\partial V}{\partial X}(X,U) + h v_1^{(U)}(X,U) \right) dt + \left( \sigma(X) + h \Sigma_1^{(U)} (X,U) \right) dW_t, \\
dX &= \left( U + h v_1^{(X)}(X,U) \right) dt + h \Sigma_1^{(X)} (X,U) dW_t.
\end{aligned}
\label{eqn:ModifiedEquations}
\end{equation}
The functions $v_1^{(U)}$, $v_1^{(X)}$, $\Sigma_1^{(U)}$ and $\Sigma_1^{(X)}$ depend on the numerical timestepping method and can be expressed explicitly in terms of $\tau(X)$ and $\sigma_U(X)$. However, these modified equations are only a first-order approximation of Eq.\ (\ref{generalisedsim}) if the $\mathcal{O}(h)$ terms in \eqref{eqn:ModifiedEquations} are indeed small. Unfortunately, this is not necessarily the case near the boundary of the domain. In particular, by expanding
\begin{xalignat*}{2}
  v_1^{(U)}(X,U) &= X^{-1}\cdot f^{(v,U)}(U) + \dots,&
  v_1^{(X)}(X,U) &= X^{-2}\cdot f^{(v,X)}(U) + \dots,\\
  \Sigma_1^{(U)}(X,U) &= X^{-3/2}\cdot f^{(\Sigma,U)}(U) + \dots,&
  \Sigma_1^{(X)}(X,U) &= X^{-1/2}\cdot f^{(\Sigma,X)}(U) + \dots
\end{xalignat*}
it can be shown that these terms diverge at the lower boundary ($X\rightarrow 0$), for both the symplectic Euler and geometric Langevin method. The functions $f^{(\cdot,\cdot)}$ are low-order polynomials in $U$ which can be written as
\begin{xalignat*}{2}
  f^{(v,U)}(U) &= a_1^{(v,U)}U, &
  f^{(v,X)}(U) &= a_1^{(v,X)}U + a_2^{(v,X)}U^2,\\
  f^{(\Sigma,U)}(U) &= a_0^{(\Sigma,U)}+a_1^{(\Sigma,U)}U, &
  f^{(\Sigma,X)}(U) &= a_0^{(\Sigma,X)}.
\end{xalignat*}
The constants $a_j^{(\cdot,\cdot)}$ depend on the timestepping method.
 Similar expressions can be derived at the upper boundary ($X\rightarrow H$), but for simplicity we do not write them down here. 

However, if we regularise the profiles $\tau(X)$ and $\sigma_U(X)$ as described in Section \ref{sec:regularisation}, then the divergent terms $X^{-\alpha}$ are replaced by the constant $\epsilon_{\text{reg}}^{-\alpha}$. Thus, the terms $v_1^{(U)}$, $v_1^{(X)}$, $\Sigma_1^{(U)}$ and $\Sigma_1^{(X)}$ can all be made $\mathcal{O}(1)$ by a suitable choice of $\epsilon_{\text{reg}}$. The existence of the modified equations then implies the quadratic variance decay required by condition (iii) in the complexity theorem, under the condition that $U$ has bounded moments. This can be seen in Fig.\ \ref{SymEulerOUReg}: unless $\epsilon_{\text{reg}}$ is chosen too small, the variance decay is indeed quadratic. 

For the BAOAB method, on the other hand, the $\mathcal{O}(h)$ terms in Eq. (\ref{eqn:ModifiedEquations}) are significantly simpler and can be written down exactly as 
\begin{xalignat*}{2}
  v_1^{(U)}(X,U) &= -\frac{1}{8} \frac{\partial}{\partial U}\left(\frac{\partial V}{\partial X}\right)^2, &
  v_1^{(X)}(X,U) &= 0,\\
  \Sigma_1^{(U)}(X,U) &= 0, &
  \Sigma_1^{(X)}(X,U) &= 0.
\end{xalignat*}
The only non-vanishing term $v_1^{(U)}$ does not diverge at the lower boundary. As $X\rightarrow 0$ it approaches the limit $v_1^{(U)}(X=0,U) = b_1^{(v,U)} U + b_3^{(v,U)}U^3$ with constant polynomial coefficients $b_j^{(v,U)}$. This explains the superior behaviour of the BAOAB method for which the variance decay is quadratic, irrespective of the regularisation height $\epsilon_{\text{reg}}$.
\section{Derivation of error bounds for smoothed indicator function} \label{weakApp}
To derive inequality \eqref{weak}, we first consider the following difference between the expected values of the indicator function $P$ and the smoothed out function $P^{r,\delta}$:
\begin{equation}
\begin{aligned}
\left| \mathbb{E}[P] - \mathbb{E}[P^{r,\delta}] \right| 
&\leq \textstyle \left| \mathbb{E} \left[ \theta (X_T-b) \right] - \mathbb{E} \left[ g_r \left( \frac{X_T - b}{\delta} \right) \right] \right| + \left| \mathbb{E} \left[ \theta (X_T-a) \right] - \mathbb{E} \left[ g_r \left( \frac{X_T - a}{\delta} \right) \right] \right| 
\leq  2 c \delta^{r + 1}.\label{eqn:rdeltaBound1}
\end{aligned}
\end{equation}
Using this expression and $P-\hat{P}^{r,\delta}_L=(P - P^{r,\delta}) + (P^{r,\delta} - \hat{P}^{r,\delta}_L)$, the bound on the bias (or weak error) in \eqref{weak} follows from the triangle inequality
\begin{align*}
\left| \mathbb{E}[P] - \mathbb{E}[\hat{P}^{r,\delta}_L] \right| &\leq \left| \mathbb{E}[P] - \mathbb{E}[P^{r,\delta}] \right| + \left| \mathbb{E}[P^{r,\delta}] - \mathbb{E}[\hat{P}^{r,\delta}_L] \right| \leq 2 c \delta^{r + 1} + \left| \mathbb{E}[P^{r,\delta}] - \mathbb{E}[P^{r,\delta}_L] \right|.
\end{align*}
To prove \eqref{strong}, we write 
\begin{equation}
\hat{P}^{r,\delta}_L-\mathbb{E}[P] = 
\left(
  \hat{P}^{r,\delta}_L - \mathbb{E}[P^{r,\delta}_L]
\right)
+
\left(
  \mathbb{E}[P^{r,\delta}_L] - \mathbb{E}[P^{r,\delta}]
\right)
+
\left(
  \mathbb{E}[P^{r,\delta}] - \mathbb{E}[P]
\right)\,.
\end{equation}
Now, using the fact that $\mathbb{E}[\hat{P}^{r,\delta}_L]=\mathbb{E}[P^{r,\delta}_L]$
it is easy to see that
\begin{equation}
\begin{aligned}
\mathbb{E} \left[ \left( \hat{P}^{r,\delta}_L - \mathbb{E}[P] \right)^2 \right] &= \text{Var}( \hat{P}^{r,\delta}_L) + \left( \mathbb{E}[P^{r,\delta}_L] - \mathbb{E}[P^{r,\delta}] \right)^2 + \left(\mathbb{E}[P^{r,\delta}]-\mathbb{E}[P]\right)^2 \\
& + 2\left(\mathbb{E}[P^{r,\delta}_L]-\mathbb{E}[P^{r,\delta}]\right)\left(\mathbb{E}[P^{r,\delta}]-\mathbb{E}[P]]\right).
\end{aligned}
\end{equation}
The desired result then follows from the bound in \eqref{eqn:rdeltaBound1}:
\begin{equation}
\mathbb{E} \left[ \left( \hat{P}^{r,\delta}_L - \mathbb{E}[P] \right)^2 \right] \leq \text{Var}( \hat{P}^{r,\delta}_L) + \left( \mathbb{E}[P^{r,\delta}_L] -\mathbb{E}[P^{r,\delta}] \right)^2 + 4 c^2 \delta^{2(r + 1)} + 4 c \delta^{r + 1} \left| \mathbb{E}[P^{r,\delta}_L] -\mathbb{E}[P^{r,\delta}] \right|
\end{equation}
\section{Data for the number of samples on each MLMC level} \label{DataOnLevels}
Tab. \ref{tab:DataOnNumberOfSamples} contains further details on the parameters and number of samples per level for two representative MLMC runs. Values are shown both for the mean particle position and the smoothed indicator function for all integrators (SE, GL and BAOAB). For a given tolerance $\epsilon$ on the root mean square error the table lists the number of coarse level time steps $M_0$, the number of levels $L$ and the time step size $h_L$ on the finest level with $h_L^{-1} = M_0 2^{L}$. The remaining rows in the table list the number of samples $N_\ell$ on level $\ell$.
\begin{table}[H]
  \centering
  \begin{tabular}{l|rrr|rrr}
    & \multicolumn{3}{|c|}{Mean particle position} & \multicolumn{3}{|c}{Concentration using $P^{r=4, \delta=0.1}_{[a,b]}(X_T)$} \\
    & SE & GL & BAOAB & SE & GL & BAOAB \\
    \hline\hline
    $\epsilon$ & $6.9\cdot 10^{-5}$ & $8.8\cdot 10^{-5}$ & $8.8\cdot 10^{-5}$ & $6.9\cdot 10^{-5}$ & $8.8\cdot 10^{-5}$ & $8.8\cdot 10^{-5}$ \\
    \hline
    $M_0$ & 40 & 40 & 40 & 80 & 80 & 40 \\
    \hline
    $L$ & 8 & 4 & 2 & 7 & 3 & 2 \\
    \hline
    $h_L$ & $9.8\cdot 10^{-5}$ & $1.6\cdot 10^{-3}$ & $6.3\cdot 10^{-3}$ & $9.8\cdot 10^{-5}$ & $1.6\cdot 10^{-3}$ & $6.3\cdot 10^{-3}$ \\
    \hline
    $N_0$ & 6787563 & 3899945 & 3432656 & 68039360 & 37703814 & 43656860 \\
    $N_1$ & 843711 & 357389 & 305644 & 10473620 & 5864870 & 11875896 \\
    $N_2$ & 190585 & 133744 & 113927 & 3220810 & 2197495 & 4735493 \\
    $N_3$ & 55965 & 48913 &  & 1055256 & 788799 & \\
    $N_4$ & 18142 & 17103 &  & 361439 &  & \\
    $N_5$ & 6064 &  &  & 126656 &  & \\
    $N_6$ & 2016 &  &  & 44076 &  & \\
    $N_7$ & 748 &  &  & 15776 &  & \\
    $N_8$ & 266 &  &  &  &  & \\
  \end{tabular}
  \caption{MLMC parameters and number of samples per level for the mean particle position and concentration.}
  \label{tab:DataOnNumberOfSamples}
\end{table}

\ifbool{PREPRINT}{ 
\bibliographystyle{elsarticle-num}

}{
\section*{References}
\bibliographystyle{elsarticle-num}
\bibliography{met_office_paper}

\begin{thebibliography}{10}
\expandafter\ifx\csname url\endcsname\relax
  \def\url#1{\texttt{#1}}\fi
\expandafter\ifx\csname urlprefix\endcsname\relax\def\urlprefix{URL }\fi
\expandafter\ifx\csname href\endcsname\relax
  \def\href#1#2{#2} \def\path#1{#1}\fi

\bibitem{Draxler2015}
R.~Draxler, D.~Arnold, M.~Chino, S.~Galmarini, M.~Hort, A.~Jones,
  S.~Leadbetter, A.~Malo, C.~Maurer, G.~Rolph, K.~Saito, R.~Servranckx,
  T.~Shimbori, E.~Solazzo, G.~Wotawa, {World Meteorological Organization's
  model simulations of the radionuclide dispersion and deposition from the
  Fukushima Daiichi nuclear power plant accident}, Journal of Environmental
  Radioactivity 139 (2015) 172 -- 184.
\newblock \href {http://dx.doi.org/10.1016/j.jenvrad.2013.09.014}
  {\path{doi:10.1016/j.jenvrad.2013.09.014}}.

\bibitem{Dacre2011}
H.~F. Dacre, A.~L. Grant, R.~J. Hogan, S.~E. Belcher, D.~Thomson, B.~Devenish,
  F.~Marenco, M.~Hort, J.~M. Haywood, A.~Ansmann, et~al., {Evaluating the
  structure and magnitude of the ash plume during the initial phase of the 2010
  Eyjafjallaj{\"o}kull eruption using lidar observations and NAME simulations},
  Journal of Geophysical Research: Atmospheres 116~(D20).
\newblock \href {http://dx.doi.org/10.1029/2011JD015608}
  {\path{doi:10.1029/2011JD015608}}.

\bibitem{Webster2012}
H.~N. Webster, D.~J. Thomson, B.~T. Johnson, I.~P.~C. Heard, K.~Turnbull,
  F.~Marenco, N.~I. Kristiansen, J.~Dorsey, A.~Minikin, B.~Weinzierl,
  U.~Schumann, R.~S.~J. Sparks, S.~C. Loughlin, M.~C. Hort, S.~J. Leadbetter,
  B.~J. Devenish, A.~J. Manning, C.~S. Witham, J.~M. Haywood, B.~W. Golding,
  {Operational prediction of ash concentrations in the distal volcanic cloud
  from the 2010 Eyjafjallaj{\"o}kull eruption}, J. Geophys. Res. 117~(D20).
\newblock \href {http://dx.doi.org/10.1029/2011JD016790}
  {\path{doi:10.1029/2011JD016790}}.

\bibitem{Jones2004}
A.~{Jones}, {Atmospheric dispersion modelling at the Met Office}, Weather 59
  (2004) 311--316.
\newblock \href {http://dx.doi.org/10.1256/wea.106.04}
  {\path{doi:10.1256/wea.106.04}}.

\bibitem{Jones2007}
A.~Jones, D.~Thomson, M.~Hort, B.~Devenish, {The U.K. Met Office's
  Next-Generation Atmospheric Dispersion Model, NAME III}, Air Pollution
  Modeling and Its Application XVII (2007) 580--589\href
  {http://dx.doi.org/10.1007/978-0-387-68854-1_62}
  {\path{doi:10.1007/978-0-387-68854-1_62}}.

\bibitem{Smith1989}
F.~B. Smith, M.~J. Clark, {The transport and deposition of airborne debris from
  the Chernobyl nuclear power plant accident with special emphasis on the
  consequences to the United Kingdom}, Meteorological Office Scientific Paper
  No. 42.

\bibitem{Redington2009}
A.~L. Redington, R.~G. Derwent, C.~S. Witham, A.~J. Manning, {Sensitivity of
  modelled sulphate and nitrate aerosol to cloud, pH and ammonia emissions},
  {Atmospheric Environment} {43}~({20}) ({2009}) {3227--3234}.
\newblock \href {http://dx.doi.org/10.1016/j.atmosenv.2009.03.041}
  {\path{doi:10.1016/j.atmosenv.2009.03.041}}.

\bibitem{Heinrich2001}
S.~Heinrich, {Multilevel Monte Carlo Methods}, in: S.~Margenov,
  J.~Wa\'{s}niewski, P.~Yalamov (Eds.), {Large-Scale Scientific Computing},
  Vol. 2179 of Lecture Notes in Computer Science, Springer Berlin Heidelberg,
  2001, pp. 58--67.
\newblock \href {http://dx.doi.org/10.1007/3-540-45346-6_5}
  {\path{doi:10.1007/3-540-45346-6_5}}.

\bibitem{Giles2008}
M.~B. Giles, {Multilevel Monte Carlo Path Simulation}, Operations Research
  56~(3) (2008) 607--617.
\newblock \href {http://dx.doi.org/10.1287/opre.1070.0496}
  {\path{doi:10.1287/opre.1070.0496}}.

\bibitem{Giles2015}
M.~B. Giles, {Multilevel Monte Carlo methods}, Acta Numerica 24 (2015) 259.
\newblock \href {http://dx.doi.org/10.1017/S096249291500001X}
  {\path{doi:10.1017/S096249291500001X}}.

\bibitem{giles2012multilevel}
M.~Giles, L.~Szpruch, \href{https://arxiv.org/abs/1212.1377}{{Multilevel Monte
  Carlo methods for applications in finance}}, arXiv preprint arXiv:1212.1377.
\newline\urlprefix\url{https://arxiv.org/abs/1212.1377}

\bibitem{graham2016mixed}
I.~G. Graham, R.~Scheichl, E.~Ullmann, {Mixed finite element analysis of
  lognormal diffusion and multilevel Monte Carlo methods}, Stochastics and
  Partial Differential Equations Analysis and Computations 4~(1) (2016) 41--75.
\newblock \href {http://dx.doi.org/10.1007/s40072-015-0051-0}
  {\path{doi:10.1007/s40072-015-0051-0}}.

\bibitem{EikeRobTony}
E.~H. M{\"u}ller, R.~Scheichl, T.~Shardlow, {Improving multilevel Monte Carlo
  for stochastic differential equations with application to the Langevin
  equation}, Proceedings of the Royal Society A 471~(2176).
\newblock \href {http://dx.doi.org/10.1098/rspa.2014.0679}
  {\path{doi:10.1098/rspa.2014.0679}}.

\bibitem{Skorokhod1961}
A.~V. Skorokhod, {Stochastic Equations for Diffusion Processes in a Bounded
  Region}, Theory of Probability \& Its Applications 6~(3) (1961) 264--274.
\newblock \href {http://arxiv.org/abs/http://dx.doi.org/10.1137/1106035}
  {\path{arXiv:http://dx.doi.org/10.1137/1106035}}, \href
  {http://dx.doi.org/10.1137/1106035} {\path{doi:10.1137/1106035}}.

\bibitem{bernardo2008piecewise}
M.~Bernardo, C.~Budd, A.~R. Champneys, P.~Kowalczyk, {Piecewise-smooth
  dynamical systems: theory and applications}, Vol. 163, Springer Science \&
  Business Media, 2008.

\bibitem{Wilson1993}
J.~D. Wilson, T.~K. Flesch, {Flow boundaries in random-flight dispersion
  models: enforcing the well-mixed condition}, Journal of Applied Meteorology
  32~(11) (1993) 1695--1707.

\bibitem{Rodean}
H.~C. Rodean, {Stochastic Lagrangian Models of Turbulent Diffusion}, Vol.~26,
  American Meteorological Society, 1996.

\bibitem{GilesNagapetyanRitter}
M.~B. Giles, T.~Nagapetyan, K.~Ritter, {Multilevel Monte Carlo Approximation of
  Distribution Functions and Densities}, SIAM/ASA Journal on Uncertainty
  Quantification 3~(1) (2015) 267--295.
\newblock \href {http://dx.doi.org/10.1137/140960086}
  {\path{doi:10.1137/140960086}}.

\bibitem{Kloeden2011}
P.~Kloeden, E.~Platen,
  \href{https://books.google.co.uk/books?id=BCvtssom1CMC}{{Numerical Solution
  of Stochastic Differential Equations}}, Stochastic Modelling and Applied
  Probability, Springer Berlin Heidelberg, 2011.
\newline\urlprefix\url{https://books.google.co.uk/books?id=BCvtssom1CMC}

\bibitem{Leimkuhler2015}
B.~Leimkuhler, C.~Matthews,
  \href{https://books.google.co.uk/books?id=-5oxrgEACAAJ}{{Molecular Dynamics:
  With Deterministic and Stochastic Numerical Methods}}, Interdisciplinary
  Applied Mathematics, Springer International Publishing, 2015.
\newline\urlprefix\url{https://books.google.co.uk/books?id=-5oxrgEACAAJ}

\bibitem{Maruyama1955}
G.~Maruyama, {Continuous Markov processes and stochastic equations}, Rendiconti
  del Circolo Matematico di Palermo 4~(1) (1955) 48--90.

\bibitem{bou2010long}
N.~Bou-Rabee, H.~Owhadi, {Long-run accuracy of variational integrators in the
  stochastic context}, SIAM Journal on Numerical Analysis 48~(1) (2010)
  278--297.
\newblock \href {http://dx.doi.org/10.1137/090758842}
  {\path{doi:10.1137/090758842}}.

\bibitem{Hoel2011}
H.~Hoel, E.~von Schwerin, A.~Szepessy, R.~Tempone, {Adaptive Multilevel Monte
  Carlo Simulation}, in: Numerical Analysis of Multiscale Computations:
  Proceedings of a Winter Workshop at the Banff International Research Station
  2009, Vol.~82, Springer Science \& Business Media, 2011, p. 217.

\bibitem{Hoel2014}
H.~Hoel, E.~Von~Schwerin, A.~Szepessy, R.~Tempone, {Implementation and analysis
  of an adaptive multilevel Monte Carlo algorithm}, Monte Carlo Methods and
  Applications 20~(1) (2014) 1--41.
\newblock \href {http://dx.doi.org/10.1515/mcma-2013-0014}
  {\path{doi:10.1515/mcma-2013-0014}}.

\bibitem{Giles2016}
M.~B. Giles, C.~Lester, J.~Whittle, Non-nested Adaptive Timesteps in Multilevel
  Monte Carlo Computations, Springer International Publishing, Cham, 2016, pp.
  303--314.
\newblock \href {http://dx.doi.org/10.1007/978-3-319-33507-0_14}
  {\path{doi:10.1007/978-3-319-33507-0_14}}.

\bibitem{MLMCLangevinCode2017}
E.~H. Mueller, G.~Katsiolides, T.~Shardlow, {MLMCLangevin code} (May 2017).
\newblock \href {http://dx.doi.org/10.5281/zenodo.580662}
  {\path{doi:10.5281/zenodo.580662}}.

\bibitem{shardlow2006modified}
T.~Shardlow, {Modified equations for stochastic differential equations}, BIT
  Numerical Mathematics 46~(1) (2006) 111--125.
\newblock \href {http://dx.doi.org/10.1007/s10543-005-0041-0}
  {\path{doi:10.1007/s10543-005-0041-0}}.

\bibitem{zygalakis2011existence}
K.~C. Zygalakis, {On the existence and the applications of modified equations
  for stochastic differential equations}, SIAM Journal on Scientific Computing
  33~(1) (2011) 102--130.
\newblock \href {http://dx.doi.org/10.1137/090762336}
  {\path{doi:10.1137/090762336}}.

\bibitem{Thomson1987}
D.~J. Thomson, {Criteria for the selection of stochastic models of particle
  trajectories in turbulent flows}, Journal of Fluid Mechanics 180 (1987)
  529--556.
\newblock \href {http://dx.doi.org/10.1017/S0022112087001940}
  {\path{doi:10.1017/S0022112087001940}}.

\bibitem{Kolmogorov1941}
A.~N. Kolmogorov, \href{http://www.jstor.org/stable/51980}{{The local structure
  of turbulence in incompressible viscous fluid for very large Reynolds
  numbers}}, in: Dokl. Akad. Nauk SSSR, Vol.~30, JSTOR, 1941, pp. 301--305.
\newline\urlprefix\url{http://www.jstor.org/stable/51980}

\bibitem{Wilson2009}
J.~Wilson, E.~Yee, N.~Ek, R.~d'Amours, {Lagrangian simulation of wind transport
  in the urban environment}, Quarterly Journal of the Royal Meteorological
  Society 135~(643) (2009) 1586--1602.
\newblock \href {http://dx.doi.org/10.1002/qj.452} {\path{doi:10.1002/qj.452}}.

\bibitem{Webster2003}
H.~Webster, D.~Thomson, N.~Morrison, {New turbulence profiles for NAME}, Met
  Office Turbulence and Diffusion note 288.

\bibitem{khasminskii2011stochastic}
R.~Khasminskii, Stochastic stability of differential equations, Vol.~66,
  Springer Science \& Business Media, 2011.

\bibitem{Lord2014}
G.~J. Lord, C.~E. Powell, T.~Shardlow, {An introduction to computational
  stochastic PDEs}, no.~50, Cambridge University Press, New York, 2014.

\bibitem{Collier2015}
N.~Collier, A.-L. Haji-Ali, F.~Nobile, E.~von Schwerin, R.~Tempone, A
  continuation multilevel monte carlo algorithm, BIT Numerical Mathematics
  55~(2) (2015) 399--432.

\bibitem{Hoel2012}
H.~Hoel, E.~Schwerin, A.~Szepessy, R.~Tempone, Numerical Analysis of Multiscale
  Computations: Proceedings of a Winter Workshop at the Banff International
  Research Station 2009, Springer Berlin Heidelberg, 2012, Ch. Adaptive
  Multilevel Monte Carlo Simulation, pp. 217--234.
\newblock \href {http://dx.doi.org/10.1007/978-3-642-21943-6_10}
  {\path{doi:10.1007/978-3-642-21943-6_10}}.

\bibitem{hoel2014implementation}
H.~Hoel, E.~Von~Schwerin, A.~Szepessy, R.~Tempone, {Implementation and analysis
  of an adaptive multilevel Monte Carlo algorithm}, Monte Carlo Methods and
  Applications 20~(1) (2014) 1--41.
\newblock \href {http://dx.doi.org/10.1515/mcma-2013-0014}
  {\path{doi:10.1515/mcma-2013-0014}}.

\bibitem{Fang2016}
W.~Fang, M.~B. Giles, {Adaptive Euler-Maruyama Method for SDEs with
  Non-globally Lipschitz Drift: Part I, Finite Time Interval}, arXiv preprint
  arXiv:1609.08101.

\bibitem{Kuo2015}
F.~Y. Kuo, R.~Scheichl, C.~Schwab, I.~H. Sloan, E.~Ullmann,
  \href{https://arxiv.org/abs/1507.01090}{{Multilevel quasi-Monte Carlo methods
  for lognormal diffusion problems}}, arXiv preprint arXiv:1507.01090.
\newline\urlprefix\url{https://arxiv.org/abs/1507.01090}

\bibitem{Kebaier2015-ie}
A.~Kebaier, J.~Lelong, Coupling importance sampling and multilevel {M}onte
  {C}arlo using sample average approximation (13~Oct. 2015).
\newblock \href {http://arxiv.org/abs/1510.03590} {\path{arXiv:1510.03590}}.

\end{thebibliography}
}
\end{document}